\documentclass{amsart}

\usepackage[all]{xy}
\usepackage{amssymb, amsmath, amsfonts, amsthm, amsbsy, amscd}
         
\newtheorem{theorem}{Theorem}[section]

\newtheorem*{theoremb}{Theorem B}

\newtheorem*{theoremai}{Theorem A}

\newtheorem*{theoremci}{Theorem C}
\newtheorem{claim}{Claim}{}
\newtheorem{proposition}{Proposition}[section]
\newtheorem{corollary}{Corollary}[section]
\newtheorem{lemma}{Lemma}[section]

\newtheorem{definition}{Definition}[section]
{}
\newtheorem{conditionb}{}{}
\newtheorem{condition}{}{}
\newtheorem{conditionai}{}{}
{}
{}

\theoremstyle{definition}
\newtheorem{example}{Example}[section]
\newtheorem{remark}{Remark}[section]

\numberwithin{equation}{section}

\newcommand{\conref}[1]{\textbf{\ref{#1}}}

\newcommand{\F}{\mathcal{F}}

\newcommand{\Reta}{R^{\eta}}
\newcommand{\taueta}{\tau^{\eta}}
\newcommand{\svol}{\Psi}

\newcommand{\sOmega}{\bar{\Omega}}
\newcommand{\C}{\mathcal{C}}
\newcommand{\At}{\mathcal{A}}
\newcommand{\B}{\mathcal{B}}

\newcommand{\ppi}{\bar{\pi}}
\newcommand{\pframe}{\bar{\mathcal{G}}}

\newcommand{\standrepz}{\mathbb{V}^{\times}}

\newcommand{\A}{O}

\newcommand{\contactop}{\mathcal{CO}^{+}}
\newcommand{\pcurv}{\mathcal{K}}

\newcommand{\etanabla}{\nabla^{\eta}}
\newcommand{\assrep}{\mathcal{E}}
\newcommand{\rep}{\mathbb{E}}
\newcommand{\standomega}{\bar{\Omega}}

\newcommand{\levi}{\mathcal{L}}

\newcommand{\G}{\mathcal{G}}

\newcommand{\spn}{span}
\newcommand{\graded}{Gr}

\newcommand{\emf}{\mathcal{E}}
\newcommand{\tcurv}{\mathcal{R}}
\newcommand{\adaptedframe}{\mathcal{G}}

\newcommand{\tR}{\mathcal{R}}

\newcommand{\tnabla}{\boldsymbol{\nabla}}
\newcommand{\tractor}{\mathcal{T}}
\newcommand{\standrep}{\mathbb{V}}

\newcommand{\contactb}{\mathsf{C}}

\newcommand{\form}{\mathsf{L}}
\newcommand{\dens}{\mathsf{L}}

\newcommand{\eul}{\mathbb{X}}

\newcommand{\candens}{\mathsf{F}}

\newcommand{\rb}{T}
\newcommand{\proj}{\mathbb{P}}

\newcommand{\symnabla}{\,^{S}\nabla}
\newcommand{\sym}{\text{Sym}}

\newcommand{\lie}{\mathfrak{L}}

\newcommand{\vect}{\text{Vec}}

\newcommand{\g}{\mathfrak{g}}
\newcommand{\f}{\mathfrak{f}}

\newcommand{\p}{\mathfrak{p}}

\newcommand{\ad}{ad}
\newcommand{\Ad}{Ad}
\newcommand{\integer}{\mathbb{Z}}

\newcommand{\tensor}{\otimes}

\newcommand{\rea}{\mathbb R}
\newcommand{\com}{\mathbb C}
\newcommand{\rsymp}{Sp(n, \mathbb{R})}

\newcommand{\tr}{\text{tr} \,}

\begin{document}

\title{Contact Projective Structures}

\author{Daniel J. F. Fox} 

\address{School of Mathematics\\ Georgia Institute of Technology\\ 686 Cherry St.\\ Atlanta, GA 30332-0160, U.S.A}

\email{fox@math.gatech.edu}

\begin{abstract}
A contact path geometry is a family of paths in a contact manifold each of which is everywhere tangent to the contact distribution and such that given a point and a one-dimensional subspace of the contact distribution at that point there is a unique path of the family passing through the given point and tangent to the given subspace. A contact projective structure is a contact path geometry the paths of which are among the geodesics of some affine connection. In the manner of T.Y. Thomas there is associated to each contact projective structure an ambient affine connection on a symplectic manifold with one-dimensional fibers over the contact manifold and using this the local equivalence problem for contact projective structures is solved by the construction of a canonical regular Cartan connection. This Cartan connection is normal if and only if an invariant contact torsion vanishes. Every contact projective structure determines canonical paths transverse to the contact structure which fill out the contact projective structure to give a full projective structure, and the vanishing of the contact torsion implies the contact projective ambient connection agrees with the Thomas ambient connection of the corresponding projective structure. An analogue of the classical Beltrami theorem is proved for pseudo-hermitian manifolds with transverse symmetry. 
\end{abstract}

\maketitle
\tableofcontents

\section{Introduction}
This paper makes a thorough study of contact projective structures, certain geometric structures supported by a contact manifold, $(M, H)$. Along with conformal, CR, and projective structures; path geometries; and almost Grassmannian geometries; contact projective structures are examples of curved geometric structures modeled on a homogeneous quotient, $G/P$, where $G$ is semi-simple and $P$ is a parabolic subgroup. For contact projective structures the model $G$ is the real symplectic linear group and $P$ is the subgroup fixing a one-dimensional subspace of the defining representation of $G$, and corresponding to marking the first node on the Dynkin diagram $C_{n}$. 

The study of the geometry of such structures is facilitated by the construction from some underlying geometric data, e.g. a section of a weighted tensor bundle (as for conformal structures) or a foliation (as for generalized path geometries), of a $P$ principal bundle supporting a $(\g, P)$ Cartan connection (for background on Cartan connections see \cite{Cap-Slovak-Soucek} or \cite{Slovak}). For this there are various approaches, e.g. E. Cartan's method of equivalence, T.Y. Thomas's ambient constructions, or N. Tanaka's Lie cohomological prolongations. Theorems associating canonical regular, normal Cartan connections to broad classes of generalized $G$-structures on filtered manifolds have been proved by N. Tanaka, \cite{Tanaka-equivalence}, T. Morimoto, \cite{Morimoto}, and A. \v{Cap} - H. Schichl, \cite{Cap-Schichl}. The projective and contact projective structures are exceptional for these theorems because in these cases a certain Lie algebra cohomology fails to vanish (see \cite{Cap-Slovak} and \cite{Yamaguchi}). The local equivalence problem for projective structures was in any case solved by various methods by E. Cartan, T.Y. Thomas, and H. Weyl in the 1920's. J. Harrison, \cite{Harrison}, and \v{C}ap-Schichl solved the local equivalence problem for structures corresponding to a subclass, characterized by the vanishing of an invariant \textbf{contact torsion}, of what are here called contact projective structures. The main result of this paper is an ambient construction in the manner of T.Y. Thomas which is used to associate to each contact projective structure a canonical regular Cartan connection which is normal if and only if the invariant contact torsion vanishes. (`Regular' and `normal' are used as in \cite{Cap-Slovak-Soucek}). The need to consider non-normal Cartan connections seems the most interesting aspect of this paper. Additionally, the basic local geometry of contact projective structures is described; as an application there is proved a pseudo-hermitian analogue of the classical Beltrami theorem.

Call a smoothly immersed one-dimensional submanifold a {\bf path}. Call a path everywhere tangent to the contact distribution a {\bf contact path}.
\begin{definition}\label{imprecisedefncontactpath}
A {\bf contact path geometry} is a $(4n-5)$ parameter family of contact paths in the $(2n-1)$-dimensional contact manifold, $(M, H)$, such that for every $x \in M$ and each $L \in \proj(H_{x})$ there is a unique contact path in the family passing through $x$ and tangent to $L$. Two such families of paths are {\bf equivalent} if there is a contactomorphism mapping the paths of one family onto the paths of the other family. A \textbf{contact projective structure} is a contact path geometry the paths of which are among the unparameterized geodesics of some affine connection. Two contact projective structures are {\bf equivalent} if and only if they are equivalent as contact path geometries. The contact paths of a contact projective structure are called \textbf{contact geodesics}.
\end{definition}
\noindent
The model contact projective structure is the projectivization of a symplectic vector space with the family of contact lines comprising the images in the projectivization of the two-dimensional isotropic subspaces. A contact projective structure is flat if it is locally equivalent to this model. Every three-dimensional contact path geometry is locally equivalent to that contact path geometry determined by the solutions of a third order ODE considered modulo contact transformations, and the local equivalence problem for these structures was solved by S.-S. Chern, \cite{Chern-Third}. The study of the higher dimensional contact path geometries is an ongoing project of the author; in this paper attention is restricted to the contact projective structures. The study made here of contact projective structures is modeled on the classical study of projective structures. 
\begin{definition}\label{projectivestructuredefined}
A {\bf path geometry} on an $n$-dimensional smooth manifold, $M$, is a $(2n-2)$ parameter family of paths in $M$, such that for every $x \in M$ and each $L \in \proj(T_{x}M)$, there is a unique path of in the family passing through $x$ and tangent to $L$. Two such families of paths are {\bf equivalent} if there is a diffeomorphism of $M$ mapping the paths of one family onto the paths of the other family. A {\bf projective structure} is a path geometry the paths of which are among the unparameterized geodesics of some affine connection. Equivalently, it is an equivalence class of torsion-free affine connections, $[\nabla]$, such that the unparameterized geodesics of any two representative connections are the same. Two projective structures are equivalent if they are equivalent as path geometries.
\end{definition}

\noindent
Using his method of equivalence, Cartan, \cite{Cartan}, attached to each projective structure a unique regular, normal Cartan connection. T.Y. Thomas's construction of an `ambient' affine connection associated to a projective structure provides an alternative to the method of Cartan. On a manifold with a projective structure there is a distinguished family of torsion free affine connections parameterized by the possible choices of a volume form. The ambient manifold is the total space of the $\rea^{\times}$ principal bundle of non-vanishing $\tfrac{1}{n+1}$-densities on the $n$-dimensional base manifold. Thomas built a functor associating to each projective structure a Ricci flat torsion-free affine connection on the ambient manifold, homogeneous in the vertical direction and making parallel a canonical volume form. 

In his thesis, \cite{Harrison}, J. Harrison called by `contact projective structures' certain projective structures on a contact manifold suitably adapted to the contact structure, and showed that Thomas's ambient connection parallelized a canonical symplectic structure, and from this built a tractor connection. It is worth emphasizing that Harrison's initial data of a full projective structure, which is a $4n-4$ parameter family of paths, comprises strictly more data than the $4n-5$ parameter family of paths herein assumed given. That the projective structures studied by Harrison correspond to the subclass of contact projective structures determined by the vanishing of the invariant contact torsion is a consequence of Theorem \ref{contactadapted}, which associates to each contact projective structure a full projective structure. In Section 3.24 of \cite{Cap-Schichl}, \v{C}ap and Schichl describe efficiently and briefly a construction of a canonical regular, normal Cartan connection in a setting which corresponds also to the case of contact projective structures with vanishing contact torsion. They formulate the input data as a choice of a section of the subbbundle comprising elements with co-closed torsion of a certain vector bundle, constructed in section 3.13 of \cite{Cap-Schichl}, over an abstractly given $P$ principal bundle. This choice determines a harmonic $P$ frame bundle of length two to which a general prolongation procedure is applied. It also determines a class of partial connections on the contact distribution, from which may be constructed a contact projective structure which necessarily has vanishing contact torsion. It is not clear whether this construction could be generalized to encompass contact projective structures with non-vanishing contact torsion. 

The rest of the introduction describes in more detail the results of this paper. Attached to each contact projective structure is the equivalence class of affine connections comprising those connections each of which has among its unparameterized geodesics the given contact geodesics. The choice of contact one-form is refered to as a choice of \textbf{scale}. Such a choice is analogous to the choice of a representative metric on a conformal manifold or a pseudo-hermitian representative on a CR manifold. Theorem A is a direct analogue of results about conformal or CR structures, associating to each choice of scale a unique affine connection. 
\begin{theoremai}
On a contact manifold, let a contact projective structure be given. Choose a contact one-form, $\theta$, and let $\rb$ denote its Reeb vector field. There exists a unique affine connection, $\nabla$, with torsion tensor, $\tau$, having among its unparameterized geodesics the given contact geodesics and satisfying:
\begin{conditionai}\label{cona1}
$\nabla \theta = 0$.
\end{conditionai}
\begin{conditionai}\label{cona2}
$\nabla d\theta = 0$.
\end{conditionai}
\begin{conditionai}\label{cona3}
$i(\rb)\tau = 0$.
\end{conditionai}
\begin{conditionai}\label{cona4}
The trace of the torsion of $\nabla$ vanishes. 
\end{conditionai}
\end{theoremai}
\noindent
As the torsion is a $\binom{1}{2}$-tensor, it may be contracted on its up index and either down index. Because the torsion is skew-symmetric in its down indices, condition \conref{cona4} does not depend of the choice of down index, and the interior multiplication in condition \conref{cona3} may be taken in either index. The connection, $\nabla$, of Theorem A preserves $H$ and the symplectic structure induced on $H$ by $d\theta$. 

A piece of the torsion of the connection associated by Theorem A to the given choice of scale, essentially the torsion in the contact directions, does not depend on the choice of scale, and is therefore an invariant of the contact projective structure. In general, this contact torsion need not vanish, although it necessarily vanishes in dimension three. In Theorem \ref{spaceofstructures} there is constructed a contact projective structure with contact torsion prescribed in a open neighborhood, and the space of contact projective structures with a given contact torsion is parameterized explicitly by sections of the bundle of completely symmetric covariant three tensors on the contact distribution. In particular there exist locally infinitely many inequivalent contact projective structures with vanishing contact torsion. When the contact torsion vanishes there is a contact projective Weyl curvature invariant under change of scale. In dimension at least $5$ the contact torsion and the contact projective Weyl tensor are the complete obstruction to local equivalence to the flat model. In dimension three the contact torsion and the contact projective Weyl tensor vanish identically and there is an invariant tensor, analogous to the Cotton tensor of three-dimensional conformal geometry, which is the complete obstruction to local flatness. 

The classical Beltrami theorem states that the projective structure determined by the geodesics of a Riemannian metric is flat if and only if the metric has constant sectional curvature. The contact geodesics of a pseudo-hermitian structure generate a contact projective structure, and Theorem \ref{phbeltrami} shows that the contact projective structure induced by an integrable pseudo-hermitian structure with transverse symmetry is flat if and only if the Webster sectional curvatures are constant.

Theorem B is an analogue for contact projective structures of the ambient construction of T.Y. Thomas. In the theorem the `ambient' manifold, $\form$, is a square-root of the bundle of positive contact one-forms on the co-oriented contact manifold, $(M, H)$; $\eul$ is the vertical vector field generating the dilations in the fibers of $\form \to M$; $\alpha$ is the tautological one-form on $\form$; and $\Omega = d\alpha$ is the canonical symplectic structure on $\form$. Any (local) section, $s:M \to \form$, determines a horizontal lift, $\hat{X} \in \Gamma(T\form)$, of a vector field $X \in \Gamma(TM)$. The curvature tensor of an affine connection, $\hat{\nabla}$, on $\form$, is denoted $\hat{R}_{IJK}\,^{L}$, and indices are raised and lowered with $\Omega_{IJ}$.
\begin{theoremb}
Let $(M, H)$ be a co-oriented contact manifold and let $\rho:\form \to M$ be a square-root of the bundle of positive contact one-forms. There is a functor associating to each contact projective structure on $M$ a unique affine connection, $\hat{\nabla}$, (the {\bf ambient connection}), on the total space of $\form$, having torsion $\hat{\tau}$, and satisfying the following conditions:
\begin{conditionb}\label{con1}$\hat{\nabla}\eul$ is the fundamental $\binom{1}{1}$-tensor on $\form$.\end{conditionb}
\begin{conditionb}\label{con2}$i(\eul)\hat{\tau} = 0$.\end{conditionb}
\begin{conditionb}\label{con3}$\hat{\nabla}\Omega = 0$.\end{conditionb} 
\begin{conditionb}\label{con4}The Ricci trace of the curvature tensor of $\hat{\nabla}$ vanishes.\end{conditionb} 
\begin{conditionb}\label{con5}There vanishes the restriction to $\ker \alpha$ of the tensor $\hat{R}_{Q}\,^{Q}\,_{IJ}$.\end{conditionb} 
\begin{conditionb}\label{con6}For any (local) section, $s:M \to \form$, the affine connection, $\bar{\nabla}$, on $M$ defined by $\bar{\nabla}_{X}Y = \rho_{\ast}(\hat{\nabla}_{\hat{X}}\hat{Y})$ represents the given contact projective stucture.\end{conditionb} 
\noindent
Moreover, the contact projective structures with vanishing contact torsion are in bijection with the torsion free affine connections satisfying conditions \conref{con1}, \conref{con3}, \conref{con4}, and \conref{con6}. In this case condition \conref{con2} is vacuous, \conref{con5} follows from \conref{con4} by the contracted first Bianchi identity, and the curvature tensor is completely trace free. 
\end{theoremb}

A projective structure on a contact manifold is said to be subordinate to a contact projective structure if the contact geodesics of the contact projective structure are among the geodesics of the projective structure. Theorem \ref{contactadapted} shows that each contact projective structure determines a canonical subordinate projective structure. If the contravariant part of the projective Weyl tensor of a subordinate projective structure takes values in the contact distribution, the projective structure is called contact adapted. Corollary \ref{contactadapted2} shows that every contact adapted projective structures arises as the canonical projective structure subordinate to a contact projective structure with vanishing contact torsion. The projective structures considered by Harrison in \cite{Harrison} are exactly the contact adapted projective structures. 

Using Theorem B and the tractor formalism of \v{C}ap - R. Gover, \cite{Cap-Gover}, there is built from the ambient connection a canonical Cartan connection associated to each contact projective structure. Each choice of $\form$ determines a tractor bundle and a canonical $P$ principal bundle, $\pi:\adaptedframe \to M$, the bundle of filtered symplectic frames in the tractor bundle. A $(\g, P)$ Cartan connection on $\adaptedframe$ is called compatible if it is regular and it induces the underlying co-oriented contact structure. In general, many compatible Cartan connections induce a given contact projective structure. The Cartan connection built from the ambient connection is a canonical representative of an isomorphism class of such connections distinguished by certain algebraic restrictions on the values taken by their curvature functions. There is defined in Section \ref{algebraic} a $P$-submodule $\pcurv\subset\Lambda^{2}((\g/\p)^{\ast})\tensor \g$ such that the following theorem holds ($\p$ is the Lie algebra of $P$).

\begin{theoremci}
On a co-oriented contact manifold, let $\form$ be a square-root of the bundle of positive contact one-forms, and let $\adaptedframe$ be the associated bundle of filtered symplectic frames in the tractor bundle. Among the compatible $(\g ,P)$ Cartan connections on $\adaptedframe$ inducing on $M$ a given contact projective structure, there is a unique isomorphism class of (regular) connections with curvature in $\pcurv$. A unique representative of this isomorphism class is distinguished by the requirement that the covariant differentiation it induces on $T\form$ be the ambient connection of the underlying contact projective structure. The underlying contact projective structure has vanishing contact torsion if and only if the associated Cartan connection is normal.
\end{theoremci} 
\noindent
The geometric data of the contact paths is regarded here as the primary structure; a canonical Cartan connection is built from them as a tool for their study. This connection is not, in general, normal, and requiring normality would not be natural because there is no obvious \textit{a priori} way to distinguish the corresponding underlying families of contact paths. When the contact torsion vanishes, the connection produced by Theorem C coincides with that constructed by \v{C}ap-Schichl in \cite{Cap-Schichl}.

Andreas \v{C}ap generously communicated his ideas about contact projective structures and he made helpful remarks about the prolongation procedure of \cite{Cap-Schichl}. I thank him for his comments which have influenced the point of view taken here. This project benefited at every stage from the guidance and criticism given the author by Robin Graham. This paper is based on part of the author's Ph.D. thesis at the University of Washington.

\section{Contact Projective Structures}\label{contactprojectivestructures}
\subsection{Preliminaries and Notational Conventions}\label{notationsection}
A {\bf contact manifold} is a smooth $(2n-1)$-dimensional manifold, $M$, with a maximally non-integrable, codimnesion $1$ distribution, $H$. A {\bf contactomorphism} is a diffeomorphism of $M$ that preserves $H$. A {\bf contact one-form} is a non-vanishing one-form, $\theta$, annihilating $H$. If an orientation is fixed on $(TM/H)^{\ast}$, $M$ is called \textbf{co-oriented}, and a choice of $\theta$ consistent with the given co-orientation is called \textbf{positive}. It will be assumed, usually without comment, that $(M, H)$ has a fixed co-orientation, and $\pi:\contactb \to M$ will denote the $\rea^{>0}$ principal bundle of positive contact one-forms. Each choice of $\theta$ determines uniquely a {\bf Reeb vector field} characterized by $\theta(\rb) = 1$ and $i(\rb)d\theta = 0$. Lowercase Latin indices will run from $1$ to $2n-2$. Lowercase Greek indices will run from $0$ to $2n-2$. A coframe, $\theta^{\alpha}$, is \textbf{$\theta$-adapted} if $\theta^{0} = \theta$ and $\theta^{i}(\rb) = 0$. An adapted coframe determines a dual frame, $E_{\alpha}$, such that $E_{0} = \rb$ and the $E_{i}$ span $H$. When a contact form is fixed, an adapted coframe and corresponding dual frame will be assumed fixed also. The notations $S_{[\alpha_{1}\dots \alpha_{k}]}$ and $S_{(\alpha_{1}\dots \alpha_{k})}$ denote, respectively, the complete skew-symmetrization and the complete symmetrization over the bracketed indices. Sometimes the abstract index notation will be used, so that equations with indices have invariant meaning. Greek abstract indices label sections of tensor bundles on $M$, while Latin abstract indices label sections of the tensor powers of $H$ and $H^{\ast}$, so that an expression such as $\tau_{[ij]}\,^{k}$ indicates a section of $\Lambda^{2}(H^{\ast})\tensor H$. Each $\theta$ determines a splitting, $TM = H \oplus \text{span}\{\rb\}$, which induces a splitting of the full tensor bundle. Using these splittings Latin abstract indices may be interpreted as the components of a tensor with respect to a $\theta$-adapted coframe and dual frame. The components of $\omega = d\theta$ are $\omega_{\alpha\beta} =  \omega_{[\alpha\beta]} = \omega(E_{\alpha}, E_{\beta})$. As $\omega_{0\alpha} = 0$, $\omega$ may be written as $\omega = \frac{1}{2}\omega_{ij}\theta^{i}\wedge \theta^{j}$. Latin indices may be raised and lowered using $\omega_{ij}$ according to the following conventions. Defining $\omega^{kl}$ by $\omega^{kl}\omega_{lj} = -\delta_{j}\,^{k}$, let $\gamma^{p} = \omega^{pq}\gamma_{q}$, and $\gamma_{p} = \gamma^{q}\omega_{qp}$. It is necessary to pay attention to which index is raised or lowered as, for instance, $\eta^{p}\gamma_{p} = -\eta_{p}\gamma^{p}$. There vanishes the trace, using $\omega^{pq}$, over any pair of indices in which a tensor is symmetric. As, for any $S^{i_{1} \dots i_{2k+1}}$, the contraction, $S_{pq} = S^{i_{1} \dots i_{2k}}\,_{p}S_{i_{1} \dots i_{2k} q}$, is symmetric, this implies $S^{i_{1} \dots i_{2k+1}}S_{i_{1} \dots i_{2k+1}} = 0$. Under a change of scale,  $\tilde{\theta} = f^{2}\theta$, ($f \neq 0$), the symplectic structure defined on $H$ by the restriction to $H$ of $\omega$ rescales by $f^{2}$, so there is a well-defined conformal symplectic structure on the contact hyperplane. The non-degeneracy of $d\theta$ on $H$ is equivalent to $\theta \wedge (d\theta)^{n-1} \neq 0$. Under rescaling this volume transforms by
\begin{equation}\label{volumetransform}
\tilde{\theta} \wedge (d\tilde{\theta})^{n-1} = f^{2n}\theta \wedge (d\theta)^{n-1}
\end{equation}
A $\tilde{\theta}$-adapted coframe and corresponding dual frame, are defined by 
\begin{align}\label{adaptedtransform}
&\tilde{\theta}^{i} = \theta^{i} - 2\gamma^{i}\theta,& &\tilde{E}_{i} = E_{i},& &f^{2}\tilde{\rb} = \rb + 2\gamma^{p}E_{p},&
\end{align}
where $\gamma = d\log{f}$. The coordinate expressions of tensors labeled with a $\tilde{\,\,}$ are taken with respect to the $\Tilde{\theta}$-adapted coframe and dual frame defined in \eqref{adaptedtransform}, unless there is made a specific indication otherwise. The expression $\tilde{\omega}_{ij}$ has an invariant meaning; it indicates the section of $\Lambda^{2}(H^{\ast})$ obtained by restricting $\tilde{\omega}$ to $H$. Alternatively, $\tilde{\omega}_{ij}$ may be viewed as the components of $\tilde{\omega}$ with respect to a $\tilde{\theta}$-adapted coframe.

Every affine connection, $\nabla$, decomposes as $\nabla_{X}Y = \symnabla_{X}Y + \tfrac{1}{2}\tau(X, Y)$, where the skew-symmetric $\binom{1}{2}$-tensor, $\tau$, is the torsion of $\nabla$, and the torsion-free affine connection, $\symnabla$, is the \textbf{symmetric part} of $\nabla$. The {\bf difference tensor}, $\Lambda$, of two affine connections, $\Bar{\nabla}$ and $\nabla$, is the $\binom{1}{2}$-tensor defined by $\Lambda(X, Y) = \Bar{\nabla}_{X}Y - \nabla_{X}Y$. The difference of the torsions of the connections is
\begin{equation}\label{skewlambda}
\Bar{\tau}_{\alpha\beta}\,^{\gamma} -\tau_{\alpha\beta}\,^{\gamma} = 2\Lambda_{[\alpha\beta]}\,^{\gamma}.
\end{equation}
Let $S^{\alpha_{1} \dots \alpha_{k}}_{\beta_{1} \dots \beta_{l}}$ be any $\binom{k}{l}$-tensor. The following equation describes the action of the difference tensor on the tensor $S^{\alpha_{1} \dots \alpha_{k}}_{\beta_{1} \dots \beta_{l}}$:
\begin{equation}\label{differenceact}
\Bar{\nabla}_{\gamma}S^{\alpha_{1} \dots \alpha_{k}}_{\beta_{1} \dots \beta_{l}} - \nabla_{\gamma}S^{\alpha_{1} \dots \alpha_{k}}_{\beta_{1} \dots \beta_{l}} = - \sum_{i = 1}^{l}\Lambda_{\gamma \beta_{i}}\,^{\sigma} S^{\alpha_{1} \dots \alpha_{k}}_{\beta_{1} \dots  \sigma \dots \beta_{l}}  + \sum_{i = 1}^{k}\Lambda_{\gamma \sigma}\,^{\alpha_{i}}S^{\alpha_{1} \dots\sigma \dots \alpha_{k}}_{\beta_{1} \dots \beta_{l}}
\end{equation}
There will be used sometimes the specialization obtained by regarding half the torsion tensor of a connection as the difference tensor of the connection and its symmetric part. Finally, for any $k$-form, $\eta$, there holds 
\begin{equation}\label{skew}
(k+1)\nabla_{[\alpha_{1}}\eta_{\alpha_{2} \dots \alpha_{k+1}]} = d\eta_{\alpha_{1} \dots \alpha_{k+1}} - \tbinom{k+1}{2}\eta_{\beta[\alpha_{1} \dots \alpha_{k-1}}\tau_{\alpha_{k} \alpha_{k+1}]}\,^{\beta}.
\end{equation}

\subsection{Canonical Representative Connections}\label{contactgeodesicsection}
A parameterized curve, $\gamma$, describes the same path as does a geodesic of the affine connection, $\nabla$, if and only if $\dot{\gamma}\wedge \nabla_{\dot{\gamma}}\dot{\gamma} = 0$. An affine connection, $\nabla$, on $(M, H)$ is said to {\bf admit a full set of contact geodesics} if every geodesic of $\nabla$ tangent to $H$ at one point is everywhere tangent to $H$.
\begin{lemma}\label{contactgeodesics}
For any contact manifold, $(M, H)$, the following are equivalent:
\setcounter{condition}{0}
\begin{condition}\label{cgcon1}
The affine connection, $\nabla$, admits a full set of contact geodesics.
\end{condition}\begin{condition}\label{cgcon2}
For any choice of one form, $\theta$, annihilating $H$, $\nabla_{(i}\theta_{j)} = 0$.
\end{condition}
\end{lemma}
\noindent
Note that the condition $\nabla_{(i}\theta_{j)} = 0$ does not depend on the choice of scale.
\begin{proof}
This follows from the observation that for any geodesic, $\gamma$, and any connection, $\nabla$, there holds $\ddot{\gamma}^{0} = \dot{\gamma}^{i}\dot{\gamma}^{j}\nabla_{(i}\theta_{j)}$. 
\end{proof}

\begin{lemma}\label{differencetensor}
Affine connections $\nabla$ and $\Bar{\nabla}$ on $(M, H)$, each admitting a full set of contact geodesics, determine the same contact projective structure if and only if their difference tensor, $\Lambda$, satisfies the following two conditions:
\setcounter{condition}{0}
\begin{condition}\label{diffcon1}
$\Lambda_{(ij)}\,^{0} = 0$.
\end{condition}\begin{condition}\label{diffcon2}
There is a smooth section, $\sigma$, of $H^{\ast}$ such that $\Lambda_{(ij)}\,^{k} = \delta_{i}\,^{k}\sigma_{j} + \delta_{j}\,^{k}\sigma_{i}$.
\end{condition}
\end{lemma}

\begin{proof}
First it will be shown that if $\nabla$ and $\Bar{\nabla}$ admit that same full set of unparameterized contact geodesics, then conditions \conref{diffcon1} and \conref{diffcon2} are satisfied by the difference tensor $\Lambda$. Choose a contact one-form $\theta$. By Lemma \ref{contactgeodesics} and \eqref{differenceact}, $\Lambda_{(ij)}\,^{0} = 0$. That the path of a parameterized curve, $\gamma$, is the path of a contact geodesic for each of $\Bar{\nabla}$ and $\nabla$ means $\dot{\gamma}\wedge \Lambda(\dot{\gamma}, \dot{\gamma}) = 0$. Since $\Bar{\nabla}$ and $\nabla$ admit the same full set of contact geodesics, this implies 
\begin{equation}\label{weylarg}
 \Lambda_{(ij}\,^{p}\delta_{k)}\,^{q} - \Lambda_{(ij}\,^{q}\delta_{k)}\,^{p} = \Lambda_{((ij)}\,^{p}\delta_{k)}\,^{q} - \Lambda_{((ij)}\,^{q}\delta_{k)}\,^{p} = 0.
\end{equation}
Tracing \eqref{weylarg} in $q$ and $k$ gives \conref{diffcon2} with $\sigma_{i} = \tfrac{1}{2n-1}\Lambda_{(ip)}\,^{p}$. For the converse, let $\gamma$ be a contact geodesic for $\nabla$ with its affine parameterization. Then $\Bar{\nabla}_{\dot{\gamma}}\dot{\gamma} = \nabla_{\dot{\gamma}}\dot{\gamma} + \Lambda(\dot{\gamma}, \dot{\gamma}) = 2\sigma(\dot{\gamma})\dot{\gamma}$, so that $\dot{\gamma} \wedge \Bar{\nabla}_{\dot{\gamma}}\dot{\gamma} = 0$, and $\gamma$ traces out the path of a contact geodesic for $\Bar{\nabla}$.
\end{proof} 
\noindent
A contact projective structure may be identified with the equivalence class, $[\nabla]$, of affine connections on $M$ having the property that each representative of the equivalence class admits the given full set of contact geodesics.

The following observations are used in the proof of Theorem A. Let $\nabla$ be an affine connection such that $\nabla \theta = 0$ and $\nabla\omega = 0$. \eqref{skew} applied to $\nabla\theta$ and to $\nabla\omega$, respectively, implies
\begin{align}
& \tau_{\alpha\beta}\,^{0} = \omega_{\alpha\beta},& &\label{torsionlemma1} \omega_{\sigma[\gamma}\tau_{\alpha\beta]}\,^{\sigma} = 0,
\end{align}
\noindent
respectively. Tracing \eqref{torsionlemma1} shows $2\tau_{ip}\,^{p} = -\tau_{p}\,^{p}\,_{i}$. This shows that there is, up to a scalar factor, only one section of $H^{\ast}$ obtainable by tracing the torsion. Applied to $\nabla \omega = 0$, \eqref{differenceact} and \eqref{torsionlemma1} give $\symnabla_{\gamma}\omega_{\alpha\beta} = \tfrac{1}{2}\omega_{\gamma \sigma}\tau_{\alpha\beta}\,^{\sigma}$. Covariantly differentiating the conditions, $\theta(\rb) = 1$ and $i(\rb)\omega = 0$, defining $\rb$ shows that $\rb$ is parallel. By definition of the torsion tensor, $\nabla_{\rb} = \lie_{\rb}$ if and only if $i(\rb)\tau = 0$.

\begin{proof}[Proof of Uniqueness in Theorem A]
Fix a $\theta$-adapted coframe. Suppose given representatives, $\nabla$ and $\Bar{\nabla}$, of $[\nabla]$ satisfying conditions \conref{cona1}-\conref{cona4}, and having difference tensor $\Lambda$. \eqref{differenceact}, $\nabla\rb = 0$, and \eqref{skewlambda} show $\Lambda_{\alpha \beta}\,^{0} = 0$, $\Lambda_{\beta 0}\,^{\alpha} = 0$, and $\Lambda_{0 \beta}\,^{\alpha} = 0$. By \eqref{differenceact}, $0 = \Bar{\nabla}_{k}\omega_{ij} - \nabla_{k}\omega_{ij} = -\Lambda_{k[ij]}$, and so also $\Lambda_{ip}\,^{p} = 0$. Since $\Bar{\tau}_{ip}\,^{p} - \tau_{ip}\,^{p} =0$, \eqref{skewlambda} implies $\Lambda_{[ip]}\,^{p} = 0$, and so $\Lambda_{(ip)}\,^{p} = \Lambda_{ip}\,^{p} - \Lambda_{[ip]}\,^{p} = 0$. By Lemma  \ref{differencetensor}, $\Lambda_{(ij)}\,^{k} = \frac{1}{2n-1}(\delta_{i}\,^{k}\Lambda_{(jp)}\,^{p} + \delta_{j}\,^{k}\Lambda_{(ip)}\,^{p}) = 0$. Any $\binom{0}{3}$-tensor skew in two indices and symmetric in two indices vanishes, so $\Lambda_{ijk} = 0$.
\end{proof}

\begin{proof}[Proof of Existence in Theorem A]
Given any representative, $\nabla \in [\nabla]$, and a choice of $\theta$, $\nabla_{(i}\theta_{j)} = 0$. A connection in $[\nabla]$ satisfying the conditions \conref{cona1}-\conref{cona4} is produced by succesive modifications of $\nabla$. In each of the following claims, a connection $\Bar{\nabla}$ is obtained from a given connection, $\nabla$, by specifying their difference tensor, $\Lambda_{\alpha\beta}\,^{\gamma}$. The proofs of the claims are straightforward computations using \eqref{skewlambda}, \eqref{differenceact}, \eqref{skew}, \eqref{torsionlemma1} and Lemma \ref{differencetensor} so are omitted.

\begin{claim}
Assume $\nabla \in [\nabla]$. Define $\Lambda_{\alpha\beta}\,^{\gamma} = \delta_{0}\,^{\gamma}\nabla_{\alpha}\theta_{\beta}$. Then $\Bar{\nabla} \in [\nabla]$ and $\bar{\nabla}\theta = 0$. Thus it may be assumed from the beginning that $\nabla\theta = 0$.
\end{claim}

\begin{claim}
Suppose given $\nabla \in [\nabla]$ such that $\nabla \theta = 0$. Define $\Lambda_{\alpha\beta}\,^{\gamma} = -\delta_{\beta}\,^{0}\nabla_{\alpha}\rb^{\gamma}$. Then $\Bar{\nabla} \in [\nabla]$; $\Bar{\nabla}\theta = 0$; and $\Bar{\nabla}\rb = 0$. 
\end{claim}

\begin{claim}\label{claim3}
Suppose given $\nabla \in [\nabla]$ with torsion, $\tau$, such that $\nabla \theta = 0$ and $\nabla \rb = 0$. Define $\Lambda_{\alpha\beta}\,^{\gamma} = -\delta_{\alpha}\,^{0}\tau_{0\beta}\,^{\gamma}$. Then $\Bar{\nabla} \in [\nabla]$; $\Bar{\nabla}\theta =0$; $\Bar{\nabla}\rb = 0$; the torsion, $\Bar{\tau}$, of $\Bar{\nabla}$ satisfies $i(\rb)\Bar{\tau} = 0$; and $\Bar{\nabla}_{\rb}\omega = 0$.
\end{claim}

\begin{claim}\label{claim4new}
Suppose given $\nabla \in [\nabla]$ such that $\nabla\theta =0$; $\nabla\rb = 0$; the torsion, $\tau$, of $\nabla$ satisfies $i(\rb)\tau = 0$; and $\nabla_{\rb}\omega = 0$. Define $\Lambda_{\alpha\beta}\,^{0} = 0$ and $\Lambda_{\alpha\beta k} = -\nabla_{k}\omega_{\alpha\beta} + \tfrac{3}{2}\tau_{[\alpha\beta}\,^{p}\omega_{k]p}$. Then $\Bar{\nabla} \in [\nabla]$; $\Bar{\nabla}\theta = 0$; $\Bar{\nabla}\omega = 0$; and the torsion, $\Bar{\tau}$, of $\Bar{\nabla}$, satisfies $i(\rb)\Bar{\tau} = 0$.
\end{claim}

\begin{claim}\label{claim5}
Suppose $\nabla \in [\nabla]$ with torsion $\tau$ is chosen so that $\nabla \theta = 0$, $\nabla \omega = 0$ and $i(\rb)\tau = 0$. Let $(2n-1)\gamma_{i} = \tau_{ip}\,^{p}$ and define $\Lambda$ by
\begin{align*}
&\Lambda_{ij}\,^{p} = \gamma_{j}\delta_{i}\,^{p} + \omega_{ij}\gamma^{p},&  &\Lambda_{\alpha \beta}\,^{0} = 0,& &\Lambda_{0\beta}\,^{\alpha} = 0 = \Lambda_{\beta 0}\,^{\alpha}.&
\end{align*} 
Then $\Bar{\nabla} \in [\nabla]$ satisfies conditions \conref{cona1}-\conref{cona4} in the statement of Theorem A.
\end{claim}
\noindent
[End of Proof of Theorem A]
\end{proof}

\subsection{Contact Torsion}\label{transformationrules}
\begin{lemma}\label{transformationruleslemma}
Given a contact projective structure, let $\Lambda$ be the difference tensor of the representatives, $\tilde{\nabla}$ and $\nabla$, associated by Theorem A to the choices of contact one-forms, $\tilde{\theta} = f^{2}\theta$ and $\theta$. With respect to a $\theta$-adapted coframe and dual frame, the components of $\Lambda$ are expressible in terms of $\gamma = d\log{f}$ as
\begin{align}
&\label{lambdatransform}\Lambda_{ij}\,^{p} = \gamma_{i}\delta_{j}\,^{p} + \gamma_{j}\delta_{i}\,^{p} + \omega_{ij}\gamma^{p},\\
&\Lambda_{\alpha \beta }\,^{0} = 2\gamma_{\alpha}\delta_{\beta}\,^{0},&
&\label{lambdai0p}\Lambda_{i0}\,^{j} = 4\gamma_{i}\gamma^{j} - 2\nabla_{i}\gamma^{j},\\
&\label{lambda0ip}\Lambda_{0i}\,^{j} = -2\gamma^{q}\tau_{qi}\,^{j} - 2\nabla_{i}\gamma^{j},&
&\Lambda_{00}\,^{i} = -2\nabla_{0}\gamma^{i} +4\gamma_{0}\gamma^{i} +4\gamma^{q}\nabla_{q}\gamma^{i}.
\end{align}
\end{lemma}

\begin{proof}
Let $\nabla$ be the connection associated to $\theta$ by Theorem A and define $\tilde{\nabla}$ by requiring that its difference tensor with $\nabla$ be given by \eqref{lambdatransform}-\eqref{lambda0ip}. Routine computations show that the $\tilde{\nabla}$ so defined satisfies the conditions of Theorem A with respect to $\tilde{\theta}$, and so is the unique connection associated to $\tilde{\theta}$ by Theorem A. 
\end{proof}

\noindent
Lemma \ref{transformationruleslemma} is useful for computing the transformation of covariant derivatives under change of scale. For instance, for a section, $\sigma_{i}$, of $H^{\ast}$. 
\begin{equation}\label{nablasigmatransform}
\tilde{\nabla}_{i}\sigma_{j} = \nabla_{i}\sigma_{j} - \sigma_{i}\gamma_{j} - \gamma_{j}\sigma_{i} - \sigma_{p}\gamma^{p}\omega_{ij}.
\end{equation}

\begin{theorem}\label{contacttorsiontheorem}
Given a contact projective structure let $\nabla$ be the affine connection associated by Theorem A to the contact one-form, $\theta$. The torsion, $\tau$, of $\nabla$ determines a section, $\tau_{ij}\,^{k}$, of $\bigwedge^{2}(H^{\ast})\tensor H$, that is independent of the choice of scale, $\theta$. The section, $\tau_{ij},^{k}$, will be called the {\bf contact torsion} of the contact projective structure.
\end{theorem}

\begin{proof}
By \eqref{adaptedtransform} and \eqref{skewlambda}, $2\Lambda_{[ij]}\,^{k} = \tilde{\tau}_{ij}\,^{k} + 2\gamma^{k}\omega_{ij} - \tau_{ij}\,^{k}$, and by \eqref{lambdatransform}, $2\Lambda_{[ij]}\,^{k} = 2\gamma^{k}\omega_{ij}$, so $\tilde{\tau}_{ij}\,^{k} = \tau_{ij}\,^{k}$.
\end{proof}

\begin{proposition}\label{threectf}
In three dimensions $\tau_{ij}\,^{k} = 0$ identically.
\end{proposition}

\begin{proof}
A trace free tensor with the symmetries of the contact torsion, $\tau_{[ij]k} = \tau_{ijk}$ and $\tau_{[ijk]} = 0$, lies in the irreducible representation of $Sp(n-1, \rea)$ corresponding to the Young diagram determined by the partition $(2, 1)$, which, by the criterion of Weyl, (\cite{Weyl-Classical}, Section $6.3$, p. $175$), is non-trivial if and only if $2n-1 > 3$.
\end{proof}

\begin{remark}
One possibly interesting class of paths determined by a contact projective structure comprises those contact paths such that for any choice of parameterization, $\gamma(t)$, and any choice of contact one-form, $\theta$, if $\nabla$ is the connection associated by Theorem A to $\theta$, then $\omega(\nabla_{\dot{\gamma}}\dot{\gamma}, \dot{\gamma}) = 0$. Using \eqref{lambdatransform} it is easy to check that this condition depends on neither the choice of $\theta$ nor the choice of parameterization. In three dimensions any such path is evidently a contact geodesic, but in higher dimensions this need not be the case, as it is easy to give examples in the flat model of such paths which are not lines.
\end{remark}

\subsection{Flat Model Contact Projective Structure}\label{flatmodel}
The symplectic group, $G = \rsymp$, is the group of linear automorphisms of a real $2n$-dimensional symplectic vector space, $(\standrep, \Omega)$, to be referred to as the {\bf standard representation}. Let $P \subset G$ be the stabilizer of the one-dimensional subspace, $\standrep^{2}\subset \standrep$, and let $\tilde{P} \subset P$ be the stabilizer of a vector, $v_{\infty}$, spanning $\standrep^{2}$. The $\Omega$-complement, $\standrep^{1} = (\standrep^{2})^{\perp}$, gives a filtration $\standrep^{2} \subset \standrep^{1} \subset \standrep$. Fixing $v_{0} \in \standrep$ such that $\Omega(v_{\infty}, v_{0}) = 1$ determines a splitting $\standrep = \standrep_{0}\oplus\standrep_{-1}\oplus\standrep_{-2}$, where $\standrep_{0} = \standrep^{2}$ and $e_{0}$ spans $\standrep_{-2}$. This grading of $\standrep$ induces on the Lie algebra, $\g$, of $G$, a $\integer$-grading $\g = \oplus_{i = -2}^{2}\g_{i}$ determined by $\g_{k} = \{g \in \g: g\cdot \standrep_{i} \subset \standrep_{i+k} \,\, \text{for}\,\, 0 \geq i \geq -2\}$. The subspace $\g_{-2}$ has dimension $1$, and the subspace $\g_{-1}$ has dimension $2n-2$. The subalgebra $\g_{-} = \g_{-2} \oplus \g_{-1}$ is the usual Heisenberg Lie algebra, and $\g_{0}$ is the conformal symplectic Lie algebra. With $v_{0}$ and $v_{\infty}$, an arbitrary basis, $v_{i}$, of $\standrep_{-1}$ determines coordinates, $u^{I}$, where $I \in \{\infty, 1, \dots, 2n-2, 0\}$, so that $\Omega = \frac{1}{2}\Omega_{IJ}du^{I}\wedge du^{J}$ with 
\begin{equation}\label{standardOmega}
\Omega_{\infty 0} = 1, \quad \Omega_{\infty i} = 0 = \Omega_{0i}, \quad \Omega_{ij} = \omega_{ij}.
\end{equation}
Raise and lower indices with $\Omega^{IJ}$ according to the convention $\Omega^{IQ}\Omega_{QJ} = -\delta_{J}\,^{I}$. Letting $v^{I}$ be a basis of $\standrep^{\ast}$ dual to $v_{I}$, each element of $G$ may be represented as $A_{I}\,^{J}v_{J}\tensor v^{I}$, where $A_{I}\,^{J}$ is a constant matrix such that $A_{I}\,^{P}A_{J}\,^{Q}\Omega_{PQ} = \Omega_{IJ}$. $P$ is a semidirect product of the subgroup, $P^{+} \subset G$, with Lie algebra $\p^{+} = \g_{1} \oplus \g_{2}$, and the subgroup, $G_{0}\subset G$, with Lie algebra $\g_{0}$. These subgroups have the forms:
\begin{align*}
 & G_{0} = \left\{ \begin{pmatrix} a & 0 & 0 \\ 0 & B_{i}\,^{j} & 0 \\ 0 & 0 & a^{-1}\end{pmatrix}: \,\, \begin{matrix}a \in \rea^{\times}, \\ B \in Sp(n-1, \rea)\end{matrix} \right\},& &P^{+} = \left\{\begin{pmatrix} 1 & \gamma_{i} & \gamma_{0} \\ 0 & \delta_{i}\,^{j} & \gamma^{j} \\ 0 & 0 & 1\end{pmatrix}\right\}. &
\end{align*}
The normal subgroup $\tilde{P} \subset P$ is distinguished by the requirement $a = 1$. $G$ acts transitively but not effectively on $\proj(\standrep)$, and its projectivization is denoted $\bar{G} = P\rsymp$. The image in $\bar{G}$ of a subgroup $H \subset G$ is denoted by $\bar{H}$. The $\rea^{\times} = P/\tilde{P}$ principal bundle, $G/\tilde{P} \to G/P$, recovers the defining bundle, $\standrepz = \standrep \smallsetminus \{0\} \to \proj(\standrep)$. The Maurer-Cartan form, $\omega_{G}$, on $G$ is the model for the Cartan connection associated by Theorem C to a general contact projective structure.

An $\Ad(P)$-invariant filtration of $\g$ is defined by $\f_{i} = \oplus_{j \geq i}\g_{j}$. The Lie algebra of $\tilde{P}$ is $\tilde{\p} = [\p, \p]$. Let $e_{\alpha}$ be a basis of $\g_{-}$, with $e_{0}$ a basis for $\g_{-2}$, and satisfying $[e_{i}, e_{j}] = -2\omega_{ij}e_{0}$ and $[e_{i}, e_{0}] = 0$. Because $\g$ is semisimple every derivation of $\g$ is inner and so there is a unique element $e_{\infty} \in \g$ such that $[e_{\infty}, \g_{i}] = i\g_{i}$. Fix an invariant bilinear form, $B$, on $\g$, let $e^{0} \in \g_{2}$ be $B$-dual to $e_{0}$, and note that $\bar{\alpha} = \tfrac{1}{2}B(e^{0}, -)$ annihilates $\f_{-1}$. View $\bar{\alpha}$ as an element of $C^{1}(\g, \rea)$, where $(C^{k}(\g, \rea), \partial)$ is the cochain complex for the cohomology of $\g$ acting (trivially) on $\rea$. Then $\bar{\omega}(e_{\alpha}, e_{\beta}) = \partial\bar{\alpha}(e_{\alpha}, e_{\beta}) = -\tfrac{1}{2}B(e^{0}, [e_{\alpha}, e_{\beta}]) = \omega_{\alpha\beta}$, and $\bar{\omega}$ restricts to give a symplectic structure on $\g_{-1}$. Because $\bar{\alpha}$ annihilates $\f_{-1}$, the symplectic form on $\f_{-1}/\p$ defined by $\bar{\omega}(e_{i} + \p, e_{j} + \p) = \bar{\omega}(e_{i}, e_{j})$ makes sense. By invariance of $B$, the adjoint action of $P$ on $\g/\p$ preserves $\bar{\omega}$ up to a positive scalar factor, so this determines a $P$-invariant positive conformal symplectic structure on $\f_{-1}/\p$. The requirement $[e_{i}, e_{j}] = -2\omega_{ij}e_{0}$ selects an $\Ad(P)$-invariant orientation on $\g/\f_{-1}$. By invariance of $B$, $\bar{\omega}(e_{I}, e_{J})$ is $\tilde{P}$-invariant, so $\bar{\omega}$ determines a symplectic form, $\bar{\Omega}$, on $\g/\tilde{\p}$. Note that $\bar{\Omega}(e_{\infty}, e_{0}) = 1$, and the map $h + \tilde{\p} \to h\cdot v_{\infty}$ defines a symplectic $\tilde{P}$-module isomorphism, $(\g/\tilde{\p}, \bar{\Omega}) \simeq (\standrep, \Omega)$, mapping $e_{I}+\tilde{\p} \to v_{I}$.

Let $E_{i}$ and $E_{0}$ be the left-invariant vector fields on $G/P$ generated by $e_{i}$ and $2e_{0}$, respectively. In the natural coordinates on $\proj(\standrep) = G/P$, these are given by $E_{i} = \frac{\partial}{\partial x^{i}} + \omega_{i p}x^{p}\frac{\partial}{\partial x^{0}}$, and $E_{0} = 2\frac{\partial}{\partial x^{0}}$. Writing $\omega_{\alpha 0} = 0 = \omega_{0 \beta}$, the Lie brackets are $[E_{\alpha}, E_{\beta}] = -\omega_{\alpha\beta}E_{0}$. The canonical contact structure on $T(G/P)$ is the left-invariant subbundle of $T(G/P)$ generated by $\g_{-1}$ and spanned by the vector fields $E_{i}$. A left-invariant section, $\Theta = \tfrac{1}{2}(dx^{0} + \omega_{pq}x^{p}dx^{q})$, of the annihilator of $H$ is determined by the requirement that $\rb = E_{0}$ be its Reeb vector field. Define a connection, $\nabla$, by requiring the left-invariant frame $E_{\alpha}$ to be parallel. Such a connection has necessarily torsion $\tau = d\Theta \tensor \rb$. It is easily checked directly that $\nabla \Theta = 0$ and $\nabla d\Theta = 0$, so $\nabla$ is a connection satisfying the conditions of Theorem A. Any left-invariant vector field is a constant coefficient linear combination of $E_{\alpha}$, and is consequently also $\nabla$-parallel. As for any affine connection the integral curves of a parallel vector field are geodesics, the integral curves of any left-invariant vector field are geodesics for $\nabla$. Explicit computation shows that these integral curves are straight lines in $\rea^{2n-1}$. The contact projective structure determined by $\nabla$ is the model for all contact projective structures. A contact projective structure is {\bf flat} if and only if it is locally equivalent to this model.

On a contact manifold with a chosen $\theta$, a coordinate chart $\psi:U \to \rea^{2n-1}$ such that $\psi(p) = 0$ and $\psi^{\ast}(\Theta) = \theta$, is called a {\bf Darboux coordinate chart} at $p$. The Darboux theorem shows that there exists always a Darboux coordinate chart at $p$. 
\begin{proposition}\label{admitscproj}
A co-oriented contact manifold admits a contact projective structure.
\end{proposition}

\begin{proof}
Fix a contact one-form, $\theta$, and cover $M$ by an atlas, $\{U_{a}\}$, of $\theta$-Darboux coordinate charts. In each $U_{a}$ let $\nabla^{a}$ be the representative described above of the flat contact projective structure associated to $\theta|_{U_{a}}$. Choose a partition of unity, $\phi_{a}$, subordinate to $\{U_{a}\}$ and define $\nabla_{X}Y = \sum_{a}\phi_{a}\nabla^{a}_{X}Y$. It is straightforward to check that this defines a connection and $\nabla_{(i}\theta_{j)} = \sum_{a}\phi_{a}\nabla^{a}\,_{(i}\theta_{j)} = 0$, so that the geodesics of $\nabla$ generate a contact projective structure.
\end{proof}

\subsection{The Affine Space of Contact Projective Structures}\label{affinespace}
There is a well defined notion of the difference tensor of two contact projective structures, $[\nabla]$ and $[\bar{\nabla}]$, as a section of $\tensor^{2}(H^{\ast})\tensor H$. Fix a contact one-form, $\theta$, and let $\Pi$ be the difference tensor of the representatives, $\bar{\nabla} \in [\bar{\nabla}]$ and $\nabla \in [\nabla]$, associated to $\theta$ by Theorem A. Each of $\bar{\nabla}$ and $\nabla$ makes parallel $\theta$ and $\rb$, and the interior multiplication of $\rb$ in the torsion of each vanishes. In conjunction with \eqref{skewlambda} and \eqref{differenceact}, these give 
\begin{equation}\label{aff0}
 \Pi_{\alpha\beta}\,^{0} = \Pi_{\alpha 0}\,^{\beta} = \Pi_{0\alpha}\,^{\beta} = 0,
\end{equation}
so that $\Pi$ may be identified with the section, $\Pi_{ij}\,^{k}$, of $\tensor^{2}(H^{\ast})\tensor H$. Let $\Tilde{\Pi}$ be the difference tensor of the representatives, $\Tilde{\bar{\nabla}} \in [\bar{\nabla}]$ and $\tilde{\nabla} \in [\nabla]$, associated to $\tilde{\theta}$ by Theorem A. Keeping in mind that $\tilde{\Pi}_{\alpha\beta}\,^{\gamma}$ denote the components of $\tilde{\Pi}$ with respect to a $\tilde{\theta}$-adapted coframe and dual frame, as in \eqref{adaptedtransform}, observe that, as in \eqref{aff0}, $\tilde{\Pi}_{\alpha\beta}\,^{0} = \tilde{\Pi}_{\alpha 0}\,^{\beta} = \tilde{\Pi}_{0\alpha}\,^{\beta} = 0$. As a consequence, the components of $\tilde{\Pi}_{ij}\,^{k}$ are the same when calculated in a $\theta$-adapted coframe and dual frame as when calculated in a $\tilde{\theta}$-adapted coframe and dual frame. It is now claimed that $\tilde{\Pi}_{ij}\,^{k} = \Pi_{ij}\,^{k}$. Let $\bar{\Lambda}$ be the difference tensor of $\tilde{\bar{\nabla}}$ and $\bar{\nabla}$, let $\Lambda$ be the difference tensor of $\tilde{\nabla}$ and $\nabla$, and observe that $\tilde{\Pi}-\Pi = \bar{\Lambda} - \Lambda$. \eqref{lambdatransform} shows that $\bar{\Lambda}_{ij}\,^{k} = \Lambda_{ij}\,^{k}$ , so that $\tilde{\Pi}_{ij}\,^{k} = \Pi_{ij}\,^{k}$. Hence $\Pi_{ij}\,^{k}$ is independent of the choice of $\theta$, and, consequently, it makes sense to speak of $\Pi_{ij}\,^{k}$ as the difference tensor of the contact projective structures. 

The facts about the irreducible representations of the symplectic group used in the sequel may be found in some form in Section 6.3 of \cite{Weyl-Classical} or in \cite{Fulton-Harris}. Let $\At$, $\B$, and $\C$, respectively, be the bundles of tensors on $H$ obtained by raising the third index of elements of, respectively, $S^{3}(H^{\ast})$; the subbundle of $\tensor^{3}(H^{\ast})$ comprising trace-free tensors satisfying $B_{i(jk)} = B_{ijk}$ and $B_{(ijk)} = 0$; and the subbundle of $\tensor^{3}(H^{\ast})$ comprising trace-free tensor satisfying $C_{i[jk]} = C_{ijk}$, and $C_{[ijk]} = 0$. Though the operation of raising an index depends on the choice of contact one-form, the bundles so defined do not. By Weyl's results the fiber over a point of any of $\At$, $\B$, or $\C$ is an irreducible $Sp(n-1, \rea)$-module. Proposition \ref{admitscproj} shows that the space of contact projective structures on a co-oriented contact manifold is non-empty. Theorem \ref{spaceofstructures} describes the affine structure on this space. 
\begin{theorem}\label{spaceofstructures}
The space of contact projective structures is an infinite-dimensional affine space modeled on $\Gamma(\At \oplus \B)$, in the sense that the difference tensor, $\Pi_{ij}\,^{k}$, of two contact projective structures on $M$ admits a direct sum decomposition, $\Pi_{ij}\,^{k} = A_{ij}\,^{k} + B_{ij}\,^{k}$, as a section of $\At \oplus \B$. Moreover,
\setcounter{condition}{0}
\begin{condition}
The difference of the contact torsions of two contact projective structures is the image of $B_{ij}\,^{k}$ under the isomorphism of tensor bundles, $\B \to \C$, determined by $B_{ij}\,^{k} \to 2\Pi_{[ij]}\,^{k} = 2B_{[ij]}\,^{k}$. In particular, the difference tensor of two contact projective structures with the same contact torsion is the section, $A_{ij}\,^{k}$, of $\At$.
\end{condition}\begin{condition}
Given $p \in M$ there is an open $U \subset M$, containing $p$, so that for any section, $\tau_{ij}\,^{k} \in \Gamma(\C)$, defined over $U$, there exists in $U$ a contact projective structure with contact torsion $\tau_{ij}\,^{k}$.
\end{condition}
\end{theorem}

\begin{proof}[Proof of Theorem \ref{spaceofstructures}]
Let $\Pi$ be the difference tensor of the contact projective structures, $[\bar{\nabla}]$ and $[\nabla]$, and let $\bar{\nabla} \in [\bar{\nabla}]$ and $\nabla \in [\nabla]$ be the representatives associated to $\theta$ by Theorem A. By \eqref{differenceact}, $0 = \bar{\nabla}_{k}\omega_{ij} - \nabla_{k}\omega_{ij} = -2\Pi_{k[ij]}$, which shows $\Pi_{ip}\,^{p} = 0$. From $ \tfrac{1}{2}(\bar{\tau}_{ijk} - \tau_{ijk}) = \Pi_{[ij]k}$ there follow $\Pi_{p}\,^{pk} = 0$ and $\Pi_{[ijk]} = 0$, so that $\Pi_{ijk}$ is completely trace-free. By complete reducibility and Weyl's description of the irreducible $Sp(n-1, \rea)$ modules, there is a direct sum decomposition, $\Pi_{ijk} = A_{ijk} + B_{ijk}$, where $A \in \Gamma(\At)$ and $B \in \Gamma(\B)$, and such that $A_{ijk} =\Pi_{(ijk)}$ and $B_{[ij]k} = \Pi_{[ij]k}$. Moreover, the map $B_{ij}\,^{k}\to 2B_{[ij]}\,^{k}$, determines an isomorphism between the fibers over a point of $\B$ and $\C$, and this shows that $B_{ij}\,^{k}$ vanishes if and only if the two contact projective structures have the same contact torsion. 

To construct, in a neighborhood of $p \in M$, a contact projective structure, $[\bar{\nabla}]$, with contact torsion $\tau_{ij}\,^{k} \in \Gamma(\C)$, proceed as follows. Fix $\theta$ and a Darboux coordinate chart, $U$, centered on $p$, and let $\nabla$ be the connection defined in section \ref{flatmodel}. This is the representative associated by Theorem A to $\theta$ of the flat contact projective structure, $[\nabla]$. Add to $[\nabla]$ a difference tensor, $\Pi_{ij}\,^{k}$, satisfying \eqref{aff0} and $\Pi_{ijk} = \frac{2}{3}\tau_{i(jk)}$. Then $\Pi_{[ij]}\,^{k} = \tfrac{1}{2}\tau_{ij}\,^{k}$, so there results a contact projective structure with contact torsion $\tau_{ij}\,^{k}$. 
\end{proof}

\begin{remark}
The group, $\contactop(M, H)$, of positive contactomorphisms of the co-oriented contact manifold acts by pullback, $\phi^{\ast}([\nabla]) = [\phi^{\ast}(\nabla)]$, on the space of contact projective structures. For a contact projective structure represented by the equivalence class of connections $[\nabla]$, define a map $\Pi:\contactop(M) \to \Gamma(\tensor^{2}(H^{\ast})\tensor H)$ by letting $\Pi(\phi)$ be the difference tensor of $[\phi^{\ast}\nabla]$ and $[\nabla]$. By definition, $\Pi(\phi\circ \psi) = \psi^{\ast}(\Pi(\phi)) + \Pi(\psi)$, so the map $\Pi$ is a $1$-cocycle of $\contactop(M, H)$ taking values in the space of sections of $\tensor^{2}(H^{\ast})\tensor H$. It may be checked that the cocycles determined in this way by different choices of $[\nabla]$ are cohomologous. 
\end{remark}

\subsection{Curvature Identities}\label{curvaturesection}
Define the curvature tensor of the affine connection, $\nabla$, by $R(E_{\alpha}, E_{\beta})E_{\gamma}= [\nabla_{E_{\alpha}}, \nabla_{E_{\beta}}]E_{\gamma} - \nabla_{[E_{\alpha}, E_{\beta}]}E_{\gamma}=  R_{\alpha\beta\gamma}\,^{\sigma}E_{\sigma}$.
The Ricci tensor is defined by contracting on the middle index, $R_{\alpha \beta} = R_{\alpha \sigma \beta}\,^{\sigma}$. Recall the Bianchi identities for an affine connection with torsion:
\begin{align}\label{bianchi1}
R_{[\alpha\beta\gamma]}\,^{\sigma} = \nabla_{[\alpha}\tau_{\beta\gamma]}\,^{\sigma} + \tau_{\delta[\alpha}\,^{\sigma}\tau_{\beta\gamma]}\,^{\delta},&
&\nabla_{[\alpha}R_{\beta\gamma]\delta}\,^{\sigma} = \tau_{[\alpha\beta}\,^{\eta}R_{\gamma]\eta\delta}\,^{\sigma}.
\end{align}
Recall also the Ricci identity:
\begin{equation*}
(2\nabla_{[\alpha}\nabla_{\beta]} + \tau_{\alpha\beta}\,^{\delta}\nabla_{\delta})S^{\alpha_{1} \dots \alpha_{k}}_{\beta_{1} \dots \beta_{l}} =  \sum_{s = 1}^{k}S^{\alpha_{1} \dots \delta  \dots \alpha_{k}}_{\beta_{1} \dots \beta_{l}}R_{\alpha\beta\delta}\,^{\alpha_{s}}   - \sum_{s = 1}^{l}S^{\alpha_{1} \dots \alpha_{k}}_{\beta_{1} \dots \delta \dots \beta_{l}} R_{\alpha\beta\beta_{s}}\,^{\delta}
\end{equation*}
Let $\nabla$ be the connection associated by Theorem A to the choice of contact one-form, $\theta$. Applying the Ricci identity to the parallel tensors, $\theta$ and $d\theta$, gives $R_{\alpha\beta\gamma}\,^{0} = 0$ and $R_{\alpha \beta kl} = R_{\alpha \beta lk}$. Because $\nabla \rb = 0$, $R_{\alpha\beta 0}\,^{\sigma} = 0$, so $R_{\alpha 0} = R_{\alpha \sigma 0}\,^{\sigma} = 0$. Since $R_{\alpha\beta 0}\,^{\sigma} = 0$, setting $\gamma = 0$ in the first Bianchi identity, \eqref{bianchi1}, and using $\tau_{0\alpha}\,^{\sigma} = 0$ gives $2R_{0 [\alpha \beta]}\,^{\sigma} = \nabla_{0}\tau_{\alpha\beta}\,^{\sigma}$, and contracting this on $\sigma$ and $\beta$ gives $-R_{0\alpha} = \nabla_{0}\tau_{\alpha\sigma}\,^{\sigma} = 0$. This shows $R_{0\alpha} = 0$. Taking all possible traces of the Bianchi identities \eqref{bianchi1} proves the following lemma.
\begin{lemma}\label{tracedcurvature}
On a contact projective manifold let $\nabla$ be the connection associated by Theorem A to the choice of contact one-form, $\theta$. Then
\begin{align}
\label{firstbianchicontracted}&R_{p}\,^{p}\,_{ij} + 2R_{ij} = 2\nabla_{p}\tau^{p}\,_{ij} - \tau^{pq}\,_{j}\tau_{pqi}\\
\label{skewtracesecondbianchi}&(2-n)R_{0ijk} = \nabla_{p}R^{p}\,_{ijk} + \tfrac{1}{2}\nabla_{i}R_{p}\,^{p}\,_{jk} + \tau_{i}\,^{pq}R_{pqjk},\\
\label{skewfirstbianchicontracted}&R_{[ij]} = -\tfrac{1}{2}\nabla_{p}\tau_{ij}\,^{p},
\qquad\qquad\qquad\qquad R_{p}\,^{p} = 0, \\
\label{firstbianchicontractedzero}&2R_{0[ij]}\,^{p} = \nabla_{0}\tau_{ij}\,^{p},
\qquad\qquad\qquad\qquad \tau^{pqr}\tau_{pqr} = 0,\\
\label{riccitracesecondbianchi}&\nabla_{p}R_{ijk}\,^{p} + 2\nabla_{[i}R_{j]k} = -\tau_{ij}\,^{p}R_{pk} + 2\tau_{p[i}\,^{q}R_{j]ql}\,^{p} + 2R_{0[ij]k},\\
\label{skewtracesecondbianchia}& \nabla_{q}R_{p}\,^{p}\,_{i}\,^{q} - 2\nabla^{p}R_{pi} = \tau^{qp}\,_{l}R_{qpi}\,^{l},\\
\label{secondbianchizero}&\nabla_{0}R_{ijkl} + 2\nabla_{[i}R_{j]0kl} = \tau_{ij}\,^{p}R_{0pkl},
\qquad \nabla_{0}R_{p}\,^{p}\,_{ij} = 2\nabla^{p}R_{p0ij},\\
\label{secondbianchizeroriccitrace}&\nabla_{0}R_{ij} + \nabla_{p}R_{0ij}\,^{p} = \tau_{iq}\,^{p}R_{0pj}\,^{q}.
\end{align}
\end{lemma}

\begin{remark}
Lemma \ref{tracedcurvature} shows that tracing $R_{ijkl}$, $\tau_{ijk}$, and $\nabla_{i}\tau_{jkl}$ gives only the two tensors, $R_{ij}$ and $R_{p}\,^{p}\,_{ij}$, and \eqref{firstbianchicontracted} shows that when $\tau_{ijk} = 0$ these tensors coincide up to a constant factor. When $2n -1 > 3$, \eqref{skewtracesecondbianchi}, \eqref{secondbianchizero}, and \eqref{firstbianchicontractedzero} show that $R_{0ijk}$, $\nabla_{0}\tau_{ijk}$, and $\nabla_{0}R_{ijkl}$ are determined completely by $R_{ijkl}$, $\tau_{ijk}$, and their covariant derivatives in the contact directions. 
\end{remark}

\noindent
The following tensors are basic in the study of contact projective structures.
\begin{align}
&\label{pdefined}P_{ij} = \tfrac{1}{n(2n-3)}\left((n-1)R_{ij} - \tfrac{1}{2n-1}R_{[ij]} + \tfrac{1}{4}R_{p}\,^{p}\,_{ij}\right),\\
&\label{qdefined}Q_{ij} = \tfrac{1}{3-2n}\left(2R_{ij} + R_{p}\,^{p}\,_{ij} -\tfrac{4}{2n-1}R_{[ij]}\right),\\
&\label{wdefined} W_{ijk}\,^{l} = R_{ijk}\,^{l} + 2\delta_{[i}\,^{l}P_{j]k} +2\omega_{k[j}P_{i]}\,^{l} + 2\omega_{ij}P_{k}\,^{l} + \omega_{ij}Q_{k}\,^{l},\\
&\label{cdefined}C_{ijk} = R_{0ijk} - \left(2\nabla_{i}P_{jk} + \nabla_{i}Q_{jk}\right)\\ \notag & + \tfrac{1}{2n-1}(2\omega_{ik}\nabla_{p}P_{j}\,^{p} + \omega_{ik}\nabla_{p}Q_{j}\,^{p} + 2\omega_{ij}\nabla_{p}P_{k}\,^{p} + \omega_{ij}\nabla_{p}Q_{k}\,^{p}).
\end{align}
\noindent
$W_{ijk}\,^{l}$ is the {\bf contact projective Weyl tensor}. $C_{ijk}$ should be regarded as an analogue of the Cotton tensor in conformal geometry, and will be called the {\bf contact projective Cotton tensor}. Direct computations verify the following identities.
\begin{align}
&\label{skewqskewp}Q_{[ij]} = -2P_{[ij]} = -\tfrac{2}{2n-1}R_{[ij]},\\
&\label{qrelatedtortraces}2(1-n)Q_{ij} + Q_{ji} = 2R_{ij} + R_{p}\,^{p}\,_{ij},& &Q_{ij} = 2R_{ij} - 4nP_{ij},
\end{align}
$R_{p}\,^{p} = 0$ and \eqref{skewqskewp} show that $P_{p}\,^{p} = 0 = Q_{p}\,^{p}$. When the contact torsion vanishes, \eqref{firstbianchicontracted} and \eqref{skewfirstbianchicontracted} show that $P_{ij} = \frac{1}{2n}R_{ij}$ and $Q_{ij} = 0$, and \eqref{cdefined} gives
\begin{equation}\label{ctfcotton}
C_{ijk} = R_{0ijk} - 2\nabla_{i}P_{jk} + \tfrac{2}{2n-1}(\omega_{ij}\nabla_{p}P_{k}\,^{p} + \omega_{ik}\nabla_{p}P_{j}\,^{p}).
\end{equation}
\noindent
The following identities are verified by direct computation using the definitions and perhaps also the Bianchi identities.
\begin{align}
&\label{wtraces}W_{p}\,^{p}\,_{ij} = 0, \quad &W_{ijp}\,^{p} = 0&, \quad &W_{ipj}\,^{p} = - \tfrac{1}{2}Q_{ij},\\
&\label{ctraces}C_{ip}\,^{p} = 0, \quad &C_{p}\,^{p}\,_{k} = 0&, \quad &C_{i[jk]} = 0.
\end{align}

Straightforward computation shows that the curvature tensors, $\tilde{R}$ and $R$, of connections, $\tilde{\nabla}$ and $\nabla$, with difference tensor, $\Lambda$, are related by
\begin{multline}\label{curvaturedifference}
\tilde{R}(X, Y)Z - R(X, Y)Z = \\
 (\nabla_{X}\Lambda)(Y, Z) - (\nabla_{Y}\Lambda)(X, Z) + \Lambda(X, \Lambda(Y, Z)) - \Lambda(Y, \Lambda(X, Z)) + \Lambda(\tau(X, Y), Z),
\end{multline}
where $\tau$ is the torsion tensor of $\nabla$. Write $\gamma_{ij} = \nabla_{i}\gamma_{j} - \gamma_{i}\gamma_{j} +\tfrac{1}{2}\gamma_{0}\omega_{ij}$. With the set-up as in Section \ref{transformationrules}, using \eqref{lambdatransform} and \eqref{lambda0ip} in \eqref{curvaturedifference} gives
\begin{align}\notag
\tilde{R}_{ijk}\,^{l} - R_{ijk}\,^{l} = 2\delta_{[j}\,^{l}\gamma_{i]k} + 2\gamma_{[i}\,^{l}\omega_{j]k} - 2\omega_{ij}\gamma_{k}\,^{l} - 2\omega_{ij}\tau_{qk}\,^{l}\gamma^{q} + \gamma_{k}\tau_{ij}\,^{l} + \gamma^{l}\tau_{ijk}.
\end{align}
Sometimes involved direct computations using this and \eqref{firstbianchicontracted}-\eqref{secondbianchizeroriccitrace} prove:
\begin{lemma}
Under change of scale there hold the following transformation rules:
\begin{align}
\label{riccitensortransform}
&\tilde{R}_{ij} - R_{ij} = 2n\gamma_{ij} + \gamma^{p}\tau_{pij},\\
\label{otherriccitransform}
&\tilde{R}_{p}\,^{p}\,_{ij} - R_{p}\,^{p}\,_{ij} = -4n\gamma_{ij} + 2\gamma^{p}\tau_{ijp} + 4(1-n)\tau_{pij},\\
\label{skewriccitransform}
&2\tilde{R}_{[ij]} - 2R_{[ij]} = (2n-1)\gamma^{p}\tau_{ijp}.&\\
&\label{ptransform} \tilde{P}_{ij} - P_{ij} = \gamma_{ij},\\
&\label{qtransform} \tilde{Q}_{ij} - Q_{ij} = 2\gamma^{p}\tau_{pij},\\
&\label{weyltransform} \tilde{W}_{ijk}\,^{l} - W_{ijk}\,^{l} = \gamma_{k}\tau_{ij}\,^{l} + \gamma^{l}\tau_{ijk},\\
&\label{ctransform}\tilde{C}_{ijk} - C_{ijk} = 2\gamma^{p}(W_{pijk} + \tau_{pij}\gamma_{k} + \tau_{pik}\gamma_{j}) + Q_{ij}\gamma_{k} + Q_{ik}\gamma_{j}.
\end{align}
In particular, if $\tau_{ij}\,^{k} = 0$, then $\tilde{W}_{ijk}\,^{l} = W_{ijk}\,^{l}$, so $W_{ijk}\,^{l}$ is an invariant of the contact projective structure.
\end{lemma}
\noindent
As is explained in Remark \ref{transformremark}, the proof of Lemma \ref{psubmodule} implies \eqref{ptransform}-\eqref{ctransform}.

\begin{lemma}\label{threeweylvanish}
In three dimensions $W_{ijk}\,^{l} = 0$ identically.
\end{lemma}

\begin{proof}
Since $\tau_{ij}\,^{k} = 0$, \eqref{bianchi1} and \eqref{wdefined} imply $W_{[ijk]l} = 0$. The Weyl tensor has the symmetries $W_{[ij]kl} = W_{ijkl}$, $W_{ij(kl)} = W_{ijkl}$, $W_{[ijk]l} = 0$, and it is completely trace free, so lies in an irreducible representation of $Sp(n-1, \rea)$, corresponding to the Young diagram determined by the partition $(3, 1)$, which, by the criterion of Weyl, (\cite{Weyl-Classical}, Section $6.3$, p. $175$), is non-trivial if and only if $2n-1 > 3$.
\end{proof}

\noindent
Proposition \ref{threectf}, Lemma \ref{threeweylvanish}, and \eqref{ctransform} show that $C_{ijk}$ is invariant when $2n -1 = 3$. 
\begin{lemma}\label{cottonweyllemma}
If the contact torsion vanishes,
\begin{align}
&\label{cottonweyl}\nabla_{p}W^{p}\,_{ijk} + \tfrac{1}{2n-1}(\nabla_{p}W_{ijk}\,^{p} - \nabla_{p}W_{kij}\,^{p}) = (2 - n)C_{ijk},&
\end{align}
so that if $2n -1 \geq 5$, then $W_{ijk}\,^{l} =0$ implies $C_{ijk} =0$. In dimension $2n-1 = 3$, $C_{ijk}$ is completely symmetric, $C_{(ijk)} = C_{ijk}$.
\end{lemma}

\begin{proof}
The first Bianchi identity implies $\nabla_{p}W_{[ijk]l} = 0$. Contracting on $p$ and $i$ gives $\nabla_{p}W_{ijk}\,^{p} = 2\nabla_{p}W^{p}\,_{[ji]k}$. Using $R_{0\alpha}= 0 = R_{\alpha0}$, \eqref{firstbianchicontractedzero}, $R_{ij} = 2nP_{ij}$, and calculating $\nabla_{p}W_{ijk}\,^{p}$ directly gives the first equality in 
\begin{align}
&\label{tracednablaweyl}\nabla_{p}W_{ijk}\,^{p} = 2(2n-1)\nabla_{[j}P_{i]k} + 2\omega_{k[i}\nabla^{p}P_{j]p} - 2\omega_{ij}\nabla^{p}P_{kp} = (2n-1)C_{[ij]k}.&
\end{align}
Since $\tau_{ij}\,^{k} = 0$, the first Bianchi identity shows $R_{0[ij]k} = 0$. By \eqref{firstbianchicontractedzero}, skewing \eqref{ctfcotton} in $ij$ gives the second equality in \eqref{tracednablaweyl}. By \eqref{wdefined}, \eqref{cdefined}, and \eqref{skewtracesecondbianchi}, 
\begin{equation}\label{cotton1}
\nabla_{p}W^{p}\,_{ijk} - (2-n)C_{ijk} =2\nabla_{[i}P_{j]k} - 2\nabla_{[k}P_{i]j}  + \tfrac{6}{2n-1}\omega_{i(j}\nabla^{p}P_{k)p}\,^{p}.
\end{equation}
Solving the first equality in \eqref{tracednablaweyl} for $2\nabla_{[i}P_{j]k}$ gives
\begin{equation*}
2\nabla_{[i}P_{j]k} - 2\nabla_{[k}P_{i]j} + \tfrac{6}{2n-1}\omega_{i(j}\nabla^{p}P_{k)p}  = \tfrac{1}{2n-1}(\nabla_{p}W_{kij}\,^{p} - \nabla_{p}W_{ijk}\,^{p}).
\end{equation*}
Substituting this into \eqref{cotton1} gives \eqref{cottonweyl}. By \eqref{cottonweyl}, when $2n-1 \geq 5$ the vanishing of $W_{ijk}\,^{l}$ implies the vanishing of $C_{ijk}$. \eqref{ctfcotton} shows that if $\tau_{ij}\,^{k} = 0$ then $C_{i(jk)} = C_{ijk}$. \eqref{tracednablaweyl} shows that if $W_{ijk}\,^{l} = 0$ then $C_{[ij]k} = 0$. When $2n-1 = 3$, both $\tau_{ij}\,^{k}$ and $W_{ijk}\,^{l}$ vanish, so $C_{ijk}$ is completely symmetric. 
\end{proof}

\section{Ambient Connection}\label{ambient}
\subsection{Thomas's Ambient Construction for Projective Structures}\label{thomasambient}
In this section the basic results concerning projective structures are reviewed, following the approach of T. Y. Thomas, \cite{Thomas}. Some version of this material can be found in various modern sources, for instance \cite{Bailey-Eastwood-Gover}, \cite{Gover}, \cite{Graham}, or \cite{Kobayashi-Nagano}.  

\subsubsection{Projective Structures}
The bundle of frames, $\mathsf{F}$, in the canonical bundle, $\wedge^{n}(T^{\ast}M)$, of the smooth $n$-dimensional manifold, $M$, is the $\rea^{\times}$ principal bundle of smooth, non-vanishing sections of the canonical bundle. When $n = 2l$ let the $\rea^{\times}$ principal bundle, $\dens$, be the unique $\frac{1}{n+1}$-root of $\mathsf{F}$. When $n = 2l-1$, assume $M$ is orientable with a fixed orientation, so that the group of $\mathsf{F}$ is reduced to $\rea^{>0}$, and let $\dens$ be a choice of $\frac{1}{n+1}$-root of $\mathsf{F}$. Denote by $\emf[\lambda]$ the line bundle associated to $\dens$ by the representation, $r\cdot s = r^{-\lambda}s$, of $\rea^{\times}$ on $\rea$, so that $\emf[-1]$ is a $\frac{1}{n+1}$-root of $\Lambda^{n}(T^{\ast}M)$. Sections of $\emf[\lambda]$ are in canonical bijection with functions $\dens \to \rea$ homogeneous of degree $\lambda$. The model for $\dens$ is the defining bundle $\standrepz \to \proj(\standrep)$, where $(\standrep, \svol)$ is an $n+1$ dimensional real vector space with volume form $\svol$. This defining bundle is the bundle of frames in the tautological line bundle over $\proj(\standrep)$. 

Fix a choice of $\rho:\dens \to M$. A tautological $n$-form, $\alpha$, is defined on $\dens$, by 
\begin{align*}
&\alpha_{s}(X_{1}, \dots, X_{n}) = s^{n+1}_{\rho(s)}(\rho_{\ast}(X_{1}), \dots, \rho_{\ast}(X_{n}))& & \text{for}\, X_{1}, \dots, X_{n} \in T_{s}\dens.&
\end{align*}
A canonical volume form, $\svol$, is defined by $\svol = d\alpha$. A canonical {\bf Euler vector field}, $\eul$, the infinitesimal generator of the fiber dilations, $\delta_{r}$, on $\dens$, is defined by $i(\eul)\svol = (n+1)\alpha$. The choice of section, $s$, induces, on the principal $\rea^{\times}$-bundle $\rho:\dens \to M$, a unique connection, $\phi$, such that $s$ is a parallel section and $\phi(\eul) = 1$.  The connection, $\phi$, determines a horizontal lift, $\hat{X}$, of each vector field, $X$, on $M$.
\begin{lemma}[Weyl, \cite{Weyl}]\label{weyllemma}
Two torsion-free affine connections have the same unparameterized geodesics if and only if there exists a one-form, $\gamma$, so that their difference tensor has the form $\gamma_{(i}\delta_{j)}\,^{k}$.
\end{lemma}
\begin{lemma}[Thomas, \cite{Thomas}]\label{thomaslemma}
Suppose given a projective structure on $M$ represented by the equivalence class of torsion-free affine connections, $[\nabla]$. For each choice of a volume form, $\mu$, on $M$ there exists a unique $\nabla \in [\nabla]$ making $\mu$ parallel.
\end{lemma}
\noindent
A torsion-free affine connection, $\hat{\nabla}$, on $\dens$, and a choice of section, $s$, determine on $M$ a torsion-free affine connection, $\nabla$, defined by $\nabla_{X}Y = \rho_{\ast}(\Hat{\nabla}_{\hat{X}}\hat{Y})$. If the condition $\hat{\nabla}\eul = \delta$ is imposed, Lemma \ref{weyllemma} may be used to show that the connection, $\tilde{\nabla}$, determined by $\tilde{s} = fs$ has the same unparameterized geodesics as does $\nabla$.
This associates a projective structure on $M$ to each torsion-free affine connection, $\hat{\nabla}$, on $\dens$, satisfying $\hat{\nabla}\eul = \delta$. It may be checked that if $\hat{\nabla}\svol = 0$ then $\nabla$ is the unique representative of $[\nabla]$ given by Lemma \ref{thomaslemma} and making $\mu= s^{\ast}(\alpha)$ parallel. Two torsion-free connections, $\hat{\nabla}$ and $\hat{\nabla}^{\prime}$, on $\dens$, satisying $\hat{\nabla}\eul = \delta = \hat{\nabla}^{\prime}\eul$ and making parallel $\svol$ determine on $M$ the same projective structure. The only freedom is in the vertical part of $\hat{\nabla}_{\hat{X}}\hat{Y}$, and requiring that $\hat{\nabla}$ be Ricci flat eliminates this freedom.

\begin{theorem}[Thomas, \cite{Thomas}]\label{thomastheorem}
Suppose fixed a choice of $\rho:\dens\to M$. There is a functor associating to each projective structure on $M$ represented by the equivalence class of torsion-free affine connections, $[\nabla]$, a unique torsion-free affine connection, $\hat{\nabla}$, (the \textbf{Thomas ambient connection}) on $\dens$ satisfying the following conditions.
\setcounter{condition}{0}
\begin{condition}\label{tcon1}
$\hat{\nabla}\eul$ is the fundamental $\binom{1}{1}$-tensor on $\dens$. 
\end{condition}
\begin{condition}\label{tcon2}
 $\hat{\nabla}\svol = 0$.
\end{condition}
\begin{condition}\label{tcon3}
The Ricci tensor of $\hat{\nabla}$ vanishes.
\end{condition}
\begin{condition}\label{tcon4}
For each section, $s$, of $\rho:\dens\to M$ the connection, $\nabla$, defined on $M$ by $\nabla_{X}Y = \rho_{\ast}(\hat{\nabla}_{\hat{X}}\hat{Y})$, represents the given projective structure. Moreover, $\nabla$ is the unique representative, given by Lemma \ref{thomaslemma}, making the volume form $s^{\ast}(\alpha)= s^{n+1}$ parallel.
\end{condition}
\end{theorem}
\noindent
The preceeding discussion shows already how to prove existence of $\hat{\nabla}$. Given $s$, let $\nabla \in [\nabla]$ be the unique representative of the given projective structure $\nabla \in [\nabla]$ making $\mu = s^{\ast}(\alpha)$ parallel. For any symmetric tensor, $P_{ij}$, conditions \eqref{tcon1}-\eqref{tcon3} of Theorem \ref{thomastheorem} are satisfied by the connection, $\hat{\nabla}$, defined by requiring it to be torsion-free and to satisfy
\begin{equation}\label{definethomasambient}
\hat{\nabla}_{V}\eul = V, \qquad\qquad \hat{\nabla}_{\hat{X}}\hat{Y} = \widehat{\nabla_{X}Y} + P(X, Y)\eul.
\end{equation}
Straightforward computation of the curvature of $\hat{\nabla}$ shows that requiring $\hat{\nabla}$ to be Ricci-flat determines $P_{ij}$ uniquely as $P_{ij} = \frac{1}{n-1}R_{ij}$, where $R_{ij}$ is the Ricci tensor of $\nabla$. Because the curvature tensor of $\hat{\nabla}$ turns out to be horizontal, its components may be regarded as tensors on $M$. The possibly non-vanishing components are the \textbf{projective Weyl tensor}, $B_{ijk}\,^{l} = R_{ijk}\,^{l} + P_{jk}\delta_{i}\,^{l} - P_{ik}\delta_{j}\,^{l}$, and the \textbf{projective Cotton tensor}, $2\nabla_{[i}P_{j]k}$. $B_{ijk}\,^{l}$ is independent of the choice of representative of $[\nabla]$, and the contracted second Bianchi identity for $\hat{\nabla}$ shows that if $n > 2$ the vanishing of $B_{ijk}\,^{l}$ implies the vanishing of $\nabla_{[i}P_{j]k}$ and consequently of the curvature of $\hat{\nabla}$. If $n = 2$, it turns out that $B_{ijk}\,^{l}$ vanishes automatically, and the vanishing of $\nabla_{[i}P_{j]k}$, implies the ambient connection is flat. 

Let $G$ be the group of linear transformations of $\standrep$ preserving $\svol$; let $P \subset G$ be the subgroup stabilizing the span of a fixed vector. Next there is sketched the construction of a unique $(\g, P)$ Cartan connection associated to each projective structure equipped with a choice of $\dens$. The terminology and notation regarding Cartan connections are the same as those in \cite{Cap-Slovak-Soucek} or \cite{Slovak}. The \textbf{tractor bundle}, $\tractor$, is the quotient of $T\dens$ by an $\rea^{\times}$ action, $P_{r}= r^{-1}\delta_{r^{-1}}^{\ast}$, covering the principal action on $\form$ and preserving the homogeneity $-1$ vector fields, $\vect_{-1}(\dens)$ on $\dens$. The homogeneity $n+1$ volume form, $\svol$, descends to a fiberwise volume form on $\tractor$, and the ambient connection, $\hat{\nabla}$, induces a covariant differentiation, $\tnabla$, on $\tractor$, the {\bf tractor connection}, defined as follows. If $t$ is a section of $\tractor$ represented by $Z\in \vect_{-1}(\dens)$, then $\tnabla_{X}t$ is the image in $\tractor$ of $\hat{\nabla}_{\hat{X}}Z\in \vect_{-1}(\dens)$. Because $\hat{\nabla}_{\eul}Z = 0$, this is independent of the choice of $s$. A $P$ principal bundle, $\pi:\adaptedframe \to M$, is defined by letting the fiber over $x \in M$ comprise all volume preserving linear isomorphisms $u:\standrep \to \tractor_{x}$ mapping the span of $e_{\infty}$ into the image of the span of $\eul$, and $\tractor$ is recovered as the associated bundle, $\adaptedframe\times_{P}\standrep$. Theorem 2.7 of \cite{Cap-Gover} shows that the tractor connection determines on $\G$ a Cartan connection, $\eta$. Conversely, any torsion-free, $(\g, P)$ Cartan connection, $\eta$, on $\adaptedframe \to M$, determines a tractor connection as the induced covariant differentiation on the associated bundle $\tractor = \adaptedframe\times_{P}\standrep$. Any $\eta$ determines a development of paths in $M$ onto $G/P = \proj(\standrep)$, (see \cite{Sharpe}), and $\eta$ induces on $M$ the projective structure comprising those paths which develop onto straight lines in $\proj(\standrep)$. Evidently gauge equivalent Cartan connections induce the same projective structure. The quotient $\adaptedframe/\tilde{P}$ recovers $\dens$ and the associated bundle $\adaptedframe \times_{\tilde{P}}\standrep \to \dens$ recovers $T\dens$; $\eta$ induces a covariant differentiation, $\etanabla$, on this associated bundle. Every torsion-free $(\g, P)$ Cartan connection is gauge equivalent to a unique $\eta$ for which the induced $\etanabla$ on $T\form$ satisifes $\etanabla \eul = \delta$. Each section $s:M \to \form$ determines from $\eta$ an affine connection on $M$ defined as in condition \conref{tcon4} of Theorem \ref{thomastheorem}. If $\etanabla \eul = \delta$, then the affine connections on $M$ determined by different choices of $s$ will be projectively related, so determine a projective structure on $M$, the paths of which develop onto lines in $\proj(\standrep)$. Among the isomorphism classes of Cartan connections inducing a given projective structure that one built from the ambient connection is the unique isomorphism class of torsion-free, normal Cartan connections. An orientation-preserving diffeomorphism of $M$ lifts to a bundle automorphism of $\form$ which commutes with $P_{r}$ and so descends to $\tractor$, and this shows that the orientation-preserving diffeomorphisms of $M$ act naturally on $\G$ and all of its associated bundles. This constructs a bijective functor between the category of projective structures and category of torsion-free, normal $(\g, P)$ Cartan connections. In the odd-dimensional case, to obtain structure independent of the choice of $\dens$, it is necessary to consider the projectivized tractor bundle, $\proj(\tractor)$, and the associated bundle of filtered projective frames, $\pframe$, the fiber over $x \in M$ of $\pframe$ comprising all projective linear isomorphisms $u:\proj(\standrep) \to \proj(\tractor_{x})$ mapping the span of $e_{\infty}$ to $\proj(\tractor^{1}_{x})$. These bundles do not depend on the choice of $\form$. Some fussing shows that the Cartan connection, $\eta$, descends to give a Cartan connection on $\pframe$ which does not depend on the choice of $\form$.

\subsection{Ambient Connection on the Symplectification}\label{symplectification}
This section gives the proof of Theorem B. On a co-oriented contact manifold, $(M, H)$, a choice of a square-root, $\form$, of the bundle of positive contact one-forms, regarded as a $\rea^{\times}$ principal bundle over $M$, will be called an {\bf ambient manifold}.

\subsubsection{The Standard Contact Structure on Odd-dimensional Projective Space}
Let $(\standrep, \Omega)$ be the standard representation of $G$ as in Section \ref{flatmodel}. The quotient of $\standrepz$ by the action on $\standrep$ of the dilations, $\delta_{r}(v) = rv$, determines an $\rea^{\times}$-principal bundle, $\rho:\standrepz \to \proj(\standrep)$, so that the dilations act as vertical principal bundle automorphisms. Let $\eul$ denote the Euler vector field, the infinitesimal generator of the fiber dilations. Let $W^{\perp}$ denote the $\Omega$-skew complement of a subspace $W \subset \standrep$. For $L \in \proj(\standrep)$ choose any non-zero $v \in L$ and set $H_{L} = \rho_{\ast}(v)(L^{\perp})$. Because any other vector spanning $L$ will have the form $\delta_{r}(v)$ for some non-zero $r$, $H_{L}$ does not depend on the choice of $v$. As $\ker \rho_{\ast}(v)$ is exactly $L$, which is contained in $L^{\perp}$, the subbundle $H$ has constant rank $2n-2$. As $G$ preserves skew complements, the subbundle $H$ is invariant under the action $G$ on $\proj(\standrep)$. Define a one-form $\alpha$ on $\standrepz$ by $\alpha = \tfrac{1}{2}i(\eul)\Omega$. A section $s:\proj(\standrep) \to \standrepz$ determines a one-form, $\theta$, by $\theta = s^{\ast}(\alpha) = s^{2}$, for which $\ker \theta = H$. Using $\theta = s^{\ast}(\alpha)$ it is easy to verify that $\theta \wedge (d\theta)^{n-1} \neq 0$, so that $\theta$ is a contact one-form, and $H$ is a contact distribution. When the section, $s$, is replaced by $\tilde{s} = fs$, $f \neq 0$, then  $\theta$ rescales to $f^{2}\theta$. This shows that the principal $\rea^{\times}$-bundle, $\rho:\standrepz \to \proj(\standrep)$, is naturally identified with a square-root of the bundle of positive contact one-forms on $\proj(\standrep)$; the flat Euclidean connection, $\hat{\nabla}$, on $\standrepz$ satisfies $\hat{\nabla}\alpha = \tfrac{1}{2}\Omega$ and $\hat{\nabla}\Omega = 0$, and is the model for the ambient connection. 

\subsubsection{Construction of the Ambient Connection}\label{ambientconstruction}
The map $\theta \to \theta \wedge d\theta^{n-1}$ determines a reduction to $\rea^{>0}$ of $\candens$, and by \eqref{volumetransform} this map identifies $\contactb^{n}$ with $\candens$ as principal $\rea^{>0}$ bundles, so that a square-root, $\form$, of $\contactb$ may be regarded simultaneously as a $\frac{1}{2n}$-root of $\candens$. The well known construction of the tautological one-form on $\contactb$ generalizes to $\form$ as follows. For $p \in \form$ view $p^{2}$ as a one-form on $M$ and define the horizontal one-form, $\alpha$, on $\form$, by $\alpha_{p}(X) = p^{2}_{\rho(p)}(\rho_{\ast}(p)(X))$, for $X\in T_{p}\form$. Each (local) section $s: M \to \form$ determines on $M$ a contact one-form $\theta = s^{\ast}(\alpha) = s^{2}$, and a fiber coordinate, $t \neq 0$, defined by $p = ts(\rho(p))$ for any $p \in \form$. Set $\Omega = d\alpha$. By definition,
\begin{align}\label{Omega}
&\alpha = t^{2}\rho^{\ast}(\theta),& &\Omega = t^{2}(2d\log{t} \wedge \rho^{\ast}(\theta) + \rho^{\ast}(d\theta)).&
\end{align}
Computing using \eqref{Omega} and $d\theta^{n} = 0$ shows that 
\begin{equation}\Omega^{n} = 2n t^{2n} d\log{t} \wedge \rho^{\ast}(\theta \wedge d\theta^{n-1}).
\end{equation}
As $dt$ is non-vanishing when restricted to the vertical subbundle of $T\form$ and $\theta \wedge (d\theta)^{n-1}$ is a volume form on $M$, this shows that $\Omega^{n} \neq 0$, so that $\Omega$ is a symplectic structure on $\form$. The \textbf{Euler vector field}, $\eul$, is the infinitesimal generator of the fiber dilations in $\form$. \eqref{Omega} shows that in coordinates $\eul = t\frac{\partial}{\partial t}$. A tensor, $S$, on $\form$ is homogeneous of degree $k$ if $\lie_{\eul}S = kS$. By non-degeneracy of $\Omega$ and definition of $\alpha$, $\lie_{\eul}\alpha = i(\eul)d\alpha = 2\alpha$, so $\alpha$ (and hence also $\Omega$) is homogeneous of degree $2$. The choice of section, $s$, induces, on the $\rea^{\times}$ principal bundle $\rho:\form \to M$, a unique connection, $\phi$, such that $s$ is a parallel section and $\phi(\eul) = 1$. Because the structure group is abelian, $\phi$ is simply a one-form on the base manifold, $M$, given in coordinates by $d\log{t}$. The section $\tilde{s}$ determines the connection $\tilde{\phi}$, related to $\phi$ by $\tilde{\phi} = \phi - \rho^{\ast}(\gamma)$. The connection, $\phi$, determines a horizontal lift, $\hat{X}$, of a vector field, $X$, on $M$; this $\hat{X}$ is the unique $\phi$-horizontal vector field on $\form$ such that $\rho_{\ast}(\hat{X}) = X$. By definition $\hat{X}$ is homogeneous of degree $0$. The horizontal lift determined by the connection associated to $\tilde{s}$ is given by $\hat{X} + \gamma(X)\eul$. A co-oriented contactomorphism of $M$ is covered by a unique bundle automorphism of $\form$, which is easily checked to be a symplectomorphism.

In this section the following indexing conventions will be employed. For $\theta = s^{\ast}(\theta)$, let $\{\theta^{\alpha}\}$ be a $\theta$-adapted coframe with dual frame $\{E_{\alpha}\}$. Define a frame, $\{F_{I}\}$, on the total space of $\form$ by setting $F_{\infty} = \eul$ and $F_{\alpha} = \hat{E_{\alpha}}$, and let $\phi^{I}$ denote the dual coframe. Thus $\phi^{\infty} = \phi$, $\phi^{\alpha} = \rho^{\ast}(\theta^{\alpha})$, and $\alpha = t^{2}\phi^{0}$. On $\form$, raise and lower indices using $\Omega$, according to the convention $\Omega^{IK}\Omega_{KJ} = -\delta_{I}\,^{J}$. Note that $F_{I}$ and $\phi^{I}$ are homogeneous of degree $0$. In general, an affine connection, $\hat{\nabla}$, on $\form$, has torsion, so there will be two distinct traces of its curvature tensor, $\hat{R}_{IJK}\,^{L}$, namely
$\hat{R}_{IJ} = \hat{R}_{IPJ}\,^{P}$ and $\hat{S}_{IJ} = \hat{R}_{Q}\,^{Q}\,_{IJ}$. Contracting the first Bianchi identity shows that when the torsion of $\hat{\nabla}$ vanishes these two traces coincide up to a factor of $-2$. 

\begin{proof}[Proof of Theorem B]
Fix a section $s:M \to \form$, let $\theta = s^{\ast}(\alpha)$, and let $\nabla$ be the connection associated to $\theta$ by Theorem A. Define $\hat{\nabla}$ by specifying its action on the frame, $\{F_{I}\}$. For conditions \conref{con1} and \conref{con2} to hold, $\hat{\nabla}$ must be defined to satisfy $\hat{\nabla}_{\hat{E}_{\alpha}}\eul = \hat{E}_{\alpha} =\hat{\nabla}_{\eul}\hat{E}_{\alpha}$, and $\hat{\nabla}_{\eul}\eul = \eul$. Then there can be written
\begin{align*}
&\hat{\nabla}_{\hat{E}_{\alpha}}\hat{E}_{\beta} = \widehat{\nabla_{E_{\alpha}}E_{\beta}} + \widehat{\A(E_{\alpha}, E_{\beta})} + P(E_{\alpha}, E_{\beta})\eul,&\end{align*}
where $P$ is an arbitrary $\binom{0}{2}$-tensor on $M$ and $\A$ is an arbitrary $\binom{1}{2}$-tensor on $M$. The proof will show that $P$ and $\A$ are defined as follows. $P_{ij}$ and $Q_{ij}$ are defined as in \eqref{pdefined} and \eqref{qdefined}, and the remaining components are
\begin{align}
& Q_{\alpha 0} = 0 = Q_{0\alpha},\\
&\label{p0idefined}P_{0i} = \tfrac{2}{2n-1}\nabla^{p}P_{ip} + \tfrac{1}{2n-1}\nabla^{p}Q_{ip},\\
&\label{pi0defined}P_{i0} = \tfrac{2}{2n-1}\nabla^{p}P_{ip} - \tfrac{1}{2(n-1)(2n-1)}\nabla^{p}Q_{ip} - \tfrac{1}{n-1}\tau_{i}\,^{pq}P_{qp},\\
&\label{p00defined}P_{00} = \tfrac{1}{n-1}\nabla^{p}P_{0p} - \tfrac{1}{n-1}(2P^{pq} + Q^{pq})P_{qp},\\
&\label{oabcdefined}\A_{\alpha\beta}\,^{\gamma} = 2(\delta_{\alpha}\,^{0}P_{\beta}\,^{\gamma} + \delta_{\beta}\,^{0}P_{\alpha}\,^{\gamma} - \delta_{\alpha}\,^{0}\delta_{\beta}\,^{0}P_{0}\,^{\gamma}) +\delta_{\alpha}\,^{0}Q_{\beta}\,^{\gamma} -\tfrac{1}{2}\delta_{0}\,^{\gamma}\omega_{\alpha\beta} .
\end{align}
where $P_{\alpha}\,^{\gamma}$ is defined by $P_{\alpha}\,^{\gamma} = P_{\alpha}\,^{q}\delta_{q}\,^{\gamma}$, and similarly for $Q_{\alpha}\,^{\gamma}$.

So that uniqueness may be shown in passing, suppose now that $P$ and $\A$ are arbitrary. It will be shown that conditions \conref{con3}-\conref{con6} determine $P$ and $\A$ uniquely as defined above. Computing $\hat{\nabla}\Omega$ using \eqref{Omega} shows that \conref{con3} holds if and only if $\A_{\alpha[jk]} = 0$, $2P_{\alpha i} = \A_{\alpha 0i}$, and $\A_{\alpha\beta}\,^{0} = -\tfrac{1}{2}\omega_{\alpha\beta}$, and these are satisfied by \eqref{oabcdefined}. Using Lemma \ref{differencetensor} it is easily verified that a connection on $T\form$ satisfying conditions \conref{con1} and \conref{con2} satisfies \conref{con6} if it satisfies \conref{con6} for a single choice of section $s:M \to \form$. Define a connection, $\bar{\nabla}$, on $M$, by $\bar{\nabla}_{X}Y = \rho_{\ast}(\hat{\nabla}_{\hat{X}}\hat{Y})$. Then the difference tensor of $\bar{\nabla}$ and $\nabla$ is $\A$, and by Lemma \ref{differencetensor} this implies  $\A_{(ij)}\,^{0} = 0$ and that there is $\gamma_{i}$ such that $\A_{(ij)k} = 2\gamma_{(i}\omega_{j)k}$. Tracing this and using $\A_{i[jk]} = 0$ gives $\A_{pi}\,^{p} = 2(2n-1)\gamma_{i}$. Skewing $\hat{\nabla}\Omega = 0$ gives $\hat{\tau}_{[IJK]} = 0$, and computing directly using \conref{con1} and \conref{con2} gives $\hat{R}_{IJ\infty}\,^{K} = \hat{\tau}_{IJ}\,^{K}$. Imposing $\hat{R}_{I\infty} = 0$ gives $0 = \hat{R}_{IQ\infty}\,^{Q} = \hat{\tau}_{IQ}\,^{Q}$, and with $\hat{\tau}_{[IJK]} = 0$, this implies $\hat{\tau}$ is completely trace free. Using \conref{con2} gives $0 = \hat{\tau}_{iQ}\,^{Q} = \hat{\tau}_{iq}\,^{q}$, which implies $0 = \tau_{iq}\,^{q} + 2\A_{[iq]}\,^{q} = 2\A_{[iq]}\,^{q}$. Since $\A_{iq}\,^{q} = 0$, this gives $\A_{qi}\,^{q} = 0$, so that $\gamma_{i} = 0$, and hence $\A_{(ij)k} = 0$. Now $\A_{ijk}$ is skew in $ij$ and symmetric in $jk$, and so $\A_{ijk} = 0$. The only components of $\A$ not determined in terms of $P$ are $\A_{0ij}$. Define $Q_{ij} = \A_{0ij} - 2P_{ij}$ and $Q_{\alpha 0} = 0 = Q_{0\alpha}$, and note that $\A$ must satisfy \eqref{oabcdefined}. Note also that $0 = \A_{0[ij]}$ implies $2P_{[ij]} = -Q_{[ij]}$. The proof will be completed if $P_{\alpha\beta}$ and $Q_{ij}$ can be determined uniquely by the remaining conditions.

Direct computation using conditions \conref{con1}, \conref{con2}, and \conref{con3} of Theorem B, $\A_{ij}\,^{k} = 0$, and the definition of the curvature tensor gives
\begin{align}
&\label{wijkl2}\hat{R}_{ijk}\,^{l} =  R_{ijk}\,^{l} + 2\delta_{[i}\,^{l}P_{j]k} +2\omega_{k[j}P_{i]}\,^{l} + 2\omega_{ij}P_{k}\,^{l} + \omega_{ij}Q_{k}\,^{l} 
= W_{ijk}\,^{l},\\
&\label{cijk2}\hat{R}_{0ij}\,^{k} = R_{0ij}\,^{k} - P_{0j}\delta_{i}\,^{k} - \omega_{ij}P_{0}\,^{k} - 2\nabla_{i}P_{j}\,^{k} - \nabla_{i}Q_{j}\,^{k} = C_{ij}\,^{k},\\
&\label{uijkdefined}\hat{R}_{ijk}\,^{\infty} = 2\nabla_{[i}P_{j]k} + \tau_{ij}\,^{p}P_{pk} + \omega_{k[j}P_{i]0} + \omega_{ij}P_{0k} = U_{ijk},\\
&\label{vijdefined}\hat{R}_{ij0}\,^{\infty} = 2\nabla_{[i}P_{j]0} + \tau_{ij}\,^{p}P_{p0} + \omega_{ij}P_{00} + 4P_{ip}P_{j}\,^{p} = 2V_{ij},\\
&\label{aijdefined} \hat{R}_{0ij}\,^{\infty} = 2\nabla_{[0}P_{i]j} - \tfrac{1}{2}\omega_{ij}P_{00} - 2P_{ip}P_{j}\,^{p} - P_{ip}Q_{j}\,^{p} = A_{ij},\\
&\label{bidefined}\hat{R}_{0i0}\,^{\infty} = 2\nabla_{[0}P_{i]0} + 4P_{0p}P_{i}\,^{p} = 2B_{i}.\\
& \label{IJinfK}\hat{R}_{IJ\infty}\,^{K} = \hat{\tau}_{IJ}\,^{K}.
\end{align}
Here \eqref{wijkl2}-\eqref{bidefined} define the tensors $W_{ijk}\,^{l}$, $C_{ijk}$, $A_{ij}$, $B_{i}$, $U_{ijk}$, and $V_{ij}$. Once it has been shown that $P_{\alpha\beta}$ and $Q_{ij}$ must be as in \eqref{pdefined}, \eqref{pi0defined}-\eqref{p00defined}, and \eqref{qdefined}, then $W_{ijk}\,^{l}$ and $C_{ijk}$ can be checked to equal the tensors defined in \eqref{wdefined} and \eqref{cdefined}. $\hat{\nabla}\Omega = 0$ implies $\hat{R}_{IJ(KL)} = \hat{R}_{IJKL}$. Using this and \eqref{IJinfK}, all the components of $\hat{R}_{IJKL}$ not listed are derivable from those listed. For instance, $\hat{R}_{ij0}\,^{k} = 2U_{ij}\,^{k}$ and $\hat{R}_{0i0}\,^{j} = 2A_{i}\,^{j}$.

Using \conref{con1} and \conref{con2} gives $\hat{R}_{I\infty J}\,^{K} = 0$, and with $\hat{\tau}_{[IJK]} = 0$, this gives $\hat{R}_{IJ} = \hat{\tau}^{\infty}_{IJ} + \hat{R}_{IqJ}\,^{q}$. Using conditions \conref{con1}-\conref{con2} it is easy to check that $\hat{R}_{i\infty} = 0 = \hat{R}_{\infty i}$ and $\hat{R}_{\infty 0}$. Direct computation shows that $\hat{R}_{0\infty} = Q_{p}\,^{p} - 2P_{p}\,^{p} = 2Q_{p}\,^{p}$, so that imposing $\hat{R}_{0\infty} = 0$ forces $Q_{p}\,^{p} = 0$, and hence also $P_{p}\,^{p} = 0$. Tracing \eqref{IJinfK} and using condition \conref{con2} implies $\hat{R}_{\infty\infty}= 0$, and using also $\hat{R}_{IJ(KL)} = \hat{R}_{IJKL}$ gives $\hat{S}_{\infty I} = \hat{S}_{I \infty} = 0$. Using \eqref{skewfirstbianchicontracted} or $\hat{R}_{IJ(KL)} = \hat{R}_{IJKL}$ shows that $\hat{S}_{[ij]} = 0$, so that requiring $\hat{R}_{ij} = 0$ and $\hat{S}_{ij} = 0$ gives the three linear equations $\hat{R}_{[ij]} = 0$, $\hat{R}_{(ij)} = 0$, and $\hat{S}_{ij} = 0$, in the three unknowns $P_{[ij]}$, $P_{(ij)}$, and $Q_{(ij)}$, and so uniquely determine these unknowns, provided the equations can be solved. Straightforward computations using \eqref{firstbianchicontracted} and \eqref{wijkl2} show that the unique solutions of this system of equations are $P_{ij}$ and $Q_{ij}$ as defined in \eqref{pdefined} and \eqref{qdefined}. Imposing $\hat{R}_{i0} = 0$ and computing using \eqref{cijk2} determines $P_{0i}$ uniquely as in \eqref{p0idefined}. Imposing $\hat{R}_{i0} = 0$ and computing using \eqref{uijkdefined} and \eqref{p0idefined} determines $P_{i0}$ uniquely as in \eqref{pi0defined}. Finally, imposing $\hat{R}_{00} = 0$ and computing using \eqref{aijdefined} determines $P_{00}$ uniquely as in \eqref{p00defined}.

There remains to verify the claims regarding the situation of vanishing contact torsion. Direct computation using $\A_{[ij]}\,^{\gamma} = -\frac{1}{2}\delta_{0}\,^{\gamma}\omega_{ij}$ gives
\begin{align}
&\label{hattauijk}\hat{\tau}_{ij}\,^{k} =  \tau_{ij}\,^{k},& &\hat{\tau}_{0i}\,^{k} = Q_{i}\,^{k},& &\hat{\tau}_{\alpha\beta}\,^{\infty} = 2P_{[\alpha\beta]},&
& \hat{\tau}_{ij}\,^{0} = 0,& &\hat{\tau}_{i0}\,^{0} = 0,&
\end{align}
When the contact torsion vanishes, $P_{ij} = \tfrac{1}{2n}R_{ij}$ is symmetric, and $Q_{ij} = 0$, and so, by \eqref{hattauijk}, if the contact torsion vanishes, the ambient connection, $\hat{\nabla}$, is torsion free. The other components of $\hat{S}$ are
\begin{align}
&\label{hatsi0}\hat{S}_{i0} = 2U_{p}\,^{p}\,_{i} =   -4\nabla^{p}P_{pi} + 2P_{i0} + 4(n-1)P_{0i} ,\\
&\label{hats00}\hat{S}_{00} = 4V_{p}\,^{p} = 4\nabla_{p}P^{p}\,_{0} + 4(n-1)P_{00} + 8P_{pq}P^{pq}.
\end{align}
When the contact torsion vanishes, \eqref{pi0defined} and \eqref{p0idefined} imply $P_{[i0]} = 0$ and using these and \eqref{p00defined} in \eqref{hatsi0} and \eqref{hats00} shows that $\hat{S}_{i0} = 0$ and $\hat{S}_{00} = 0$. Alternatively, applying conditions \conref{con1} and \conref{con2} shows that $\hat{\nabla}^{Q}\hat{\tau}_{QIJ}= \hat{\nabla}^{q}\hat{\tau}_{qIJ}$. Using this, contracting the first Bianchi identity, using $\hat{\tau}_{[IJK]} = 0$ and that $\hat{\tau}$ is trace free shows
\begin{align}
\label{skewtracefirstbianchirefined}
2\hat{R}_{IJ} + \hat{S}_{IJ} = -2\hat{\nabla}^{q}\hat{\tau}_{qIJ} - \hat{\tau}^{pq}\,_{I}\hat{\tau}_{pqJ}.
\end{align}
and the vanishing of $\hat{S}_{i0}$ and $\hat{S}_{00}$ then follows from the vanishing of $\hat{\tau}$ and the Ricci curvature of $\hat{\nabla}$.
\end{proof}

The following remarks are made for use in Section \ref{algebraic}. \eqref{uijkdefined}, \eqref{hatsi0} and \eqref{hats00} give immediately
\begin{align}
&\label{uipp}U_{ip}\,^{p}  = P_{[i0]},& 
&A_{p}\,^{p} = \tfrac{1}{2}\hat{R}_{00} = 0,&
&V_{p}\,^{p} = \tfrac{1}{4}\hat{S}_{00}.&
\end{align}
Using \eqref{skewtracesecondbianchia}, \eqref{wdefined}, and \eqref{qrelatedtortraces} gives $2(1-n)\nabla^{p}Q_{pi} + \nabla^{p}Q_{ip} = \tau^{pqr}R_{pqri} = \tau^{pqr}W_{pqri} -2P^{pq}\tau_{ipq}$. Using this, \eqref{p0idefined}, and \eqref{pi0defined}, $\hat{S}_{i0}$ may be computed from \eqref{hatsi0}. 
\begin{align*}
&(1-n)\hat{S}_{i0} = 2(1-n)(\nabla^{p}Q_{pi} + 2P_{[i0]}) = \tau^{pqr}(W_{pqri} + 2\omega_{ri}P_{pq}).&
\end{align*}
The quantity $\tau^{pqr}(W_{pqri} + 2\omega_{ri}P_{[pq]})$ is independent of the choice of scale. This can be verified by direct computation using \eqref{ptransform} and \eqref{weyltransform}. Alternatively, it follows by computing $\hat{S}$ with respect to frames on $\form$ determined by different choices of scale and using the vanishing of $\hat{S}_{ij}$, $\hat{S}_{\infty I}$, and $\hat{S}_{I \infty}$. Similarly, $\hat{S}_{00}= 8\nabla^{p}P_{[0p]} - 4P^{pq}Q_{pq}$ follows from \eqref{hats00}. If $\hat{S}_{i0}$ vanishes then $\hat{S}_{00}$ is independent of the choice of scale.

Also note for later the following. Let $\bar{\nabla}$ be defined by \conref{con6}. Using \eqref{Omega}, condition \conref{con3}, and \eqref{oabcdefined} gives 
\begin{align}
&\label{hatnablaalphaOmega}
\bar{\nabla}_{\gamma}\omega_{ij} = 0,&
&\bar{\nabla}_{\alpha}\theta_{\beta} =\tfrac{1}{2}\omega_{\alpha\beta},& 
&\hat{\nabla}\alpha = \tfrac{1}{2}\Omega. 
\end{align}

\begin{remark}
The formulas for the transformations of $P$ and $\A$ under a change of scale are derivable consequences of Theorem B, rather than necessary constituents of its proof. It is possible also to prove Theorem B by defining $\hat{\nabla}$ as above, simply defining $P$ and $\A$ by \eqref{pdefined}, \eqref{qdefined}, and \eqref{pi0defined}-\eqref{p00defined}. Verifying that the $\hat{\nabla}$ so defined does not depend on the choice of scale requires computing explicitly how $P$ and $\A$ transform under a change of scale. This is computationally intensive, requiring, for instance, the verification of the following transformation rules.
\begin{align}
&\label{pi0transform} f^{2}\tilde{P}_{i0} - P_{i0} = \nabla_{i}\gamma_{0} - 2\gamma_{i}\gamma_{0} + \gamma^{p}(2\nabla_{i}\gamma_{p} +4P_{ip}),\\
&\label{p0itransform}f^{2}\tilde{P}_{0i} - P_{0i} =  \nabla_{0}\gamma_{i} - 2\gamma_{0}\gamma_{i}  + \gamma^{p}(2\nabla_{i}\gamma_{p} + 4P_{pi} + Q_{pi} + 2\tau_{piq}\gamma^{q}),\\
&\label{p00transform} f^{4}\tilde{P}_{00} - P_{00} = \nabla_{0}\gamma_{0} - \gamma_{0}^{2} + \gamma^{p}(4\nabla_{0}\gamma_{p} + 4\gamma^{q}\nabla_{p}\gamma_{q} + 4P_{0p} + 2P_{p0})\\\notag &\qquad\qquad\qquad + \gamma^{p}(12\gamma^{q}P_{pq} + 2\gamma^{q}Q_{pq}).
\end{align}
\end{remark}

\begin{remark}
An isotropic submanifold, $N$, of the contact projective manifold, $(M, H)$, is {\bf totally geodesic} if every contact geodesic of $M$ tangent to $N$ at one point lies on $N$. The contact geodesics of $M$ lying on such an $N$ determine on $N$ a path geometry. $N$ is, in the usual sense, a totally geodesic submanifold of $(M, \nabla)$, for any $\nabla$ representing the contact projective structure on $M$; consequently it makes sense to restrict to $N$ such $\nabla$, and, moreover, the contact geodesics of $M$ lying on $N$ are the unparameterized geodesics of this restricted connection. If any other connection representing the contact projective structure is restricted to $N$, the geodesics of the restriction will be also the contact geodesics of $M$ lying on $N$, and so $M$ induces on $N$ a well-defined projective structure. Formally, this is evident from \eqref{lambdatransform}. Because the restriction of $\omega$ to a path vanishes, \eqref{lambdatransform} implies also that every contact geodesic of a contact projective manifold acquires a flat projective structure. Associating to the totally geodesic isotropic submanifold, $N \subset (M, H)$, the submanifold $\rho^{-1}(N)$ gives a bijection between totally geodesic isotropic submanifolds of the contact projective manifold $(M, H)$, and totally geodesic isotropic submanifolds of $(\form, \hat{\nabla})$ containing the vertical.
\end{remark}

\subsection{Canonical Subordinate Projective Structures}\label{adaptedsection}
The contact lines of the flat model contact projective structure are a subset of the family of all lines in projective space, the latter family constituting the flat model projective structure. Theorem \ref{contactadapted} shows that something similar happens for any contact projective structure, namely there are canonically determined paths transverse to the contact structure and filling out, with the given $4n-5$ parameter family of contact geodesics, a full projective structure. The contravariant part of the projective Weyl tensor of the projective structure obtained in this way is a constant multiple of the contact torsion.

\begin{definition}
A projective structure on $(M, H)$ having among its geodesics a full set of contact geodesics is said to be {\bf subordinate} to the contact projective structure determined by those contact geodesics.
\end{definition}
\noindent
There may be many projective structures subordinate to a contact projective structure. 

\begin{theorem}\label{contactadapted}
To each contact projective structure there is associated a canonical subordinate projective structure the projective Weyl tensor of which satisfies $B_{ijk}\,^{0} = - \tfrac{1}{4}\tau_{ijk}$. 
\end{theorem}

\begin{proof}
For each choice of $s$, the ambient connection determines on $M$ a connection defined by $\bar{\nabla}_{X}Y = \rho_{\ast}(\hat{\nabla}_{\hat{X}}\hat{Y})$. The difference tensor of the connection determined in this way by $\tilde{s} = fs$ and $\bar{\nabla}$ is $\Lambda_{\alpha\beta}\,^{\sigma} = \gamma_{(\alpha}\delta_{\beta)}\,^{\sigma}\,$, where $\gamma = d\log{f}$, and so these connections determine on $M$ a projective structure. It will next be shown that the geodesics of this projective structure are the projections to $M$ of the unparameterized geodesics of $\hat{\nabla}$ transverse to the vertical. For any $X \in \Gamma(TM)$ tangent to the image, $\rho(L)$, of the unparameterized geodesic, $L$, of $\hat{\nabla}$, there is $f$ so that $Z = \hat{X} + f\eul$ is tangent to $L$. That $L$ is a geodesic means that $Z \wedge \hat{\nabla}_{Z}Z$ vanishes along $L$, and computing $\rho_{\ast}(Z \wedge \hat{\nabla}_{Z}Z )$ shows that $X \wedge \bar{\nabla}_{X}X$ vanishes along $\rho(L)$, so that $\rho(L)$ is an unparameterized geodesic of $\bar{\nabla}$. By \eqref{hatnablaalphaOmega}, $\bar{\nabla}\theta = \tfrac{1}{2}\omega$, so $\bar{\nabla}$ admits a full set of contact geodesics, and it is evident that these are the images in $M$ of geodesics of $\hat{\nabla}$ tangent to $\ker \alpha$. If $\nabla$ is the connection associated to $\theta$ by Theorem A, then the difference tensor of $\bar{\nabla}$ and $\nabla$ is $\A_{\alpha\beta}\,^{\gamma}$. Lemma \ref{differencetensor} and the explicit formula, \eqref{oabcdefined}, for $\A_{\alpha\beta}\,^{\gamma}$ show that $\bar{\nabla}$ has the same contact geodesics as does $\nabla$. This implies that the projective structure determined on $M$ by $\hat{\nabla}$ is subordinate to the given contact projective structure. 

The condition $B_{ijk}\,^{0} = -\tfrac{1}{4}\tau_{ijk}$ is an invariant condition, as both sides rescale in the same way under a change of scale. In order to compute $B_{ijk}\,^{0}$ it is necessary to work with a symmetric representative of the projective structure. Let $\nabla^{\prime}$ be the symmetric part of $\bar{\nabla}$ and let $R^{\prime}_{\alpha\beta\gamma}\,^{\delta}$ be its curvature tensor. The difference tensor of $\nabla^{\prime}$ and $\nabla$ is $\Lambda_{\alpha\beta}\,^{\gamma} = \A_{(\alpha\beta)}\,^{\gamma} - \tfrac{1}{2}\tau_{\alpha\beta}\,^{\gamma}$. Direct computation gives
\begin{align}\label{subinter}
&\nabla^{\prime}_{\alpha}\theta_{\beta} = \tfrac{1}{2}\omega_{\alpha\beta},& &\nabla^{\prime}_{\alpha}\omega_{ij} = 2\delta_{\alpha}\,^{0}P_{[ji]} - \tfrac{1}{2}\tau_{ij}\,^{p}\omega_{p\alpha},& &\nabla^{\prime}_{\alpha}\omega_{i0} = 2P_{\alpha i} + \tfrac{1}{2}Q_{\alpha i}.
\end{align}
Because $\nabla^{\prime}$ is torsion free the Ricci identity and \eqref{subinter} give 
\begin{align*}
&R^{\prime}_{ijk}\,^{0} = -2\nabla^{\prime}_{[i}\nabla^{\prime}_{j]}\theta_{k} = - \nabla^{\prime}_{[i}\omega_{j]k} = \tfrac{1}{2}\nabla^{\prime}_{j}\omega_{ki} = - \tfrac{1}{4}\tau_{ijk}.
\end{align*}
By definition, $B_{\alpha\beta\gamma}\,^{\sigma} = R^{\prime}_{\alpha\beta\gamma}\,^{\sigma} + \delta_{[\alpha}\,^{\sigma}P^{\prime}_{\beta]\gamma} - 2P^{\prime}_{[\alpha\beta]}\delta_{\gamma}\,^{\sigma}$, where $R^{\prime}_{\alpha\beta}$ is the Ricci curvature of $\nabla^{\prime}$ and $P^{\prime}_{\alpha\beta} = \tfrac{1}{2(n-1)}(R^{\prime}_{\alpha\beta} - \tfrac{1}{n}R^{\prime}_{[\alpha\beta]})$, and so $B_{ijk}\,^{0} = R^{\prime}_{ijk}\,^{0} = - \tfrac{1}{4}\tau_{ijk}$.
\end{proof}

When the contact torsion vanishes the canonical subordinate projective structure is linked to the given contact projective structure in a stronger way; the Thomas ambient connection and the contact projective ambient connection are the same. By \eqref{volumetransform}, a square-root, $\form$, of the bundle of positive contact one-forms is canonically regarded as a $\frac{1}{2n}$th-root of $\candens$. Under this identification the top exterior power, $\Omega^{n}$, of the ambient symplectic structure on $\form$ is naturally identified with the canonical ambient volume form. If the contact torsion vanishes, then the ambient connection $\hat{\nabla}$ is torsion free and parallelizes the volume form, $\Omega^{n}$. By Theorem B, $\hat{\nabla}$ satisfies all the conditions of Theorem \ref{thomasambient}. The volume form determined on $M$ by the choice of section, $s$, is identified with $\theta \wedge d\theta^{n-1}$, where $\theta$ is the contact one-form determined by $s$. The proof of Theorem \ref{thomasambient} shows that $\hat{\nabla}$ determines on $M$ a projective structure for which the connection $\bar{\nabla}$ is the unique torsion-free affine connection representing the projective structure and parallelizing $\theta \wedge d\theta^{n-1}$.

On a contact manifold $(M, H)$, a projective structure is called \textbf{contact adapted}, if among the unparameterized geodesics of the projective structure there is a full set of contact geodesics, and if the projective Weyl tensor satisfies $B_{\alpha\beta\gamma}\,^{0} = 0$.

\begin{corollary}\label{contactadapted2}
To each contact projective structure with vanishing contact torsion there is associated a canonical subordinate contact adapted projective structure, and every contact adapted projective structure arises in this way.
\end{corollary}

\begin{proof}
In this case $\nabla^{\prime} = \bar{\nabla}$ and \eqref{subinter} shows $\bar{\nabla}_{\gamma}\omega_{\alpha\beta} = 4P_{\gamma[\alpha}\delta_{\beta]}\,^{0}$ and computing as above shows that the curvature satisfies $- \bar{R}_{\alpha\beta\gamma}\,^{0} = 2\delta_{[\alpha}\,^{0}P_{\beta]\gamma} + 2\delta_{\gamma}\,^{0}P_{[\alpha\beta]}$. The identification of the projective and contact projective ambient connections shows that the same tensor $P_{\alpha\beta}$ is determined by either $\bar{\nabla}$ or $\nabla$, and from this there follows $B_{\alpha\beta\gamma}\,^{0} = 0$. Moreover, given a contact adapted projective structure it is easy to see that it must be the canonical projective structure subordinate to the contact projective structure that it determines.
\end{proof}
\noindent
In \cite{Harrison} it was observed that the ambient connection associated by Theorem \ref{thomasambient} to a contact adapted projective structure makes parallel the canonical symplectic structure determined on the ambient space by the contact structure on the base. As was remarked in the introduction, even in the case of vanishing contact torsion the results here are stronger than those of \cite{Harrison} because the existence of a full projective structure is derived rather than assumed as a hypothesis.

\begin{remark}
Theorem \ref{contactadapted} could be proved directly without using the ambient connection. Let $\nabla$ be the connection associated to $\theta$ by Theorem A. Define $\nabla^{\prime}$ by letting its difference tensor with $\nabla$ be $\A_{(\alpha\beta)}\,^{\gamma} -\tfrac{1}{2}\tau_{\alpha\beta}\,^{\gamma}$, (so $\nabla^{\prime}$ is the symmetric part of $\bar{\nabla}$), and observe that the torsion free connection, $\nabla^{\prime}$, admits the given full set of contact geodesics. Direct computation using \eqref{lambdatransform}-\eqref{lambda0ip}, \eqref{ptransform}, \eqref{qtransform}, and \eqref{p0itransform} shows that the connections associated by this construction to different choices of scale have the same unparameterized geodesics, so determine a projective structure subordinate to the given contact projective structure. It is straightforward to check directly that when the contact torsion vanishes, $\bar{\nabla} = \nabla^{\prime}$ is the unique representative of this projective structure making parallel the volume $\theta \wedge d\theta^{n-1}$.
\end{remark}

\begin{remark}\label{distinguishedgeodesicsprep}
A Cartan connection gives a notion of development of paths. Once there is in hand the canonical Cartan connection of Theorem C this can be used to recover the existence of the canonical subordinate projective structure of Theorem \ref{contactadapted}, as will be explained in Section \ref{inducedsection}.
\end{remark}

\subsubsection{Obstruction to Flatness for a Contact Projective Structure}\label{obstructionsection}
\begin{theorem}\label{vantheorem}
If the contact torsion vanishes the components of the curvature of the ambient connection are all expressible in terms of $P_{ij}$, $C_{ijk}$, $W_{ijkl}$, and their covariant derivatives, as follows.
\begin{align}
&\label{A4} A_{ij} = \tfrac{1}{2(1-n)}\left(\nabla_{p}C_{ij}\,^{p} + 2W_{pijq}P^{pq}\right),\\
&\label{A5}U_{ijk} = -C_{[ij]k} =  \tfrac{1}{1-2n}\nabla_{p}W_{ijk}\,^{p},\\
&\label{A6}B_{i} = \tfrac{1}{2(n-1)}\left(\nabla_{p}\nabla_{q}C^{qp}\,_{i}- 2C_{ipq}P^{pq}\right),\\
&\label{A7}V_{ij} = -A_{[ij]} = \tfrac{1}{2(n-1)}\left(\nabla_{p}C_{[ij]}\,^{p} - W_{ijpq}P^{pq}\right).
\end{align}
A contact projective structure of dimension at least $5$ is flat if and only if the contact torsion and the contact projective Weyl tensor vanish, and a three-dimensional contact projective structure is flat if and only if the contact projective Cotton tensor vanishes.
\end{theorem}
\noindent
Theorem \ref{vantheorem} follows also by applying, to the Cartan connection of Theorem C, Proposition \ref{normalcharacter} and Theorem 2.9 of \cite{Tanaka-equivalence} on the vanishing of harmonic curvature components. 

\begin{proof}
Because $\tau_{ij}\,^{k} = 0$, the first Bianchi identity for $\hat{\nabla}$ gives $\hat{R}_{ij0P} = -2\hat{R}_{0[ij]P}$, which implies $U_{ijk} = -C_{[ij]k}$ and the truth of \eqref{A4} coupled with \eqref{vijdefined} implies that $A_{[ij]} = -V_{ij}$. Alternatively, $A_{[ij]} = -V_{ij}$ follows from \eqref{aijdefined} and \eqref{vijdefined} using \eqref{p0idefined}, \eqref{pi0defined}, and \eqref{p00defined}, and $(2n-1)C_{[ij]k} = \nabla_{p}W_{ijk}\,^{p}$ follows from \eqref{tracednablaweyl}. There remains to prove \eqref{A4} and \eqref{A6}. 
By \eqref{wdefined} and the Ricci identity,
\begin{equation}\label{wp}
 W_{pij}\,^{q}P_{q}\,^{p} = 2\nabla_{[i}\nabla_{p]}P_{j}\,^{p} - 2(n-1)P_{ip}P_{j}\,^{p} - \omega_{ij}P_{p}\,^{q}P_{q}\,^{p} - \nabla_{0}P_{ij}.
\end{equation}
Expanding the definition, \eqref{aijdefined}, of $A_{ij}$, using \eqref{p0idefined}-\eqref{p00defined}, and computing $2(n-1)A_{ij} + \nabla_{p}C_{ij}\,^{p}$ using \eqref{cdefined} gives twice the right hand side of \eqref{wp}, proving \eqref{A4}. To prove \eqref{A6} proceed as follows. Tracing the Ricci identity applied to $\nabla_{[p}\nabla_{q]}P_{ij}$ in $p$ and $q$ and using \eqref{firstbianchicontracted} gives $\nabla_{p}\nabla^{p}P_{ij} = (1-n)\nabla_{0}P_{ij}$, which with \eqref{cdefined}, \eqref{secondbianchizero}, and  \eqref{secondbianchizeroriccitrace}, gives
\begin{equation}\label{nablacfirstindex}
-\tfrac{1}{2}\nabla_{p}C^{p}\,_{ij} = \nabla_{0}P_{ij} + \tfrac{1}{2n-1}(\nabla_{i}\nabla_{p}P_{j}\,^{p} + \nabla_{j}\nabla_{p}P_{i}\,^{p}).
\end{equation}
Applying the Ricci identity gives
\begin{equation}\label{tracedricid1}
2\nabla_{[q}\nabla_{k]}\nabla_{p}P_{i}\,^{p} = -\omega_{qk}\nabla_{0}\nabla_{p}P_{i}\,^{p} - R_{qki}\,^{s}\nabla_{p}P_{s}\,^{p}.
\end{equation}
Tracing \eqref{tracedricid1} on $q$ and $k$ and subsituting into the traced covariant derivative of \eqref{nablacfirstindex} gives
\begin{equation}\label{nablanablacprelim}
-\tfrac{1}{2}\nabla_{q}\nabla_{p}C^{pq}\,_{i} = \nabla_{p}\nabla_{0}P_{i}\,^{p} + \tfrac{1}{2n-1}(\nabla_{q}\nabla_{i}\nabla_{p}P^{qp} + (1-n)\nabla_{0}\nabla_{p}P_{i}\,^{p} - 2nP_{iq}\nabla_{p}P^{pq}).
\end{equation}
The Ricci identity gives $\nabla_{p}\nabla_{0}P_{i}\,^{p} = \nabla_{0}\nabla_{p}P_{i}\,^{p} + R_{0pi}\,^{q}P_{q}\,^{p}$. Using this in \eqref{nablanablacprelim}; using \eqref{cdefined} to rewrite $R_{0ijk}$ in terms of $C_{ijk}$; and substituting into the result the trace of \eqref{tracedricid1} in $q$ and $i$ gives
\begin{align}
&\label{nablanablac}\nabla_{q}\nabla_{p}C^{pq}\,_{i} - 2C_{pqi}P^{pq}=\\ 
&\notag 4P^{pq}\nabla_{p}P_{iq} + \tfrac{1}{2n-1}(2(1-n)\nabla_{0}\nabla_{p}P_{i}\,^{p} + 4(n+2)P_{iq}\nabla_{p}P^{pq}-2\nabla_{i}\nabla_{p}\nabla_{q}P^{pq}).
\end{align}
Expanding \eqref{bidefined} using \eqref{p0idefined}-\eqref{p00defined} and substituting into \eqref{nablanablac} shows
\begin{equation}\label{almostA6}
2(1-n)B_{i} + \nabla_{q}\nabla_{p}C^{pq}\,_{i} =  2C_{pqi}P^{pq} + 8P^{pq}\nabla_{[i}P_{p]q} - \tfrac{12}{2n-1}P_{ip}\nabla_{q}P^{pq}.
\end{equation}
The proof of \eqref{tracednablaweyl} shows that the right hand side of \eqref{almostA6} is expressible in terms of $P^{pq}\nabla_{s}W_{ipq}\,^{s} = (2n-1)P^{pq}C_{[ip]q}$, and this gives \eqref{A6}.

If the ambient connection is flat the proof of Theorem \ref{contactadapted} shows that it is the Thomas ambient connection of a flat subordinate projective structure, and that the Thomas ambient connection parallelizes a symplectic form. From this there follows that $(L, \Omega, \hat{\nabla})$ is locally equivalent to $(\standrep, \bar{\Omega})$ with the flat Euclidean connection, and by functoriality of the ambient construction this shows that the given contact projective structure is flat, so that the vanishing of the curvature of the ambient connection implies the contact projective structure is flat. The preceeding shows that the vanishing of $\tau_{ij}\,^{k}$ and $W_{ijk}\,^{l}$ (or of $C_{ijk}$ in three dimensions) implies the ambient connection is flat. 
\end{proof}

The existence of non-flat contact projective structures with vanishing contact torsion may be demonstrated using the set-up of Section \ref{affinespace}. By Theorem \ref{spaceofstructures} and the accompanying discussion, the difference tensor, $\Lambda$, of the flat contact projective structure, $[\bar{\nabla}]$, described in Section \ref{flatmodel}, and the most general contact torsion-free contact projective structure, $[\nabla]$, is an arbitrary section of $\At$. The curvature of the representative $\nabla \in [\nabla]$ associated to $\theta$ by Theorem A may be computed explicitly in terms of $\Lambda$. Using the basic facts about irreducible representations of the symplectic group found in \cite{Weyl-Classical} or \cite{Fulton-Harris}, the following facts may be proved. Given a fixed point, $p$, on a contact manifold of dimension at least $5$, and an arbitrary trace-free tensor $a_{ijkl}$ such that $a_{[ijk]l} = 0$, $a_{[ij]kl} = a_{ijkl} = a_{ij(kl)}$, there exists a contact projective structure with vanishing contact torsion such that the value of $W_{ijkl}$ at the point $p$ equals $a_{ijkl}$. Given a fixed point, $p$, on a three-dimensional contact manifold, and an arbitrary tensor such that $a_{(ijk)} = a_{ijk}$, there exists a contact projective structure such that the value of $C_{ijk}$ at the point $p$ equals $a_{ijk}$. Finally there may be constructed a contact projective structure for which the component, $\hat{S}_{i0}$, of the traced curvature of the ambient connection does not vanish at a point. The proof uses  Weyl's First Main Theorem for the invariants of the symplectic group, (Theorem 6.1.A, p. 167, \cite{Weyl-Classical}), the complete reducibility of representations of the symplectic group, and the Littlewood-Richardson rules for the symplectic group.

\section{Canonical Cartan Connection}\label{cartancsection}
The canonical Cartan connection of Theorem C is constructed from the ambient connection by using the tractor formalism developed by \v{C}ap-Gover in \cite{Cap-Gover} and following the pattern of the conformal case treated in \cite{Cap-Gover-Conformal}. This approach was suggested to the author by Rod Gover. Section \ref{algebraic} describes the algebraic normalizations on the curvature of the Cartan connection of Theorem C.

\subsection{Definition of Tractor Connection and Canonical Cartan Connection}\label{tractorbundlesection}
Choose a square-root, $\form$, on the co-oriented contact manifold, $(M, H)$. Let $\delta_{r}:\form \to \form$ denote the principal $\rea^{\times}$-action of $r$. Denote by $\emf[k]$ the real line bundle associated to $\form$ by the representation of $\rea^{\times}$ on $\rea$ given by $s \cdot t = s^{-k}t$. Sections $h \in \Gamma(\emf[k])$ are in canonical bijection with functions $\tilde{h}:\form \to \rea$ of homogeneity $k$, in the sense that $\tilde{h}(\delta_{r}(p)) = r^{k}\tilde{h}(p)$. It will be convenient to write $\emf_{i_{1}\dots i_{p}}^{j_{1}\dots j_{q}}[k] = \emf[k]\tensor (\tensor^{q}(H))\tensor(\tensor^{p}(H^{\ast}))$. The Jacobi identity shows that the Lie algebra, $\vect(\form)$, of vector fields on $\form$ is graded as a Lie algebra by the homogeneity degree, $\vect(\form) = \oplus_{k \in \integer}\vect_{k}(\form)$. More generally, if $S$ is a tensor of homogeneity $l$ and $X \in \vect_{k}(\form)$, then $\lie_{X}S$ has homogeneity $k + l$. If $X \in \vect_{k}(\form)$ and $Y \in \vect_{l}(\form)$, expanding $2\hat{\nabla}\alpha(X, Y) = \Omega(X, Y)$ shows that $\hat{\nabla}_{X}Y \in \vect_{k+l}(\form)$.

Let $(\standrep, \sOmega)$ be the standard representation of $G$, as in Section \ref{flatmodel}, and choose coordinates so that $\sOmega$ satisfies \eqref{standardOmega}. Call a frame, $\{F_{I}\}$, on $\form$, {\bf symplectic}, if $\Omega(F_{I}, F_{J}) = \sOmega_{IJ}$. Call a symplectic frame {\bf adapted} if $F_{\infty} = \eul$. Because $\Omega$ is homogeneous of degree $2$, and $\eul \in \vect_{0}(\form)$, it must be that $F_{0} \in \vect_{-2}(\form)$ and all the $F_{i} \in \vect_{-1}(\form)$. By definition, $2\alpha(F_{0}) = 1$. An adapted symplectic frame is easily constructed as follows. For $\theta = s^{\ast}(\alpha)$, choose on $M$ a $\theta$-adapted coframe and dual frame such that $\omega = \frac{1}{2}\sOmega_{ij}\theta^{i}\wedge\theta^{j}$, and set $F_{\infty} = \eul$, $F_{i} = t^{-1}\hat{E}_{i}$, and $F_{0} = \tfrac{1}{2}t^{-2}\hat{\rb}$. $T\form$ is canonically filtered by the vertical subbundle, $T^{2}\form$, and its $\Omega$-skew complement, $T^{1}\form = \ker \alpha$. 

The \textbf{tractor bundle}, $\tractor$, is the rank 2n quotient of $T\form$ by an $\rea^{\times}$ action, $P_{r}$, on $T\form$ covering $\delta_{r}$ and leaving invariant $\vect_{-1}(\form)$, and defined by $P_{r}(Z) = r^{-1}\delta_{r^{-1}}\,^{\ast}(Z)$. By construction, the space of sections, $\Gamma(\tractor)$, is identified with $\vect_{-1}(\form)$. Because $\Omega$ has homogeneity $2$, it descends to $\tractor$ as a fiberwise symplectic form, also denoted $\Omega$. The $P_{r}$-invariant filtration $T^{2}\form \subset T^{1}\form \subset T\form$ descends to a filtration $\tractor^{2} \subset \tractor^{1} \subset \tractor$. 
\begin{lemma}\label{gradedtractorlemma}
There is a canonical isomorphism of graded vector bundles
\begin{equation}\label{gradedtractor}
\graded\tractor  = \tractor^{2} \oplus \tractor^{1}/\tractor^{2} \oplus \tractor/\tractor^{1}\simeq \emf[-1]\oplus\emf^{i}[-1]\oplus\emf[1].
\end{equation}
\end{lemma}

\begin{proof}
Sections of $\tractor^{2}$ are in bijection with vector fields on $\form$ having the form $\tilde{h}\eul$, where $h\in \Gamma(\emf[-1])$. This shows $\emf[-1] \simeq \tractor^{2}$. If $Z \in \vect_{-1}$, then the homogeneity $1$ function, $\Omega(\eul, Z)$, corresponds to a section of $\emf[1]$, and it vanishes if and only if $Z$ corresponds to a section of $\tractor^{1}$. This shows $\tractor/\tractor^{1} \simeq \emf[1]$. For $X \in \Gamma(TM)$ and $h\in \Gamma(\emf[-1])$, $\tilde{h}\hat{X}\in\vect_{-1}(\form)$ so corresponds to a section of $\tractor$. The image in $\tractor/\tractor^{2}$ of the section of $\tractor$ so determined does not depend on the choice of horizontal lift, $\hat{X}$, as any two horizontal lifts of $X$ differ by a vector field of the form $g\eul$, where $g$ has homogeneity $0$, and as $g\tilde{h}\eul$ corresponds to a section of $\tractor^{2}$. The section of $\tractor/\tractor^{2}$ determined by $X$ and $h$ lies in $\tractor^{1}/\tractor^{2}$ if and only if $X$ is a section of $H$, and this shows $H \tensor \emf[-1] \simeq \tractor^{1}/\tractor^{2}$. 
\end{proof}

The \textbf{bundle of filtered symplectic frames} in $\tractor$ is the principal $P$-bundle, $\pi:\adaptedframe \to M$, the fiber over $x$ of which comprises all filtration preserving symplectic linear isomorphisms $u:(\standrep, \sOmega) \to (\tractor_{x}, \Omega_{x})$. The fibers of the vertical subbundle $V\adaptedframe \subset T\adaptedframe$ are linearly isomorphic to $\p$, so $V\adaptedframe$ carries a filtration induced by that of $\p$; namely for $i = 0, 1, 2$, $V^{i}_{u}\adaptedframe \simeq \f_{i}$. The tangent bundle $T\adaptedframe$ carries a canonical filtration modeled on the filtration, $\f_{i}$, of $\g$ and uniquely determined by the requirements that the projection onto $TM$ induced by $\pi_{\ast}$ is a filtered surjection and that inclusion of $V\adaptedframe$ is a filtered injection. By construction the quotient $T\G/T^{-1}\G$ has the orientation induced via $\pi_{\ast}$ from the chosen orientation on $TM/H$. The {\bf projective tractor bundle}, $\proj(\tractor) \to M$, does not depend on the choice of $\form$, so is intrinsically associated to the co-oriented contact structure on $M$. The fiber $\proj(\tractor_{x})$ is a copy of $\proj(\standrep)$ osculating $M$ at $x$, and the point of contact is the distinguished point, $\proj(\tractor^{2}_{x})$. A principal $\bar{P}$ bundle, $\bar{\adaptedframe} \to M$, the {\bf bundle of filtered projective symplectic frames}, is defined by letting the fiber over $x \in M$ comprise all projective symplectic linear isomorphisms $u:\proj(\standrep) \to \proj(\tractor_{x})$ mapping $\proj(\standrep^{i})$ to $\proj(\tractor^{i}_{x})$ for $i = 1, 2$. Explicit computation shows that the differential of the lift to $\form$ of an element of $\contactop(M)$ commutes with $P_{r}$, and so descends to a vector bundle automorphism of $\tractor$. The induced fiber bundle automorphism on $\proj(\tractor)$ does not depend on the choice of square-root $\form$, so that $\contactop(M)$ acts naturally on $\tractor$ and $\proj(\tractor)$ and consequently also on $\G$ and $\bar{\adaptedframe}$ and their associated bundles.

Any representation $\nu:P \to GL(\rep)$ determines an associated bundle, $\assrep = \adaptedframe\times_{\nu(P)}\rep$, and sections, $t$, of $\assrep$ are in canonical bijection with $P$-equivariant functions, $\tilde{t}:\adaptedframe \to \rep$. The projection defined by mapping $(u, v) \in \adaptedframe \times \standrep$ to $u(v) \in \tractor_{\pi(u)}$, where $u \in \adaptedframe$ is viewed as an isomorphism $u:\standrep \to \tractor_{\pi(u)}$, descends to give an isomorphism of filtered vector bundles between the associated bundle, $\adaptedframe \times_{P}\standrep$, and $\tractor$. The equivariant function corresponding to the section, $t$, of $\tractor$, is defined by $\tilde{t}:u \to u^{-1}(t(\pi(u)))$. The symplectic form on $\tractor$ is recovered by $\Omega_{\pi(u)}(s, t) = \standomega(\tilde{s}(u), \tilde{t}(u))$. Because $\tilde{P}$ acts trivially on $\standrep^{2}$, the bundle of frames in $\tractor^{2}$ may be identified with the $P/\tilde{P} = \rea^{\times}$ principal bundle $\adaptedframe/\tilde{P} \to M$. By Lemma \ref{gradedtractorlemma}, the subbundle $\tractor^{2}$ is canonically identified with $\emf[-1]$, and the bundle of frames of $\emf[-1]$ is by definition $\form$. This recovers $\form$ from $\adaptedframe$ as the quotient $\adaptedframe/\tilde{P}$. $\adaptedframe$ is a principal $\tilde{P}$-bundle over $\form$, and it could have been realized directly as a reduction of the frame bundle of $T\form$, the fiber  over $p \in \form$ of which comprises all linear symplectic isomorphisms $u:\standrep \to T_{p}\form$ such that $u(v_{\infty}) = \eul_{p}$. A useful by-product of this observation is the identification of $T\form$ and the associated bundle $\adaptedframe \times_{\tilde{P}}\standrep$ as filtered vector bundles. Because $\tilde{P}$ acts trivially on $\standrep^{2}$, the constant function $v_{\infty}$ on $\adaptedframe$ is $\tilde{P}$-equivariant, and it is the $\tilde{P}$-equivariant function corresponding to $\eul \in \Gamma(T\form)$. 

Given $X \in \Gamma(TM)$ and $t \in \Gamma(\tractor)$ represented by $Z \in \vect_{-1}(\form)$, define $\tnabla_{X}t$ to be the section of $\tractor$ represented by $\hat{\nabla}_{\hat{X}}Z\in \vect_{-1}(\form)$. The connection $\tnabla$ defined on $\tractor$ is independent of the choice of $s$ as any two horizontal lifts of $X$ differ by a vector field of the form $f\eul$, where $f$ has homogeneity $0$, and $\hat{\nabla}_{\eul}Z = 0$ for any $Z \in \vect_{-1}(\form)$. It is easily checked that the connection, $\tnabla$, on $\tractor$ parallelizes $\Omega$. The connection, $\tnabla$, is non-degenerate in the sense that for any section, $t$, of $\tractor^{2}$ there exists a vector field, $X$, on $M$ such that $\tnabla_{X}t \notin \tractor^{2}$. In the terminology of \cite{Cap-Gover}, $\tractor$ is a standard tractor bundle, $\adaptedframe$ is an adapted frame bundle, and the connection, $\tnabla$, on $\tractor$ is a tractor connection. Theorem 2.7 of \cite{Cap-Gover} shows that a tractor connection, $\tnabla$, determines on $\pi:\adaptedframe\to M$ a Cartan connection, $\eta$. The construction of $\eta$ is recalled here, though \cite{Cap-Gover} should be consulted for the verifications. For each $u \in \adaptedframe$, $X \in \Gamma(T\adaptedframe)$, and $t \in \Gamma(\tractor)$, the linearity of the map $\standrep \to \standrep$ defined by $\tilde{t}(u) \to (\widetilde{\tnabla_{\pi_{\ast}(X)}t})(u) - \lie_{X}\tilde{t}(u)$ follows from the Leibniz rules for $\lie_{X}$ and $\tnabla_{\pi_{\ast}(X)}$. It is easily checked that this linear map preserves $\bar{\Omega}$, so lies in $\g$.  Because $\g$ is simple, it acts effectively on $\standrep$, so $\eta(X)$ is defined uniquely by 
\begin{align}\label{etadefined}
&\eta(X)\cdot\tilde{t} = \widetilde{\tnabla_{\pi_{\ast}(X)}t} - \lie_{X}\tilde{t},& &X \in \Gamma(T\adaptedframe),\,\, t \in \Gamma(\tractor).&
\end{align}
Note that the right hand side of \eqref{etadefined} is a $P$-equivariant function if and only if $X$ is invariant under the action of $P$. That $\eta$ is a Cartan connection is the content of Theorem 2.7 of \cite{Cap-Gover}. The non-degeneracy of the tractor connection is the key point in showing that $\eta$ induces a linear isomorphism $T_{u}\adaptedframe\simeq\g$. The other properties of a Cartan connection follow from careful unraveling of the definition of $\eta$.

Let $\rep$ be any $P$-module for which the induced action of $\p$ extends to an action of $\g$. The Cartan connection, $\eta$, induces on the associated bundle, $\assrep = \adaptedframe \times_{P}\rep$, a covariant differentiation defined by 
\begin{align}\label{assrepdefined}
&\widetilde{\etanabla_{X}t} = \lie_{\bar{X}}\tilde{t} + \eta(\bar{X})\cdot \tilde{t},& &X \in \Gamma(TM), \quad t \in \Gamma(\assrep),&
\end{align}
where $\bar{X}$ is any horizontal lift of $X$ to a vector field on $\adaptedframe$. Because $\tilde{t}$ is $P$-equivariant, $\lie_{V}\tilde{t} + \eta(V)\cdot \tilde{t} = 0$ for any vertical vector field, $V$, on $\adaptedframe$, and this shows that $\bar{X}$ may be chosen to be invariant under the action of $P$, in which case the right hand side of \eqref{assrepdefined} is $P$-equivariant, and this shows that \eqref{assrepdefined} is well-defined. When $\rep = \standrep$, the induced covariant differentiation, $\etanabla$, is the tractor connection $\tnabla$. It may be checked that $\eta$ is a $(\g, \tilde{P})$ Cartan connection on the $\tilde{P}$ principal bundle $\adaptedframe \to \form$. In this case, the induced covariant differentiation on the associated bundle $T\form \simeq \adaptedframe\times_{\tilde{P}}\standrep$ is the ambient connection, $\hat{\nabla}$. 

The curvature of $\etanabla$ is defined by $R^{\eta}(X, Y)(t) = [\etanabla_{X},\etanabla_{Y}]t- \etanabla_{[X, Y]}t$. Recall that the curvature of $\eta$ is defined by $K = d\eta + \eta \wedge \eta$. It is straightforward to check (as in Proposition 2.9 of \cite{Cap-Gover}) that the curvature of $\etanabla$ acts via the curvature of $\eta$,
\begin{equation}\label{asscurv}
R^{\eta}(X, Y)(t) = \kappa(\eta(\bar{X}), \eta(\bar{Y}))\cdot t = K(\bar{X}, \bar{Y})\cdot t.
\end{equation}
Define the curvature function, $\kappa \in C^{\infty}(\adaptedframe, \bigwedge^{2}(\g/\p)^{\ast} \tensor \g)$, of $\eta$ by
\begin{equation}\label{curvaturefunction}
\kappa(u)(h_{1}, h_{2}) = K_{u}(\eta^{-1}(h_{1}), \eta^{-1}(h_{2})) = [h_{1}, h_{2}] - \eta_{u}([\eta^{-1}(h_{1}), \eta^{-1}(h_{2})]).
\end{equation}
The $\integer$-grading of $\g$ induces a $\integer$-grading of $\Lambda^{2}(\g_{-}^{\ast})\tensor \g$ by homogeneity degree. If $\sigma \in \Lambda^{2}(\g_{-}^{\ast})\tensor \g$, then $\sigma^{(i)}$ maps $\g_{j} \times \g_{k}$ into $\g_{j+k+i}$. A Cartan connection is \textbf{regular} if $\kappa^{(i)} = 0$ for $i \leq 0$. 

\subsection{Explicit Description of Tractor Connection}
In this section, $\eta$ is the Cartan connection constructed from the ambient connection. The tractor connection is described explicitly for use in the algebraic characterization of the curvature normalizations imposed on $\eta$.

Via the horizontal lift, the choice of $s$ determines a $P_{r}$-invariant splitting $T\form = T^{2}\form \oplus \hat{H} \oplus \spn\{\hat{\rb}\}$. As $\vect_{-1}(\form)$ is naturally identified with $\Gamma(\tractor)$, this splitting induces an isomorphism (depending on $s$), $\tractor \simeq \emf[-1] \oplus \emf^{i}[-1] \oplus \emf[-1] = \emf^{A}$. If $z^{A} \in \Gamma(\emf^{A})$, the corresponding $Z \in \vect_{-1}(\form)$ has the form 
\begin{equation}\label{tractorcoord}
Z = \tilde{z}^{\infty}\eul + \tilde{z}^{i}\hat{E}_{i} + \tfrac{1}{2}\tilde{z}^{0}\hat{\rb},
\end{equation} 
where the $\tilde{z}^{P}$ are the corresponding homogeneity $-1$ functions on $\form$. Though it would be natural to work with the isomorphism $\tractor \simeq  \emf[-1] \oplus \emf^{i}[-1] \oplus \emf[1]$ induced by $s$ and \eqref{gradedtractor}, to do so would require systematic utilization of weighted tensors, e.g. the canonical symplectic pairing on $\emf^{i}[-1]$ induced by the Levi form $H \times H \to TM/H \simeq \emf[2]$. Because of the way the ambient connection was described it is more convenient to work with the isomorphism $\tractor \simeq \emf^{A}$ induced directly by the splitting. $\Gamma(\tractor^{\ast})$ is identified with the space of one-forms on $\form$ homogeneous of degree $1$. The splitting of $T\form$ given by $s$ induces a dual isomorphism $\tractor^{\ast} \simeq \emf[1]\oplus\emf_{i}[1]\oplus\emf[1] = \emf_{A}$. Raise and lower indices using the symplectic form, $\Omega$, on $\tractor$, and let $\emf_{A}$ be the vector bundle dual to $\emf^{A}$. If $z_{A} = \begin{pmatrix} z_{\infty} & z_{j} & z_{0} \end{pmatrix}$ is a section of $\emf_{A}$, the corresponding homogeneity $1$ one-form, $\beta$, on $\form$ may be written as 
\begin{equation}\label{dualtractorcoords}
\beta = \tilde{z}_{\infty}d\log{t} + \tilde{z}_{i}\rho^{\ast}(\theta^{i}) + 2\tilde{z}_{0}\rho^{\ast}(\theta).
\end{equation}
The connection, $\nabla$, associated to $\theta = s^{\ast}(\alpha)$ by Theorem A induces on $\emf[k]$ and $\emf^{i}[k]$ connections which will also be denoted by $\nabla$. By definition of the horizontal lift, $\widetilde{\nabla_{\alpha}h} = d\tilde{h}(\hat{E}_{\alpha})$. Direct computation of $\hat{\nabla}_{\hat{E}_{i}}Z$ and $\hat{\nabla}_{\hat{\rb}}Z$ using \eqref{tractorcoord} yields
\begin{equation}\label{tnablavup}
 \tnabla_{\alpha}\begin{pmatrix} z^{\infty} \\ z^{p} \\ z^{0}\end{pmatrix} = \begin{pmatrix} \nabla_{\alpha}z^{\infty}\\\nabla_{\alpha}z^{p} \\ \nabla_{\alpha}z^{0}\end{pmatrix} + \begin{pmatrix} 0 & P_{\alpha q} & \frac{1}{2}P_{\alpha 0}\\ \delta_{\alpha}\,^{p} & \delta_{\alpha}\,^{0}(2P_{q}\,^{p} + Q_{q}\,^{p}) & P_{\alpha}\,^{p}\\ 2\delta_{\alpha}\,^{0}& -\omega_{\alpha q}& 0\\\end{pmatrix}\begin{pmatrix}z^{\infty} \\ z^{q} \\ z^{0}\end{pmatrix}.
\end{equation}
The curvature of $\tnabla$ is two-form on $M$ taking values in the endomorphism bundle of $\tractor$. It is computed by $[\tnabla_{\alpha}, \tnabla_{\beta}]z^{A} + \tau_{\alpha\beta}\,^{\sigma}\tnabla_{\sigma}z^{A} = \tR_{\alpha\beta B}\,^{A}z^{B}$, where $\tau_{\alpha\beta}\,^{\sigma}$ is the torsion of $\nabla$. For $h \in \emf[k]$ and $z^{p} \in \emf^{i}[k]$ the Ricci identity gives
\begin{align}
&\label{vric1}2\nabla_{[\alpha}\nabla_{\beta]}h + \tau_{\alpha\beta}\,^{\sigma}\nabla_{\sigma}h = 0,&
&2\nabla_{[\alpha}\nabla_{\beta]}z^{p} + \tau_{\alpha\beta}\,^{\sigma}\nabla_{\sigma}z^{p} = R_{\alpha\beta q}\,^{p}z^{q}.
\end{align} 
Direct computation using \eqref{tnablavup}, \eqref{vric1}, \eqref{uipp}, and \eqref{wtraces} gives
\begin{align}\label{trijab}
&\tR_{ijB}\,^{A} = \begin{pmatrix} -Q_{[ij]} & U_{ijq} & V_{ij}\\ \tau_{ij}\,^{p} & W_{ijq}\,^{p} & U_{ij}\,^{p}\\ 0 & -\tau_{ijq} & Q_{[ij]} \end{pmatrix},&
&\tR_{0iB}\,^{A}  = \begin{pmatrix} 2U_{si}\,^{s} &  A_{iq} & B_{i}\\ Q_{i}\,^{p}& C_{iq}\,^{p} & A_{i}\,^{p}\\ 0 & -Q_{iq}& -2U_{si}\,^{s} \end{pmatrix},&
\end{align}
where the tensors $U$, $V$, $A$, and $B$ are defined in \eqref{uijkdefined}-\eqref{bidefined}, and it should be recalled that $Q_{ij} = -2W_{ipj}\,^{p}$. If the contact torsion vanishes, these simplify to
\begin{align*}
\tR_{ijB}\,^{A} =  \begin{pmatrix}0 & -C_{[ij]q} & -A_{[ij]}\\ 0 & W_{ijq}\,^{p} & -C_{[ij]}\,^{p}\\ 0 & 0 & 0 \end{pmatrix},&
&\tR_{0iB}\,^{A} = \begin{pmatrix} 0 &  A_{iq} & B_{i}\\ 0 & C_{iq}\,^{p} & A_{i}\,^{p}\\ 0 & 0 & 0\end{pmatrix}.&
\end{align*}

\begin{example}[Contact Hessian]
\noindent
For $h \in \Gamma(\emf[k])$, the exterior derivative $d\tilde{h}$ is a one-form on $\form$, homogeneous of degree $k$, so may be regarded as a section of $\tractor^{\ast}$. Using \eqref{dualtractorcoords}, the components of $d\tilde{h}$, viewed as a section of $\emf_{A}$, may be written as
\begin{equation}\label{1jeth}
D_{A}h = \begin{pmatrix} kh & \nabla_{j}h & \frac{1}{2}\nabla_{0}h\end{pmatrix}.
\end{equation}
Taking $h \in \emf[1]$ and applying the tractor connection to \eqref{1jeth} gives
\begin{align}
&\label{invop2} \tnabla_{i}D_{A}h = \left(0 \quad L_{ij}h \quad \tfrac{1}{2}\nabla_{i}\nabla_{0}h - P_{i}\,^{p}\nabla_{p}h - \tfrac{1}{2}P_{i0}h \right),
\end{align}
\begin{align}
\notag \tnabla_{0}D_{A}h = 
 \left(0 \quad \nabla_{0}\nabla_{j}h + (2P_{pj} + Q_{pj})\nabla^{p}h  - P_{0j}h \quad  \tfrac{1}{2}\nabla_{0}\nabla_{0}h + P_{0p}\nabla^{p}h - \tfrac{1}{2}P_{00}h \right).
\end{align}
\eqref{invop2} shows that the operator, $L:\Gamma(\emf[1]) \to \Gamma(\emf_{ij}[1])$, defined by 
\begin{equation}\label{invop1}
L_{ij}h = \nabla_{i}\nabla_{j}h + \tfrac{1}{2}\omega_{ij}\nabla_{0}h - P_{ij}h,
\end{equation}
is an invariant differential operator. That is, $\tilde{L}_{ij}h = L_{ij}h$. Alternatively, the invariance $L_{ij}h$ follows from the transformation rules for the covariant derivatives of a $-\lambda/2n$-density, $h \in \Gamma(\emf[\lambda])$,  
\begin{align}\label{densitytransform}
&\tilde{\nabla}_{i}h - \nabla_{i}h = \lambda\gamma_{i}h,&
&f^{2}\tilde{\nabla}_{0}h - \nabla_{0}h = \lambda\gamma_{0}h + 2\gamma^{p}\nabla_{p}h.&
\end{align}
Because the connection induced by $\nabla$ on the canonical bundle is flat, the Ricci identity implies $2\nabla_{[i}\nabla_{j]}h + \tau_{ij}\,^{\alpha}\nabla_{\alpha}h = 0$. Tracing this gives $\nabla_{p}\nabla^{p}h = (1-n)\nabla_{0}h$, so that the trace free part of $\nabla_{i}\nabla_{j}h - P_{ij}h$ is $L_{ij}h$. The skew part of $L_{ij}h$, $L_{[ij]}h$ is first order and vanishes when the contact torsion vanishes, so it makes sense to focus attention on the \textbf{contact Hessian}, $L_{(ij)}h = \nabla_{(i}\nabla_{j)}h - P_{(ij)}h$.


The projection $T\form \to \tractor$ induces a projection $\tensor^{\ast}(T\form) \to \tensor^{\ast}(\tractor)$ on the full tensor bundles. Denote by $h_{A_{1} \dots A_{k}}^{B_{1} \dots B_{l}}$ a section of $\tensor^{\ast}(\tractor)$. The tractor operator $D_{I}h_{A_{1} \dots A_{k}}^{B_{1} \dots B_{l}}$ is defined by lifting the section $h_{A_{1} \dots A_{k}}^{B_{1} \dots B_{l}}$ to an appropriately weighted tensor on $\form$ to which the ambient connection is applied. The resulting tensor is projected down to give a section of $\tensor^{\ast}(\tractor)$. This allows the iteration of the tractor $D$ operator. 
\end{example}

\subsection{Contact Projective Structures Induced by Cartan Connections}\label{inducedsection}
Let $(M, H)$ be a co-oriented contact manifold and $\bar{\pi}:\bar{\adaptedframe} \to M$ the corresponding $\bar{P}$ principal bundle of filtered projective symplectic frames. The purpose of this section is to explain how a Cartan connection, $\eta$, on $\bar{\G}$, satisfying an appropriate compatibility condition induces on $M$ a contact projective structure. Each choice of a square-root, $\rho:\form \to M$, determines a $P$ principal bundle $\pi:\adaptedframe \to M$ and a double covering $\Phi:\adaptedframe \to \pframe$. The quotient $\adaptedframe/\tilde{P}$ recovers $\form$, and $\tilde{\pi}:\adaptedframe \to \form$ is the induced $\tilde{P}$ principal bundle. A Cartan connection, $\eta$, on $\ppi:\pframe \to M$ pulls back via $\Phi$ to a Cartan connection, also denoted $\eta$, on $\pi:\adaptedframe \to M$, and will also be viewed as a Cartan connection on $\tilde{\pi}:\adaptedframe \to \form$. The group of principal bundle automorphisms of $\adaptedframe$ acts on the Cartan connections on $\adaptedframe$ by pullback. A principal bundle automorphism covering the identity on $M$ will be called a \textbf{change of gauge}. A change of gauge $\psi:\adaptedframe \to \adaptedframe$ has the form $\psi(u) = u\cdot \tilde{p}(u)$ where $\tilde{p}:\adaptedframe \to P$ is the $P$-equivariant function corresponding to a section, $p$, of the bundle $\adaptedframe\times_{\Ad(P)}P$. The pullback Cartan connection has the form $\psi^{\ast}(\eta) = \Ad(\tilde{p}^{-1})(\eta) + \tilde{p}^{\ast}(\omega_{P})$, where $\omega_{P}$ is the Maurer-Cartan form of $P$. The only case in which $\psi^{\ast}(\eta) = \eta$ is when $\tilde{p}$ is a constant function and $\Ad(\tilde{p})$ acts trivially on $\g$. A choice of $\form$ determines a unique lift of a bundle automorphism of $\bar{\G}$ to a bundle automorphism of $\G$. 

Call an $\eta$, on $\ppi:\pframe \to M$, \textbf{compatible} with the co-oriented contact structure on $M$ if $\eta$ is regular; if the isomorphisms, $\eta:T_{u}\pframe \to \g$, are filtration-preserving; and if the induced isomorphisms, $T_{u}\pframe/T^{-1}_{u}\pframe \simeq \g/\f_{-1}$, are orientation-preserving. It is evident that compatibility is a gauge invariant condition. Compatibility implies $\ppi_{\ast}(\eta^{-1}(\f_{-1})) = H$ and that the image of $\ppi_{\ast}(\eta^{-1}(e_{0}))$ is consistent with the given orientation on $TM/H$, and so, with regularity, that $\eta$ induces the given co-oriented contact distribution and its conformal symplectic structure. A $(\g, P)$ Cartan connection $\eta$ on $\pi:\adaptedframe \to M$ will be called compatible if it is induced by a compatible $(\g, \bar{P})$ Cartan connection on $\ppi:\pframe \to M$. Henceforth assume all $\eta$ are compatible. 

Viewing $\eta$ as a Cartan connection on the $\tilde{P}$ principal bundle, $\tilde{\pi}:\adaptedframe \to \form$, denote by $\etanabla$ the affine connection induced on $T\form = \adaptedframe\times_{\tilde{P}}\standrep$ by \eqref{assrepdefined}. Call \textbf{contact non-slip} a compatible Cartan connection for which $\etanabla_{X}\eul - X \in \Gamma(T^{2}\form)$ for all $X \in \Gamma(T^{1}\form)$. Using the definition of $\etanabla$ it may be checked that $\eta$ is contact non-slip if and only if $(\etanabla_{X}\alpha)(Y) = \frac{1}{2}\Omega(X, Y)$ for all $X, Y \in \Gamma(T^{1}\form)$. Likewise $\etanabla\eul = \delta$ if and only if $\etanabla\alpha = \tfrac{1}{2}\Omega$. If $\eta$ is the Cartan connection defined from the ambient connection, $\hat{\nabla}$, by \eqref{etadefined}, then $\etanabla \eul = \delta$ is condition \conref{con1} of Theorem B. Lemma \ref{gaugelemma} shows that this condition should be interpreted as fixing a gauge. 
\begin{lemma}\label{gaugelemma}
A given compatible $(\g, P)$ Cartan connection on $\pi:\G \to M$ is equivalent by a $G_{0}$-valued change of gauge to a contact non-slip Cartan connection. Two gauge equivalent contact non-slip $(\g, P)$ Cartan connections are related by a $P^{+}$-valued change of gauge. Moreover, a given contact non-slip Cartan connection is equivalent by a unique $P^{+}$-valued change of gauge to a Cartan connection for which $\etanabla \eul = \delta$; consequently, every compatible $(\g, P)$ Cartan connection on $\pi:\adaptedframe \to M$ is gauge equivalent to a unique Cartan connection for which $\etanabla \eul = \delta$. 
\end{lemma}
\begin{proof}
For a compatible $\eta$ let $\sigma:U\subset \form \to \adaptedframe$ be a local section and set $\phi = \sigma^{\ast}(\eta)$. Under a change of gauge, $\psi:\G \to \G$, as above, $\phi$ is replaced by $\Ad(b^{-1})(\phi) + b^{\ast}(\omega_{P})$, where $b = \tilde{p}\circ \sigma:U \to P$. Write $\phi(X) = X^{Q}\phi_{Q}$ with respect to an adapted symplectic frame, $\{F_{I}\}$, on $\form$, and represent $\phi_{Q} \in \g$ by a matrix,
\begin{align*}
\phi_{Q} = \begin{pmatrix}A_{Q}\,^{\infty} & \ast &\ast \\ A_{Q}\,^{p} & \ast & \ast \\ A_{Q}\,^{0} & -A_{Q q} & -A_{Q}\,^{\infty} \end{pmatrix}.
\end{align*}
Because $e_{\infty} \in \p$, $\eta(X^{e_{\infty}}) = e_{\infty}$, and because $\tilde{\pi}_{\ast}(X^{e_{\infty}}) = \eul$, $\sigma_{\ast}(\eul) - X^{e_{\infty}} \in \eta^{-1}(\tilde{\p})$, so $\phi(\eul)\cdot v_{\infty} = e_{\infty}\cdot v_{\infty} = v_{\infty}$. This shows $A_{\infty}\,^{\infty} = 1$, $A_{\infty}\,^{p} = 0$, and $A_{\infty}\,^{0} = 0$. From the compatibility of $\eta$ there follows $\phi(X) \in \f_{-1}$ if and only if $X \in T^{1}\form$, and this implies $A_{i}\,^{0} = 0$. Regularity, $A_{i}\,^{0} = 0$, and $A_{\infty}\,^{0} = 0$, imply that $A_{i}\,^{p}A_{j}\,^{q}\omega_{pq} = \bar{\Omega}(\phi(F_{i})\cdot v_{\infty}, \phi(F_{j})\cdot v_{\infty}) = \tfrac{1}{2}\bar{\Omega}([\phi(F_{i}), \phi(\F_{j})]\cdot v_{\infty}, v_{\infty}) = -\bar{\alpha}(\phi([F_{i}, F_{j}])\cdot v_{\infty}) = A_{0}\,^{0}\omega_{ij}$. Compatibility implies $A_{0}\,^{0} = c^{2}> 0$. Let $B_{q}\,^{p}$ be the unique symplectic matrix such that $cB_{s}\,^{j}A_{i}\,^{s} = \delta_{i}\,^{j}$. The change of gauge given by $b:\form \to G_{0}$ defined by $b^{-1} = \exp(c^{-1}v_{\infty}\tensor v^{\infty} + B_{p}\,^{q}v_{q}\tensor v^{p} + cv_{0}\tensor v^{0})$ produces a contact non-slip connection. As the image in $P/\tilde{P}$ of $b$ is non-trivial, the contribution of $b^{\ast}(\omega_{P})$ under this change of gauge is non-trivial, effecting the components $A_{I}\,^{\infty}$. Because the adjoint action of $P$ on $P/\tilde{P}$ is trivial, $b^{\ast}(\omega_{P})(\eul) \in \tilde{\p}$, and the identity $A_{\infty}\,^{\infty} = 1$ persists. Thus every compatible Cartan connection is gauge equivalent to a contact non-slip Cartan connection. Now suppose given a contact non-slip Cartan connection and set $\gamma^{p} = A_{0}\,^{p}$ and $\gamma_{0} = A_{0}\,^{\infty} - \gamma_{s}A_{0}\,^{s}$. For $b:U \to P^{+}$ defined by $b = \exp(\gamma_{p}v_{\infty}\tensor v^{p} + \gamma^{p}v_{p}\tensor v^{0} + \gamma_{0}v_{\infty}\tensor v^{0})$, $\Ad(b^{-1})(\phi) + b^{\ast}(\omega_{P})$ satisfies $\etanabla\eul = \delta$. Similar easy computations show that gauge equivalent contact non-slip Cartan connections must be related by a change of gauge of the form $b = \exp(\gamma_{p}v_{\infty}\tensor v^{p} + \gamma^{p}v_{p}\tensor v^{0} + \gamma_{0}v_{\infty}\tensor v^{0}) \in P^{+}$. The same computation shows that a change of gauge preserving the condition $\etanabla\eul = \delta$ must be trivial. 
\end{proof}

Recall Cartan's notion of development as described in section $4$ of chapter $5$ of \cite{Sharpe}. If $\tilde{\gamma}:[a,b] = I \to \G$ is any lift of $\gamma:[a, b] = I \to M$, there is a unique $c:I \to G$ such that $c^{\ast}(\omega_{G}) = \tilde{\gamma}^{\ast}(\eta)$ and $c(a) = e \in G$, where $\omega_{G}$ is the Maurer-Cartan form of $G$. While the curve $c$ depends on the choice of lift, its projection $\Psi \circ c :I \to G/P$ is well-defined up to the action of $P$ on $G/P$, and is called the \textbf{development} of $\gamma$. (Note that in the proof of Proposition 4.13 in \cite{Sharpe} Lemma 4.12 of \cite{Sharpe} is misapplied to incorrectly conclude that the projection into $G/P$ of $c$ does not depend on the choice of lift; the issue is the choice of a basepoint). Two parameterized curves with the same path develop onto curves with the same path, so it makes sense to say that the development of a path on $M$ is a member of a $P$-invariant family of paths in $\proj(\standrep)$, e.g. it makes sense to say that a path develops onto a line. Since the family of contact lines (resp. all lines) in $G/P$ is invariant under the action of $P$ on $G/P$, the statement that a path in $M$ develops onto a contact line (resp. line) in $G/P$ is independent also of the choice of Cartan connection representing the isomorphism class of $\eta$. If $\eta$ is the Cartan connection built from the ambient connection, direct computation as in the proof of Proposition 8.3.5 of \cite{Sharpe} and using the explicit description of $\eta$ implicit in \eqref{tnablavup} shows that the paths of $\eta$ which develop onto isotropic lines are exactly the contact geodesics and the paths which develop onto lines are exactly the geodesics of the subordinate projective connection of Theorem \ref{contactadapted}. Consequently, for an arbitrary compatible $\eta$, it makes sense to define the induced contact projective structure to comprise those contact paths in $M$ which develop onto contact lines in $G/P$. There are in general many Cartan connections on $\adaptedframe$ inducing a given contact projective structure. 

A compatible $\eta$ induces on the $G$ principal bundle, $\G\times_{\tilde{P}}G \to \form$, the $G$ principal connection, $\mu$, induced by the connection $\mu_{u, b} = \Ad(b^{-1})(\pi_{\G}^{\ast}(\eta)) + \pi_{G}^{\ast}(\omega_{G})$ on $\G \times G$, where $\pi_{\G}$ and $\pi_{G}$ are the projections onto the factors of $\G \times G$, and $\omega_{G}$ is the Maurer-Cartan form of $G$. $\G\times_{\tilde{P}}G$ is naturally identified with the symplectic frame bundle of $T\form$ and the induced $\etanabla$ is just the affine connection on $\form$ induced by $\mu$. Moreover, if $c:I \to G$ is associated as above to the lift $\tilde{\gamma}$ of $\gamma:I \to M$, then $\bar{\gamma}(t) = [\tilde{\gamma}(t), c^{-1}(t)]$ is a horizontal curve, $\bar{\gamma}:I \to\G\times_{\tilde{P}}G$, covering $\gamma$. Using this it may be checked that the paths in $\form$ which develop via $\eta$ onto lines in $\standrep$ are exactly the unparameterized geodesics of $\etanabla$. Evidently these are the horizontal lifts of the paths in $M$ which develop onto lines in $\proj(\standrep)$.

Using \eqref{assrepdefined}, it is straightforward to show that, for compatible $\eta$, $\etanabla \Omega = 0$ and $\etanabla_{\eul}\eul = \eul$. Let $\tau^{\eta}$ and $R^{\eta}$ be the torsion and curvature of $\etanabla$. Using that $\eta$ is contact non-slip and $\etanabla_{\eul}\eul = \eul$, and computing directly the curvature shows that, for all $X, Y \in \Gamma(T^{1}\form)$, $\tau^{\eta}(X, Y) - R^{\eta}(X, Y)\eul \in T^{2}\form$. The curvature function, $\kappa(h_{1}, h_{2})$, of $\eta$ vanishes whenever either $h_{1}$ or $h_{2}$ is in $\p$, so, by \eqref{asscurv}, $R^{\eta}(\eul, X)\eul = \kappa(\eta(\bar{\eul}), \eta(\bar{X}))\cdot v_{\infty} = 0$, and consequently $\tau^{\eta}(\eul, X) \in T^{2}\form$ for all $X \in \Gamma(T^{1}\form)$. Next the induced contact projective structure is described more analytically. For a section, $s:M \to \form$, let $\theta = s^{\ast}(\alpha)$ and define an affine connection, $\nabla$, on $M$ by $\nabla_{X}Y = \rho_{\ast}(\etanabla_{\hat{X}}\hat{Y})$. For a contact non-slip $\eta$, the restrictions to $T^{1}\form \times T^{1}\form$ of $\etanabla\alpha$ and $\tfrac{1}{2}\Omega$ agree, and this implies that $\sym \nabla \theta$ vanishes when restricted to $H \times H$. By Lemma \ref{contactgeodesics}, the geodesics of $\nabla$ determine a contact projective structure. Because $\eta$ is contact non-slip, 
\begin{align}\label{induceintermed}
\etanabla_{\eul}\hat{Y} - \hat{Y} = (\etanabla_{\hat{Y}}\eul- \hat{Y}) + \tau^{\eta}(\eul, \hat{Y}) \in T^{2}\form \quad \text{for}\quad  Y \in \Gamma(H). 
\end{align}
In particular, $\rho_{\ast}(\etanabla_{\eul}\hat{Y}) = Y$. If $\tilde{s} = fs$ and direct computation using \eqref{induceintermed} shows that the connection, $\tilde{\nabla}$, induced on $M$ by $\tilde{s}$ determines the same contact projective structure as does $\nabla$. Moreover it is straightforward to check, as in the proof of Theorem \ref{contactadapted}, that the contact geodesics of the induced contact projective structure are the images in $M$ of the paths of the geodesics of $\etanabla$ transverse to the vertical and tangent to $\ker\alpha$. By the remarks above the geodesics of $\etanabla$ tangent to $\ker\alpha$ develop onto contact lines in $\standrep$, and hence their images in $M$ develop onto contact lines in $\proj(\standrep)$, so the induced contact projective structure is the same as the contact projective structure defined using development. To show via the tractor connection that gauge equivalent contact non-slip Cartan connections, $\eta$ and $\mu$, induce the same contact projective structure, choose $s:M \to \form$ and compute the difference tensor, $\Lambda$, of the connections $\rho^{\ast}(\nabla^{\mu}_{\hat{X}}\hat{Y})$ and $\rho^{\ast}(\etanabla_{\hat{X}}\hat{Y})$. Computing directly using \eqref{assrepdefined} and Lemma \ref{gaugelemma} shows that $\Lambda_{ij}\,^{k} = \gamma^{k}\omega_{ij} + \gamma_{j}\delta_{i}\,^{k}$ and $\Lambda_{ij}\,^{0} = \omega_{ij}$. By Lemma \ref{differencetensor}, the induced contact projective structures are the same.

\subsection{Algebraic Characterization of Curvature of Canonical Cartan Connection}\label{algebraic}
The purpose of this section is to reformulate conditions \conref{con4} and \conref{con5} of Theorem B as algebraic conditions on the curvature function of the associated Cartan connection, and to prove Theorem C. 

There is a co-boundary operator, $\partial:\bigwedge^{k}((\g_{-})^{\ast})\tensor\g \to\bigwedge^{k+1}((\g_{-})^{\ast})\tensor\g$, (see \cite{Yamaguchi}), defined when $k = 2$ by 
\begin{align*}
&\partial \phi(x_{1}, x_{2}, x_{3}) = \text{Cycle}\left[x_{1}\cdot\phi(x_{2}, x_{3}) - \phi([x_{1}, x_{2}], x_{3})\right]. 
\end{align*} 
Any invariant, symmetric bilinear form, $B$, on $\g$ is a constant multiple of the Killing form, and so $B$ gives an identification $(\p^{+})^{\ast} \simeq \g_{-}$. Let $e_{\alpha}$ be a basis for $\g_{-}$ and let $e^{\alpha}$ be the $B$-dual basis for $\p^{+}$. The formal adjoint with respect to $B$ of the co-boundary $\partial$ is a boundary $\partial^{\ast}:\bigwedge^{2}(\g_{-}^{\ast})\tensor \g \to g_{-}^{\ast}\tensor \g$ defined by $B(\partial^{\ast}\psi, \phi) = -B(\psi, \partial\phi)$ for $\psi \in \partial^{\ast}:\bigwedge^{2}(\g_{-}^{\ast})\tensor \g$ and $\phi \in g_{-}^{\ast}\tensor \g$. Via the linear isomorphism, $\g/\p \simeq \g_{-}$, induced by $B$, the adjoint action makes $\g_{-}$ a $P$-module; it can be checked that $\partial^{\ast}$ is $P$-equivariant. Recall that a Cartan connection is {\bf normal} if $\partial^{\ast}\kappa = 0$. 

Because the tractor connection, $\tnabla = \etanabla$, is the induced connection on the associated bundle, $\G\times_{P}\standrep$, its curvature, $\tR$, may be computed from the curvature of $\eta$, as in \eqref{asscurv}. By \eqref{trijab}, $\kappa_{-2} = 0$; in particular, $\eta$ is regular. The isomorphism $\eta:T_{u}\adaptedframe \to \g$ induces an adjoint map $\eta^{t}:\g^{\ast} \to T^{\ast}_{u}\adaptedframe$. Because $\ker B(e^{0}, -) = \f_{-1}$, and $\Ad(P)(\g_{-2}) = \g_{-2}$, the image under $\pi_{\ast}$ of the horizontal subbundle of $T^{\ast}\G$ spanned by $\eta^{t}(B(e^{0}, -))$ is the annihilator of $H$. Let $e_{I}$ and $v_{I}$ be as in Section \ref{flatmodel}, and define a frame on $M$ by $E_{i} = \pi_{\ast}(\eta^{-1}(e_{i}))$ and $\rb = \pi_{\ast}(\eta^{-1}(2e_{0}))$. Choose $\theta$ such that $\theta(\rb) = 1$ and $\theta$ annihilates $H$. Because $\pi^{\ast}(\theta)$ is a multiple of $\eta^{t}(B(e^{0},- ))$ and $\kappa_{-2} = 0$, $d\theta(E_{\alpha}, E_{\beta}) = 2\omega_{\alpha\beta}$, so that $\theta$ is a contact one-form with Reeb field $\rb$. Choose also a $\theta$-adapted coframe dual to the given frame. By \eqref{asscurv} the components, $\kappa(e_{i}, e_{j})$, of the curvature function are identified with the components, $\tR_{ij}$, of the curvature of the tractor connection, and $2\kappa(e_{0}, e_{i}) = \tcurv_{0i}$. Choose $B$ so that with respect to the standard representation, $\rho:\rsymp \to GL(\standrep)$, $B(x, y) = \tfrac{1}{2}\tr(\rho(x)\rho(y))$. With respect to the basis $v_{I}\tensor v^{J}$ of $\g\subset \standrep\tensor \standrep^{\ast}$, 
\begin{align*}
&e_{i} = v_{i}\tensor v^{\infty} - \omega_{iq}v_{0}\tensor v^{q},& &e_{0} = v_{0}\tensor v^{\infty},& &e_{\infty} = v_{\infty}\tensor v^{\infty} - v_{0}\tensor v^{0},&\\
 &e^{i} = v_{\infty}\tensor v^{i} + \omega^{pi}v_{p}\tensor v^{0},& &e^{0} = 2v_{\infty}\tensor v^{0}.&
\end{align*}
\begin{proposition}\label{normalcharacter}
The Cartan connection, $\eta$, is normal if and only if the underlying contact projective structure has vanishing contact torsion.
\end{proposition}

\begin{proof}
In general $\partial^{\ast}\kappa(x) = [e^{\alpha}, \kappa(e_{\alpha}, x)] - \tfrac{1}{2}\kappa(e_{\alpha},[e^{\alpha}, x]_{-})$, which gives
\begin{align}
&\notag \partial^{\ast}\kappa(e_{j}) = [e^{0},\kappa(e_{0}, e_{j})] + [e^{i}, \kappa(e_{i}, e_{j})],&
&
\partial^{\ast}\kappa(e_{0}) = [e^{i},\kappa(e_{i}, e_{0})] - \tfrac{1}{2}\kappa(e_{i}, [e^{i}, e_{0}]_{-}).
\end{align}
Direct computation of $\partial^{\ast}\kappa(e_{j})$ using \eqref{uipp}, and direct computation of $\partial^{\ast}\kappa(e_{0})$ using  $[e^{i}, e_{0}] = -\omega^{ij}e_{j}$, \eqref{wtraces}, \eqref{ctraces}, and \eqref{uipp} give
\begin{align*}
&\partial^{\ast}\kappa(e_{j}) = \begin{pmatrix}0 & -Q_{jq} + \frac{1}{2}Q_{qj} & -2U_{sj}\,^{s} \\ 0 & \tau_{j}\,^{p}\,_{q} + \tau_{jq}\,^{p} & -Q_{j}\,^{p} + \frac{1}{2}Q^{p}\,_{j}\\ 0 & 0 & 0  \end{pmatrix}&\\ &\partial^{\ast}\kappa(e_{0}) = \begin{pmatrix}0 & U_{sq}\,^{s} + \frac{1}{2}U_{s}\,^{s}\,_{q} & \frac{1}{2}V_{s}\,^{s}\\ 0& \frac{1}{2}(Q^{p}\,_{q} + Q_{q}\,^{p})& U_{s}\,^{ps} + \frac{1}{2}U_{s}\,^{sp} \\0 &0 &0 \end{pmatrix} &
\end{align*}
If $\partial^{\ast}\kappa(e_{j}) = 0$, then $\tau_{j(pq)} = 0$, from which follows $0 = 2\tau_{i(jk)} - 2\tau_{j(ik)} = 3\tau_{ijk}$, so that $\partial^{\ast}\kappa = 0$ implies the contact torsion vanishes. On the other hand, $\tau_{ijk} = 0$ implies $Q_{ij} = 0$ and $P_{[\alpha\beta]} = 0$, and, by \eqref{uipp}, it implies also $V_{s}\,^{s} = 0$ and $U_{s}\,^{ps} = 0$. \eqref{p0idefined}, \eqref{uipp}, and \eqref{hatsi0} show $U_{s}\,^{s}\,_{q} + 2U_{sq}\,^{s} = -2\nabla^{p}P_{pq} + (2n-1)P_{0q} = \nabla^{p}Q_{pq}$, so that $U_{s}\,^{s}\,_{q} + 2U_{sq}\,^{s}$ vanishes also. Thus $\tau_{ij}\,^{k} = 0$ implies $\partial^{\ast}\kappa(e_{\alpha}) = 0$.
\end{proof}

The goal of the remainder of this subsection is to describe the $P$ submodule, $\pcurv \subset C^{2}(\g_{-}, \g)$, appearing in the statement of Theorem C. Given a compatible $(\g, P)$-Cartan connection on the principal $P$-bundle $\adaptedframe \to M$, let $\tnabla = \etanabla$ be the induced covariant differentiation on the associated tractor bundle $\tractor = \adaptedframe\times_{P}\standrep$. Let $\tcurv_{\alpha\beta}$ be the components with respect to $E_{\alpha}$ of the $\g$-valued curvature two-form of $\tnabla$, so $\tcurv_{ij} = \kappa(e_{i}, e_{j})$ and $\tcurv_{0i} = 2\kappa(e_{0}, e_{i})$. Because the $\eta$ constructed from the ambient connection always satisfies $\kappa_{-2} = 0$, it will be assumed that $\kappa_{-2} = 0$, which implies that $\tcurv$ has the form
\begin{align}
&\label{gentcurv0i}\tcurv_{ijA}\,^{B} = \begin{pmatrix} a_{ij} & c_{ijq} & e_{ij}\\ b_{ij}\,^{p} & d_{ijq}\,^{p} & c_{ij}\,^{p} \\ 0 & -b_{ijq} & -a_{ij}\end{pmatrix},& &\tcurv_{0iA}\,^{B} = 2\begin{pmatrix} a_{i} & c_{iq} & e_{i} \\ b_{i}\,^{p} & d_{iq}\,^{p} & c_{i}\,^{p} \\ 0 & -b_{iq} & -a_{i}\end{pmatrix},
\end{align}
where $d_{ij[qp]} = 0$ and  $d_{i[qp]} = 0$ are forced by the requirement that these matrices are in $\g$. Rewritten in terms of $\kappa$, the Bianchi identity, $dK = [K, \eta]$, becomes
\begin{equation}\label{kappabianchi}
\partial \kappa (h_{1}, h_{2}, h_{3}) + \text{Cycle}\left[\kappa(\kappa_{-}(h_{1}, h_{2}), h_{3}) + \lie_{\eta^{-1}(h_{1})}(\kappa(h_{2}, h_{3}) )\right] = 0.
\end{equation}
In \eqref{kappabianchi}, $ \lie_{\eta^{-1}(h_{1})}(\kappa(h_{2}, h_{3}))$ is the Lie derivative of the function $\kappa(h_{2}, h_{3}):\adaptedframe \to \g$. Because $\kappa_{-2} = 0$, \eqref{kappabianchi} implies $(\partial\kappa)_{-2}(e_{\alpha}, e_{\beta}, e_{\gamma}) = 0$, giving 
\begin{align}
&\label{park21}\text{Cycle}\left( [e_{i},\kappa_{-1}(e_{j}, e_{k})] +\omega_{ij}\kappa_{-2}(e_{0}, e_{k})\right) = 0,\\
&\label{park22}[e_{0}, \kappa_{0}(e_{i}, e_{j})] + [e_{i}, \kappa_{-1}(e_{j}, e_{0})] + [e_{j}, \kappa_{-1}(e_{0}, e_{i})] = 0.
\end{align}
As $[e_{k},\kappa_{-1}(e_{i}, e_{j})] = b_{ijk}e_{0}$, \eqref{park21} shows $b_{[ijk]} = 0$. Direct computation using \eqref{park22} gives $a_{ij} = -2b_{[ij]}$ so that $\tcurv_{\alpha\beta B}\,^{A}$ has the form
\begin{align}&\label{gentcurvij}
\tcurv_{ijA}\,^{B} = \begin{pmatrix} -2b_{[ij]} & c_{ijq} & e_{ij}\\ b_{ij}\,^{p} & d_{ijq}\,^{p} & c_{ij}\,^{p} \\ 0 & -b_{ijq} & 2b_{[ij]}\end{pmatrix},& &\tcurv_{0iA}\,^{B} = 2\begin{pmatrix} a_{i} & c_{iq} & e_{i} \\ b_{i}\,^{p} & d_{iq}\,^{p} & c_{i}\,^{p} \\ 0 & -b_{iq} & -a_{i}\end{pmatrix}.&
\end{align}
Because of the curvature and torsion normalizations of Theorem B, the curvature of the Cartan connection, $\eta$, built from the ambient connection has the form \eqref{trijab} and satisfies additionally,
\begin{align*} 
&\tau_{[ijk]} = 0,& &\tau_{pi}\,^{p} = 0 = \tau_{p}\,^{p}\,_{i},& &P_{[ij]} = -2Q_{[ij]},& &Q_{p}\,^{p} = 0,& \\ 
&W_{p}\,^{p}\,_{i}\,^{j} = 0,& & W_{ipj}\,^{p} = -\tfrac{1}{2}Q_{ij},& &U_{pi}\,^{p} = P_{[0i]},& &C_{p}\,^{pi} = 0,& & A_{p}\,^{p} = 0,&
\end{align*}
which motivate Definition \ref{pcurvdefined}. In Definition \ref{pcurvdefined} and Lemma \ref{psubmodule}, $\kappa$ should be understood as an element of the $P$-module $C^{2}(\g_{-}, \g)$.
\begin{definition}\label{pcurvdefined}
The subspace $\pcurv \subset C^{2}(\g_{-}, \g)$ comprises $\kappa$ of the form \eqref{gentcurvij}, i.e. 
\begin{align}
&\label{kappa-2}\kappa_{-2} = 0,&&a_{ij} = -b_{[ij]},&&b_{[ijk]} = 0,& 
\end{align}
and satisfying the following additional requirements
\begin{align} 
&\label{bijktraces} b_{pi}\,^{p} = 0,& &b_{p}\,^{p}\,_{i} = 0,& &b_{p}\,^{p} = 0,&\\
&\label{pcurvcond} c_{pi}\,^{p} = a_{i},& & d_{ipj}\,^{p} = -b_{ij},& &d_{p}\,^{pi} = 0,& & d_{p}\,^{p}\,_{i}\,^{j} = 0,& & c_{p}\,^{p} = 0.&
\end{align}
\end{definition}

\begin{lemma}\label{psubmodule}
$\pcurv$ is a $P$-submodule of $C^{2}(\g_{-}, \g)$.
\end{lemma}
\begin{proof}
For $b \in P$, let $\tilde{e}_{\alpha} = \ad(b)(e_{\alpha})$, let $\tilde{\kappa} = b\cdot \kappa$, and write $\tilde{\kappa}_{\alpha\beta} = \tilde{\kappa}(e_{\alpha}, e_{\beta}) = \Ad(b^{-1})(\kappa(\tilde{e}_{\alpha}, \tilde{e}_{\beta}))$. The adjoint action of $P$ on $\g_{-2}$ is trivial which suffices to show the $P$-invariance of the condition $\kappa_{-2} = 0$. Assuming \eqref{kappa-2}, direct computation 
\begin{align}
&\partial^{\ast}\kappa(e_{j}) = \begin{pmatrix} b_{ij}\,^{i} & -2b_{jq} - d_{jiq}\,^{i} + b_{qj} & -2a_{j} + 2c_{ij}\,^{i} \\ 0 & b_{j}\,^{p}\,_{q} + b_{jq}\,^{p} & -2b_{j}\,^{p} + b^{p}\,_{j} - d_{ji}\,^{pi} \\ 0 & 0 & -b_{ij}\,^{i}\end{pmatrix},\\
&\partial^{\ast}\kappa(e_{0}) = \begin{pmatrix} -b_{l}\,^{l} & a_{q} - d_{sq}\,^{s} + \tfrac{1}{2}c_{s}\,^{s}\,_{q} & \tfrac{1}{2}e_{s}\,^{s} - 2c_{s}\,^{s} \\ \tfrac{1}{2}b_{s}\,^{sp} & b_{q}\,^{p} + b^{p}\,_{q} + \tfrac{1}{2}d_{s}\,^{s}\,_{q}\,^{p} & -d_{s}\,^{sp} + a^{p} + \tfrac{1}{2}c_{s}\,^{sp} \\ 0 & -\tfrac{1}{2}b_{s}\,^{s}\,_{q} & b_{s}\,^{s} \end{pmatrix}.
\end{align}
From these it is apparent that in the presence of \eqref{kappa-2} the requirement that $\partial^{\ast}\kappa \in C^{1}(\g_{-}, \tilde{\p})$ is equivalent to \eqref{bijktraces}. 
Representing $b \in P^{+}$ by $\exp(\gamma_{q}e^{q} + \frac{1}{2}\gamma_{0}e^{0})$ and computing directly gives $\tilde{e}_{i} = e_{i}$ and $\tilde{e}_{0} = e_{0} + \gamma^{p}e_{p}$, from which follow $\tilde{\kappa}_{ij} = \ad(b^{-1})(\kappa_{ij})$, and $\tilde{\kappa}_{0i} = \ad(b^{-1})(\kappa_{0i} + \gamma^{s}\kappa_{si})$, for $b \in P^{+}$. Further direct computation gives the components of $\tilde{\kappa}_{\alpha\beta} - \kappa_{\alpha\beta}$,
\begin{align}
&\label{skewbijgauge}\tilde{b}_{[ij]} - b_{[ij]} = -\tfrac{1}{2}\gamma^{p}b_{ijp} ,\\
&\label{bijpgauge}\tilde{b}_{ij}\,^{p} - b_{ij}\,^{p} = 0 ,\\
&\label{cijqgauge}\tilde{c}_{ijk} - c_{ijk} = -2\gamma_{k}b_{[ij]} + \gamma_{0}b_{ijk} + \gamma_{k}\gamma^{p}b_{ijp} + \gamma^{p}d_{ijkp},\\
&\label{dijqpgauge}\tilde{d}_{ijk}\,^{l} - d_{ijk}\,^{l} =  b_{ij}\,^{l}\gamma_{k} + \gamma^{l}b_{ijk} ,\\
&\label{eijgauge}\tilde{e}_{ij} - e_{ij} =  -4\gamma_{0}b_{[ij]} - 2\gamma_{0}\gamma^{p}b_{ijp} + 2\gamma^{p}c_{ijp} - \gamma^{p}\gamma^{q}d_{ijpq},\\
&\label{aigauge}\tilde{a}_{i} - a_{i} = 2\gamma^{p}b_{ip} -\gamma^{p}b_{pi} + \gamma^{q}\gamma^{p}b_{qip},\\
&\label{bijgauge}\tilde{b}_{ij} - b_{ij} = \gamma^{p}b_{pij},\\
&\label{diqpgauge}\tilde{d}_{ij}\,^{k} - d_{ij}\,^{k} =b_{i}\,^{k}\gamma_{j} + \gamma^{k}b_{ij} + \gamma^{p}\left(d_{pij}\,^{k} + \gamma_{j}b_{pi}\,^{k} + \gamma^{k}b_{pij} \right) ,\\
&\label{ciqgauge}\tilde{c}_{ij} - c_{ij} =  a_{i}\gamma_{j} + \gamma_{0}b_{ij}\\\notag &\qquad\qquad + \gamma^{p}\left(\gamma_{j}b_{ip} + d_{ijp} + c_{pij} - \gamma_{j}\gamma_{q}b_{pi}\,^{q} + \gamma^{q}d_{pijq} \right),\\
&\label{eigauge}\tilde{e}_{i} - e_{i} = 2\gamma_{0}a_{i}   + \gamma^{p}\left(2c_{ip} - 2 \gamma_{0}b_{ip} +  e_{pi}\right) \\\notag &\qquad\qquad - \gamma^{p}\gamma^{q}\left(d_{ipq} + 2\gamma_{0}b_{piq} - 2c_{piq} + \gamma^{r}d_{piqr}\right).
\end{align}
Represent $b \in G_{0}$ by $b = \exp(f^{-1}v_{\infty}\tensor v^{\infty} + F_{q}\,^{p}v_{p}\tensor v^{q} + fv_{0}\tensor v^{0})$, where $ F_{i}\,^{p}F_{j}\,^{q}\omega_{pq} = \omega_{ij}$. For such $b$ direct computation gives $\tilde{e}_{i} = fF_{i}\,^{p}e_{p}$, and $\tilde{e}_{0} = f^{2}e_{0}$, from which follow
\begin{align*}
&\tilde{\kappa}_{ij} = \ad(b^{-1})(f^{2}F_{i}\,^{p}F_{j}\,^{q}\kappa_{pq}),&
&\tilde{\kappa}_{0i} = \ad(b^{-1})(f^{3}F_{i}\,^{p}\kappa_{0p}),& 
&\text{for} \,\, b \in G_{0}.&
\end{align*}
Defining $G_{p}\,^{i}$ by $F_{i}\,^{p}G_{p}\,^{j} = \delta_{i}\,^{j}$, further direct computation gives
\begin{align}
&\label{kappaijg0} \tilde{\kappa}_{ij} = f^{2}F_{i}\,^{u}F_{j}\,^{v}\begin{pmatrix} 
a_{uv} & fc_{uvp}F_{q}\,^{p} & f^{2}e_{uv} \\
f^{-1}G_{q}\,^{p}b_{uv}\,^{q} & d_{uvr}\,^{s}G_{s}\,^{p}F_{q}\,^{r} &   fc_{uvs}F^{ps}\\
0 & -f^{-1}G_{sq}b_{uv}\,^{s} & -a_{uv}
\end{pmatrix},
\end{align}
\begin{align}
&\label{kappa0ig0} \tilde{\kappa}_{0i} =  f^{3}F_{i}\,^{u}\begin{pmatrix} a_{u} & fc_{up}F_{q}\,^{p} & f^{2}e_{u} \\
f^{-1}G_{q}\,^{p}b_{u}\,^{q} & G_{s}\,^{p}F_{q}\,^{r}d_{ur}\,^{s} & fG_{q}\,^{p}c_{u}\,^{q}\\
0 & -f^{-1}b_{up}F_{q}\,^{p} & - a_{u} 
\end{pmatrix}.
\end{align}
\eqref{skewbijgauge}, \eqref{bijgauge}, \eqref{bijpgauge}, \eqref{kappaijg0}, and \eqref{kappa0ig0} show directly the $P$-invariance of the conditions $a_{ij} = -b_{[ij]}$, $b_{[ijk]} = 0$, and \eqref{bijktraces}. The $P$-invariance of the conditions \eqref{pcurvcond} follows from \eqref{cijqgauge}, \eqref{dijqpgauge}, \eqref{aigauge}, \eqref{bijgauge}, \eqref{diqpgauge}, \eqref{ciqgauge}, \eqref{kappaijg0}, and \eqref{kappa0ig0} by direct computation.
\end{proof}
It is desirable to find a conceptual description of $\pcurv$ that makes clear its $P$-invariance. The conditions in \eqref{kappa-2} other than $\kappa_{-2} = 0$ are purely algebraic consequences of the Bianchi identities, \eqref{park21} and \eqref{park22}, and so must be $P$-invariant. As $\partial^{\ast}$ is $P$-equivariant and $\tilde{\p}$ is invariant under the adjoint action of $P$ on $\g$, the condition that $\partial^{\ast}\kappa$ take values in $\tilde{\p}$ is $P$-invariant, and so the equivalent conditions \eqref{bijktraces} are $P$-invariant as well. The formulation of the remaining conditions is {\it ad hoc} and so demonstration of their $P$-invariance requires explicit computations.

\begin{proof}[Proof of Theorem C]
The construction from the ambient connection of the Cartan connection demonstrates the existence of a regular $(\g, P)$ Cartan connection inducing the given contact projective structure, having curvature in $\pcurv$, and satisfying $\etanabla\eul = \delta$. 
Now let $\eta$ be a compatible $(\g, P)$ Cartan connection with curvature in $\pcurv$ and inducing the given contact projective structure. By Lemma \ref{gaugelemma},  $\eta$ is gauge equivalent to a Cartan connection for which the associated covariant differentiation, $\etanabla$, on $T\form = \adaptedframe \times_{\tilde{P}}\standrep$, satisfies $\etanabla \eul = \delta$. Let $\Reta$ and $\taueta$ denote the curvature and torsion of $\etanabla$, and recall that $\etanabla \Omega = 0$. From $\etanabla \eul = \delta$ there follows $\Reta(\eul, X)\eul = \taueta(\eul, X)$, and $\Reta(\eul, X)\eul = 0$ by \eqref{asscurv}, so that $\taueta(\eul, -) = 0$. This shows that $\etanabla$ satisfies conditions \conref{con1}-\conref{con3} of Theorem B. 

Because $\ker B(e^{0}, -) = \f_{-1}$, the one-form $\psi^{0} = \eta^{t}(B(e^{0}, -))$ is horizontal and annihilates $\eta^{-1}(\f_{-1})$. Because $\Ad(b)(e^{0}) = e^{0}$ for all $b \in \tilde{P}$, though not for all $b \in P$, $\psi^{0}$ is $\tilde{P}$-invariant, though it is not $P$-invariant, and hence $\psi^{0}$ is the pullback of some one-form on $\form$ which, by compatibility, annihilates $\ker \alpha$. Because $\etanabla\eul =\delta$, the frame defined in $T\form$ by $F_{I} = \tilde{\pi}_{\ast}(\eta^{-1}(e_{I}))$, is an adapted symplectic frame, and from this there follows $\tilde{\pi}^{\ast}(\alpha) = \tfrac{1}{2}\psi^{0}$. By \eqref{asscurv}, $\Reta_{PQR}\,^{S}e_{S} = \kappa(e_{P}, e_{Q})\cdot e_{R}$. It was already observed that $\Reta_{\infty QR}\,^{S} = 0$ and $\Reta_{PQ\infty}\,^{S} = \taueta_{PQ}\,^{S}$. Because $\etanabla \Omega = 0$, the Ricci identity gives $\Reta_{PQ(AB)} = \Reta_{PQAB}$. Consequently, $\Reta_{ABC}\,^{0} = -\Reta_{ABC\infty} = - \Reta_{AB\infty C} = -\taueta_{ABC}$. By the preceeding, $\Reta_{IJ} = \Reta_{IQJ}\,^{Q} = \Reta_{IqJ}\,^{q} - \taueta_{I0J}$. Now explicit computations using the conditions defining $\pcurv$ show \conref{con4} and \conref{con5}, so that $\etanabla$ satisfies all the conditions of Theorem B. The uniqueness statement in Theorem B shows $\etanabla$ must be the ambient connection of the given contact projective structure, from which there follow the identifications
\begin{align}
&\label{etaij}b_{ijk} = \tau_{ijk},& &c_{ijk} = U_{ijk},& &d_{ijkl} = W_{ijkl},& &e_{ij} = V_{ij},&\\
&\label{eta0i}2b_{ij} = Q_{ij},& &2c_{ij} = A_{ij},& &2d_{ijk} = C_{ijk},& &2e_{i} = B_{i}.&
\end{align}
The identification $\tau_{ijk} = b_{ijk}$ and Proposition \ref{normalcharacter} show that $\eta$ is normal if and only there vanishes the contact torsion of the underlying contact projective structure. Direct computation shows that a change of gauge preserving condition \conref{con1} must be trivial, so that this condition singles out the ambient connection as a distinguished representative of the given isomorphism class of Cartan connections. 
\end{proof}

\begin{remark}\label{transformremark}
By \eqref{etaij} and \eqref{eta0i}, the application to the Cartan connection constructed from the ambient connection of the transformation formulae in \eqref{skewbijgauge}-\eqref{eigauge} gives the rules for the transformations of the components of the curvature of the ambient connection under change of scale. For example, \eqref{bijpgauge} shows the invariance of the contact torsion, \eqref{dijqpgauge} recovers \eqref{weyltransform}, and \eqref{diqpgauge} recovers \eqref{ctransform}.
\end{remark}

\section{Applications}
\subsection{Lifting Projective Structures on Integral Symplectic Manifolds}
An integral symplectic manifold admits a canonical contactification, and projective structures on the symplectic manifold can be lifted to contact projective structures on the contactification. This construction provides many examples of contact projective structures. It is necessary to take the slightly unusual perspective of allowing the affine connections representing a projective structure to have torsion.

\begin{proposition}\label{adaptedrepresentative}
Given a projective structure on the symplectic manifold $(N, \omega)$, there exists a unique affine connection, $\nabla$, representing the projective structure, having torsion, $\tau$, and such that
\setcounter{condition}{0}
\begin{condition}\label{q2}
$\nabla \omega = 0$.
\end{condition}\begin{condition}\label{q3}
The trace of the torsion of $\nabla$ vanishes.
\end{condition} 
\end{proposition}

\begin{proof}
Raise and lower indices with $\omega$. Let $\nabla$ be a connection with torsion, $\tau$, representing the given projective structure, and define $\bar{\nabla}$ by letting its difference tensor with $\nabla$ be $\Lambda_{ijk} = -\nabla_{k}\omega_{ij} - \tfrac{3}{2}\tau_{[ijk]}$. As $\Lambda_{(ij)k} = 0$, $\bar{\nabla}$ represents the given projective structure. \eqref{differenceact} and \eqref{skew} show $\bar{\nabla}_{k}\omega_{ij} = 0$. Consequently $\nabla$ may be assumed to make parallel $\omega$. If $\bar{\nabla}$ is another representative of the given projective structure satisfying Condition \conref{q2}, there must be a one-form, $\gamma$, so that the difference tensor, $\Lambda_{ij}\,^{k}$, of $\bar{\nabla}$ and $\nabla$, satisfies $\Lambda_{(ij)k} = \tfrac{1}{2n+1}(\omega_{ik}\gamma_{j} + \omega_{jk}\gamma_{i})$, and $\Lambda_{k[ij]} = 0$. Because a three tensor skew in two indices and symmetric in two indices must vanish, the most general possible choice for $\Lambda_{ijk}$ is $\Lambda_{ijk} = \tfrac{2}{2n+1}(\omega_{ij}\gamma_{k} + \omega_{ik}\gamma_{j})$. Tracing the difference of the torsion tensors shows that the unique choice of $\gamma$ so that \conref{q3} holds also is $\gamma_{i} = \frac{1}{2}\tau_{ip}\,^{p}$.
\end{proof}

On a symplectic manifold, $(N, \omega)$, such that $\omega$ represents an integral cohomology class $[\omega] \in H^{2}(N, \integer)$, there is a principal $S^{1}$-bundle $\pi:M \to N$ with connection one-form, $\theta$, such that the curvature $d\theta = \pi^{\ast}\omega$. The kernel of $\theta$ is a contact structure on $M$ and the infinitesimal generator, $\rb$, of the principal $S^{1}$ action is the Reeb vector field of $\theta$. The collection of all horizontal lifts of the paths in a projective structure on $N$ determines on $M$ a contact path geometry, which will now be shown to be a contact projective structure. Fix $\nabla$ representing the given projective structure and making $\omega$ parallel. Define on $M$ a connection, $\hat{\nabla}$, by
\begin{align*}
&\hat{\nabla}_{\hat{X}}\hat{Y} = \widehat{\nabla_{X}Y},& &\hat{\nabla}_{\hat{X}}\rb = 0,& &\hat{\nabla}_{\rb}\hat{X} = \lie_{\rb}\hat{X},
\end{align*}
where $\hat{X}$ is the horizontal lift determined by $\theta$ of the vector field $X$ on $N$. It is straightforward to check that $\hat{\nabla}\theta = 0$, $\hat{\nabla}d\theta = 0$, and $i(\rb)\hat{\tau} = 0$. Moreover, if $\nabla$ is the unique representative given by Proposition \ref{adaptedrepresentative}, then the trace of the torsion of $\hat{\nabla}$ vanishes, so that $\hat{\nabla}$ is the unique connection associated to $\theta$ by Theorem A. 

\subsection{Pseudo-Hermitian Beltrami Theorem}
Recall the classical 
\begin{theorem}[Beltrami Theorem]
The projective structure determined by the geodesics of a Riemannian metric is flat if and only if the Riemannian metric has constant sectional curvatures.
\end{theorem}
\noindent
The proof, which is an exercise in using the Bianchi identities, is in \cite{Chern-Griffiths}. Next there is proved an analogous theorem for pseudo-hermitian manifolds. 

A {\bf pseudo-hermitian} structure is a contact manifold equipped with a distinguished contact one-form, $\theta$, and an almost complex structure, $J$, on $H$. Here $J: H \to H$ is a real automorphism such that $J_{i}\,^{p}J_{p}\,^{j} = - \delta_{i}\,^{j}$ and assumed to satisfy the compatibility condition
\begin{equation}\label{Jcompatibility}
[J(X), Y] + [X, J(Y)] \in \Gamma(H) \quad \forall \, X, Y \in \Gamma(H).
\end{equation}
If in addition $J$ satisfies
\begin{equation}\label{Jintegrable}
[J(X), J(Y)] - [X, Y] = J([J(X), Y] + [X, J(Y)]) \quad \forall \, X, Y \in \Gamma(H)
\end{equation}
then the pseudo-hermitian structure is said to be {\bf integrable}. By assumption \eqref{Jcompatibility}, the inner product defined on $H$ by $g_{ij} = -J_{ij}$ is symmetric, and because $d\theta$ is non-degenerate on $H$, $g$ is non-degenerate on $H$. The almost complex structure $J$ determines an $(n-1)$-dimensional complex subbundle of the complexified tangent bundle, $T_{1, 0} \subset \com TM$, comprising vector fields of the form $X + iJ(X)$ for $X \in \Gamma(H)$. Define the {\bf Levi form}, $\levi(U, V) = -i\omega(U, \bar{V})$ for $U, V \in \Gamma(T_{1, 0})$. 

In \cite{Tanaka-homogeneous} and \cite{Tanaka-Book}, Tanaka constructed a canonical affine connection associated to a pseudo-hermitian structure. In the integrable case his construction specializes as in the following theorem, the statement of which is taken from \cite{Tanaka-CR}.

\begin{theorem}[N. Tanaka]\label{tanakaconnection}
On an integrable pseudo-hermitian manifold there exists a unique affine connection, $\Bar{\nabla}$, the \textbf{pseudo-hermitian connection}, having torsion $\tau$ and satisfying:
\setcounter{condition}{0}
\begin{condition}\label{tanakacon1}
$\bar{\nabla} \theta = 0$.
\end{condition}
\begin{condition}\label{tanakacon2}
$\bar{\nabla} d\theta = 0$.\end{condition}
\begin{condition}
$\tau_{ij}\,^{\gamma} = \omega_{ij}\delta_{0}\,^{\gamma}$.\end{condition}
\begin{condition}
$J_{i}\,^{p}A_{p}\,^{j} = -A_{i}\,^{p}J_{p}\,^{j}$ where $A_{\alpha}\,^{\beta} = \tau_{0\alpha}\,^{\beta}$ is the \textbf{pseudo-hermitian torsion}.\end{condition}
\end{theorem}
\noindent
By definition $A_{0}\,^{\beta} = 0$; by Condition \conref{tanakacon2}, $A_{\alpha}\,^{0} = 0$; and by \eqref{torsionlemma1}, $A_{[ij]} = 0$. By Condition \conref{tanakacon1}, Lemma \ref{contactgeodesics} implies that the geodesics of $\bar{\nabla}$ induce a contact projective structure. Define a connection, $\nabla$, by $\nabla_{X}Y - \bar{\nabla}_{X}Y = \Lambda(X, Y)$, where $\Lambda_{\alpha\beta}\,^{\gamma} = -\delta_{\alpha}\,^{0}A_{\beta}\,^{\gamma}$. Because $\Lambda_{(ij)}\,^{\gamma} = 0$, $\nabla$ admits the same full set of contact geodesics as does $\bar{\nabla}$. Moreover, $\nabla$ is the representative associated to $\theta$ by Theorem A. That $\nabla\theta = 0$ and $\nabla \omega = 0$ follow by direct computation using \eqref{differenceact}. The torsions are related by $\tau_{\alpha\beta}\,^{\gamma} - \bar{\tau}_{\alpha\beta}\,^{\gamma} = 2\delta_{[\alpha}\,^{0}A_{\beta]}\,^{\gamma}$, so that $\tau_{0\beta}\,^{\gamma} = A_{\beta}\,^{\gamma} - A_{\beta}\,^{\gamma} = 0$, and $\tau_{ij}\,^{k}= 0$ (so that $\tau_{ip}\,^{p} = 0$). This shows $\nabla$ satisfies the conditions of Theorem A and proves:

\begin{proposition}\label{pshcp}
An integrable pseudo-hermitian structure induces a contact projective structure with vanishing contact torsion.
\end{proposition}

\noindent
Note that, while $\nabla_{k}J_{i}\,^{j} = 0$, $\nabla_{0}J_{i}\,^{j} = 2A_{i}\,^{p}J_{p}\,^{j}$. 
Let $\bar{R}_{\alpha\beta\gamma}\,^{\sigma}$ denote the curvature of $\bar{\nabla}$. Lemma \ref{curvaturedifference} shows that
\begin{align}
&\bar{R}_{ijk}\,^{l} - R_{ijk}\,^{l} = \omega_{ij}A_{k},^{l},& &\bar{R}_{0ij}\,^{k} - R_{0ij}\,^{k} = - \bar{\nabla}_{i}A_{j}\,^{k} = -\nabla_{i}A_{j}\,^{k}.
\end{align}
\noindent
Define a quadratic form on $T_{1, 0}$ by $Q(Z) = \levi(\bar{R}(Z, \bar{Z})Z, Z)$ for $Z \in \Gamma(T_{1, 0})$. For $Z \in \Gamma(T_{1, 0})$, Webster, \cite{Webster}, defined by $Q(Z) = K(Z)\levi(Z, Z)^{2}$ the {\bf holomorphic sectional curvature}, $K(Z)$, of the subspace spanned by $Z$ and $\bar{Z}$. The pseudo-hermitian structure has {\bf constant sectional curvature $\kappa$} if there is a constant $\kappa$ so that $K(Z) = \kappa$ for all $Z \in \Gamma(T_{1, 0})$. The pseudo-hermitian connection is said to have \textbf{transverse symmetry} if $A_{ij} = 0$. The following theorem is the desired pseudo-hermitian analogue of the Beltrami theorem.

\begin{theorem}\label{phbeltrami}
The contact projective structure induced by an integrable pseudo-hermitian structure having transverse symmetry is flat if and only if the pseudo-hermitian structure has constant sectional curvatures. 
\end{theorem}

\begin{proof}
Extend $J_{ij}$ to $J_{\alpha\beta}$ by setting $J_{\alpha 0} = 0 = J_{0\beta}$. Because $A_{ij} = 0$, the preceeding discussion shows that the pseudo-hermitian connection, $\bar{\nabla} = \nabla$, is the representative associated to $\theta$ by Theorem A of the contact projective structure induced by the given pseudo-hermitian structure. By Proposition \ref{pshcp} the induced contact projective structure has vanishing contact torsion and so, by the first Bianchi identity, \eqref{bianchi1}, $R_{[\alpha\beta\gamma]}\,^{\delta} = 0$. By the Ricci identity and $\nabla J = 0$, 
\begin{equation}\label{jrskew}
R_{\alpha\beta\gamma}\,^{p}J_{p\delta} = - R_{\alpha\beta\delta}\,^{p}J_{p\gamma},
\end{equation}
so that the tensor $R_{\alpha\beta\gamma}\,^{q}J_{q\delta}$ has the same symmetries as does the curvature of a Levi-Civita connection, and a standard argument implies that
\begin{equation}\label{polar0}
R_{\alpha\beta\gamma}\,^{q}J_{q\delta} = R_{\gamma\delta\alpha}\,^{q}J_{q\beta}.
\end{equation}
By the non-degeneracy of $J_{ij}$, $R_{0ij}\,^{q}J_{qk} = R_{jk0}\,^{q}J_{qi} = 0$ implies $R_{0ijk} = 0$. Since $\tau_{ijk} = 0$ and $W_{ijkl} = 0$, expanding \eqref{jrskew} gives
\begin{align}
\label{Rdecompose3} 0 =  &R_{ik}J_{jl} - R_{jk}J_{il} + \omega_{jk}R_{i}\,^{p}J_{pl} - \omega_{ik}R_{j}\,^{p}J_{pl} - 2\omega_{ij}R_{k}\,^{p}J_{pl} +&\\  
&\notag R_{il}J_{jk} - R_{jl}J_{ik} + \omega_{jl}R_{i}\,^{p}J_{pk} - \omega_{il}R_{j}\,^{p}J_{pk} - 2\omega_{ij}R_{l}\,^{p}J_{pk}.&
\end{align}
Set $R = R_{pq}J^{pq}$. Tracing \eqref{jrskew} in $ij$ and using \eqref{firstbianchicontracted} gives $J_{i}\,^{q}R_{q}\,^{p}J_{pk} = -R_{ik}$, so that tracing \eqref{Rdecompose3} with $J^{jl}$ gives
\begin{align}\label{einstein}
&2(n-1)R_{ij} = J_{ij}R,&\\
&\label{omegajR} 4n(n-1)R_{ijkl} = R(\omega_{jl}J_{ik} - \omega_{il}J_{jk} + \omega_{jk}J_{il} - \omega_{ik}J_{jl} - 2\omega_{ij}J_{kl}).&
\end{align}

\noindent
When $2n-1 = 3$ the hypothesis is $C_{ijk} = 0$, and when $2n-1 \geq 5$ Lemma \ref{cottonweyllemma} shows that $W_{ijkl} = 0$ implies $C_{ijk} = 0$. So $C_{(ijk)} = 0$, which, by $0 = R_{0ijk}$, \eqref{ctfcotton}, and \eqref{einstein}, implies $0 = \tfrac{1}{2n(n-1)}J_{(ij}\nabla_{k)}R$. Tracing this with $J^{ij}$ gives $\nabla_{k}R = 0$. By the traced second Bianchi identity, \eqref{secondbianchizero}, and the vanishing of $R_{0ijk}$, $\nabla_{0}R_{ij} = -\nabla_{p}R_{0ij}\,^{p} = 0$, so, by \eqref{einstein}, $0 = J_{ij}\nabla_{0}R$, and tracing this shows $\nabla_{0}R = 0$. Thus $R$ is a constant function. Using \eqref{jrskew} in \eqref{omegajR} gives $n(1-n)R_{ipkl}J_{j}\,^{p}x^{i}x^{j}x^{k}x^{l} = Rg_{ij}g_{kl}x^{i}x^{j}x^{k}x^{l}$, which implies $n(n-1)Q(Z) = R\levi(Z, Z)^{2}$ for all $Z \in \Gamma(T_{1, 0})$, 
so that the pseudo-hermitian structure has constant sectional curvatures.

Assume given the contact projective structure induced by a pseudo-hermitian structure with transverse symmetry and constant sectional curvatures. Using \eqref{polar0} and the other symmetries of $R_{ijkl}$ shows $R_{pqkl}J_{i}\,^{p}J_{j}\,^{q} = R_{ijkl}$. Polarizing the equality of quartic forms $n(1-n)R_{ipkl}J_{j}\,^{p}x^{i}x^{j}x^{k}x^{l} = Rg_{ij}g_{kl}x^{i}x^{j}x^{k}x^{l}$ given by the assumption of constant sectional curvature gives $R_{p(ijk}J_{l)}\,^{p} = cg_{(ij}g_{kl)}$, and using \eqref{polar0}, $R_{pqkl}J_{i}\,^{p}J_{j}\,^{q} = R_{ijkl}$, and relabeling gives
\begin{equation}\label{polar3}
R_{ljki} - R_{jkli} + R_{plqi}J_{j}\,^{q}J_{k}\,^{p} = c(\omega_{ji}g_{kl} + \omega_{jl}g_{ik} + \omega_{jk}g_{il}).
\end{equation}
Tracing \eqref{polar3} in $ij$ and using $R_{pqkl}J_{i}\,^{p}J_{j}\,^{q} = R_{ijkl}$ and \eqref{polar0} gives \eqref{omegajR}, for  a suitable constant, $R$. This implies \eqref{einstein} and hence $W_{ijkl} = 0$, so that if the dimension is at least $5$, the contact projective structure is flat. In three dimensions, $R_{0ijk} = 0$, so $C_{ijk} = \tfrac{1}{n}\nabla_{(i}R_{jk)}$, which vanishes by \eqref{einstein}.
\end{proof}

\begin{example}
In \cite{Matveev}, V. Matveev studies the space of Riemannian metrics projectively equivalent to a given Riemannian metric. An analogous problem is to describe the space of pseudo-hermitian structures determining the same contact projective structure. Two pseudo-hermitian structures are {\bf homothetic} if their contact one-forms differ by a constant. Two pseudo-hermitian structures represent the same CR structure if their contact one-forms differ by a non-vanishing scale factor. It is easy to prove that conformally equivalent Riemannian metrics determine the same projective structure if and only if they are homothetic. Similarly, using \eqref{lambdatransform} and the analogous formula for the transformation of the pseudo-hermitian connection under a change of scale, it may be shown that, if the dimension is at least five, two integrable pseudo-hermitian structures representing the same CR structure generate the same contact projective structure if and only if they are homothetic. The idea for the following example of non-homothetic, contact projectively equivalent pseudo-hermitian structures comes from Matveev's description, \cite{Matveev}, of Beltrami's example of projectively equivalent non-homothetic Riemannian metrics. Any linear conformal symplectic automorphism, $A$, of $(\rea^{2n}, \Omega)$, determines the contactomorphism, $x \to \psi(x) = |Ax|^{-1}Ax$, of $S^{2n-1}$. As $A$ preserves linear isotropic subspaces, $\psi$ preserves the intersections with $S^{2n-1}$ of isotropic subspaces. These intersections are equatorial spheres each of which is everywhere tangent to the contact hyperplane. In particular, $\psi$ maps contact lines on $S^{2n-1}$ to contact lines, so $\psi$ is an automorphism of the standard flat contact projective structure on $S^{2n-1}$. Equipping $\com^{n}$ with the standard Hermitian form the imaginary part of which is $\Omega$ induces on the unit sphere, $S^{2n-1}$, the standard pseudo-hermitian structure. Because $\psi$ is a contactomorphism, the standard pseudo-hermitian structure on $S^{2n-1}$ may be pulled back via $\psi$. As $\psi$ preserves contact circles, the contact geodesics of the pulled back pseudo-hermitian structure are necessarily the same as those of the standard pseudo-hermitian structure. $\psi$ is a homothety of the standard pseudo-hermitian structure if and only if it agrees with the restriction to $S^{2n-1}$ of a scalar multiple of a unitary linear transformation, and it is easily checked that this happens if and only if $A$ is a scalar multiple of a unitary transformation. Consequently, if $A$ is not a constant multiple of a unitary transformation, then these pseudo-hermitian structures are not homothetic, though they induce equivalent contact projective structures.
\end{example}

\nocite{Cap-Slovak}

\def\cprime{$'$} \def\cprime{$'$} \def\cprime{$'$} \def\cprime{$'$}
  \def\cprime{$'$}
\providecommand{\bysame}{\leavevmode\hbox to3em{\hrulefill}\thinspace}

\bibliographystyle{amsplain}

\end{document}